\documentclass[reqno,11pt]{amsart}
\usepackage{amsthm,amsfonts, amssymb, amscd}
\usepackage[all,cmtip]{xy}
\usepackage{euscript, mathrsfs}

\usepackage{multirow}

\usepackage{url}

\usepackage{tikz}
\usetikzlibrary{matrix,arrows}
\usetikzlibrary{positioning}
\usetikzlibrary{cd}

\usepackage{tikz-cd}

\usepackage[normalem]{ulem}

\usepackage{hyperref}
\hypersetup{%
    linktoc=page
}

\usepackage{cellspace} 
\let\oldtocsection=\tocsection
\let\oldtocsubsection=\tocsubsection
\renewcommand{\tocsection}[2]{\hspace{0em}\oldtocsection{#1}{#2}}
\renewcommand{\tocsubsection}[2]{\hspace{1em}\oldtocsubsection{#1}{#2}}

\usepackage[margin=1.25in]{geometry}

\newtheorem{Thm}{Theorem}[section]
\newtheorem{Lem}[Thm]{Lemma}
\newtheorem{Cor}[Thm]{Corollary}
\newtheorem{Prop}[Thm]{Proposition}

\theoremstyle{remark}
\newtheorem{Rem}[Thm]{Remark}
\theoremstyle{remark}
\newtheorem{example}[Thm]{Example}
\theoremstyle{definition}

\theoremstyle{definition}

\theoremstyle{definition}

\newtheorem{Def}[Thm]{Definition}

\numberwithin{equation}{section}

\newcommand{\A}{\mathbb{ A}}

\newcommand{\R}{\mathbb{ R}}           
\newcommand{\C}{\mathbb{C}}           
\newcommand{\Z}{\mathbb{ Z}}           

\newcommand{\QQ}{\mathbb{ Q}}        
\newcommand{\Ad}{\operatorname{Ad }}             
\newcommand{\ad}{\operatorname{ad}}             

\newcommand{\Ind}{\operatorname{Ind}}

\renewcommand{\ker}{\operatorname{ker }}
\newcommand{\diag}{\operatorname{diag}}

\newcommand{\tors}{\operatorname{tors}}

\newcommand{\SL}{\operatorname{SL}}
\newcommand{\SO}{\operatorname{SO}}
\newcommand{\Sp}{\operatorname{Sp}}
\newcommand{\Spin}{\operatorname{Spin}}

\newcommand{\fb}{{\mathfrak b}}

\newcommand{\fg}{{\mathfrak g}}

\newcommand{\fl}{{\mathfrak l}}
\newcommand{\fm}{{\mathfrak m}}

\newcommand{\fp}{{\mathfrak p}}

\newcommand{\ft}{{\mathfrak t}}
\newcommand{\fu}{{\mathfrak u}}

\newcommand{\ga}{\alpha}

\newcommand{\gD}{\Delta}

\newcommand{\gl}{\lambda}

 \newcommand{\ci}{\mathcal{I}}

 \newcommand{\cl}{\mathcal{L}}
 \newcommand{\cm}{\mathcal{M}}
 \newcommand{\cn}{\mathcal{N}}
 \newcommand{\co}{\mathcal{O}}
 \newcommand{\cp}{\mathcal{P}}
 \newcommand{\cq}{\mathcal{Q}}

\newcommand{\cz}{\mathcal{Z}}

\newcommand{\toric}{\mathcal{V}}

\renewcommand{\mod}{\;\mathrm{mod}\,}

\newcommand{\vc}{\mathbf{c}}

\newcommand{\shc}{\underline{\C}}

\newcommand{\nilp}{\nu}

\renewcommand{\tilde}{\widetilde}

\renewcommand{\bar}[1]{\overline{#1}}

\renewcommand{\hat}{\widehat}

\def\v0{\mathbf{0}}

\def\ibarT{\overline{\iota}_T}

\newcommand{\reg}{\mathrm{reg}}
\newcommand{\sfl}{\mathsf{L}}

\newcounter{myenumi}
\renewcommand{\themyenumi}{$(\arabic{myenumi})$}

\title[Extended Springer Resolution]{Lusztig sheaves, local systems, and the extended Springer resolution}

\author{William Graham}
\address{Department of Mathematics\\ University of Georgia\\ Boyd Research and Education Center\\ Athens, GA\\ 30602\\ USA
}
\email{wag@math.uga.edu}

\author{Martha Precup}
\address{Department of Mathematics\\ Washington University in St. Louis \\ One Brookings Drive\\ St. Louis, Missouri\\ 63130\\ USA  }
\email{martha.precup@wustl.edu}

\author{Amber Russell}
\address{Department of Mathematical Sciences\\ Butler University \\ 4600 Sunset Avenue \\ Indianapolis, Indiana 46208\\ USA}
\email{acrusse3@butler.edu}

\date{\today}

\begin{document}
\parskip=4pt \baselineskip=14pt

\begin{abstract} The Springer resolution of the nilpotent cone is used to give a geometric construction of the Springer correspondence by applying tools in the derived category of perverse sheaves, particularly the decomposition theorem. In this paper we continue earlier work of the authors on the extended Springer resolution, a map introduced by the first author that factors through the Springer resolution.  Our first main result shows that applying the decomposition theorem to this map yields a sum of Lusztig sheaves in all Lie types, each corresponding to a triple of so-called minimal cuspidal datum.  We also study the simple perverse sheaves that occur as summands of the resulting sheaf. Our second main result proves that every simple perverse sheaf corresponding to orbits that are Levi-regular must occur, generalizing earlier work in the type $A$ setting to all Lie types.  
\end{abstract}

\maketitle

\tableofcontents

\section{Introduction}\label{sec.intro}
Let $G$ be a complex simply connected semisimple algebraic group, and let $\cn$ be the nilpotent cone in the Lie algebra $\fg$.
The Springer resolution $\mu: \tilde{\cn} \to \cn$ plays a central role in the representation theory of $G$.
This paper studies the extended Springer resolution
\[
\psi: \tilde{\cm} \to \tilde{\cn} \to \cn,
\]
Lusztig's generalized Springer correspondence, and the relationship between them.  The paper
contains two main results, which extend previous work of the authors for $G = SL_n(\C)$ (see ~\cite{GPR23}).  The first
shows that the pushforward $\psi_* \shc_{\tilde{\cm}}[\dim \cn]$ is the direct sum of the minimal cuspidal Lusztig sheaves.
The second shows that for any irreducible $G$-equivariant local system $\cl$ on a Levi-regular nilpotent orbit $\co$, the intersection
cohomology sheaf $IC(\overline{\co}, \cl)$ appears as a summand in a minimal cuspidal Lusztig sheaf.  The second result
is motivated by the first: the Levi-regular nilpotent orbits are the only nilpotent orbits for which the extended Springer
fibers can differ from the Springer fibers.  However, the proof of the second result is independent of the first; instead,
it is based on a study of the combinatorics of Lusztig's ``symbols."

We now explain our results in more detail.  The Springer resolution has long played a prominent role in geometric representation theory, arising from 
Springer’s striking discovery that there is a close relationship (which is a bijection in type $A$) between Weyl group representations and G-orbits on $\cn$~\cite{Springer1978}.   
Borho and MacPherson showed that this relationship can be understood by considering the pushforward $\mu_* \shc_{\tilde{\cn}}[\dim \cn]$
of the constant sheaf on $\tilde{\cn}$ under the map $\mu$~\cite{Borho1981}.  
Lusztig extended Springer's results to what is called the generalized Springer correspondence, via his construction of the Lusztig sheaves on the 
nilpotent cone~\cite{Lu84}.  The Springer sheaf $\mu_* \shc_{\tilde{\cn}}[\dim \cn]$
is the particular Lusztig sheaf that realizes Springer's original correspondence. 
Unlike Springer's original construction, Lusztig's construction does not produce the Lusztig sheaves as direct summands of the pushforward of the constant sheaf on a particular variety.

In type $A$, the extended Springer resolution plays a role with respect to the generalized Springer
correspondence similar to the role of the Springer resolution in the usual Springer correspondence, in that the pushforward $\psi_* \shc_{\tilde{\cm}}[\dim \cn]$
is the direct sum of all the Lusztig sheaves.  However, in other types, not all Lusztig sheaves are obtained in this way.

We briefly describe the Lusztig sheaves.  We consider triples $\vc = (L, \co_L, \cl^L)$, where $L$ is a Levi subgroup of $G$,
$\co_L$ is a nilpotent orbit for $L$, and $\cl_L$ is an $L$-equivariant local system on $\co_L$.  A cuspidal datum is a triple $\vc$
satisfying certain conditions.  Lusztig constructs an induction functor which takes a cuspidal datum $\vc$ and produces the Lusztig sheaf
$A_{\vc}$, which is a direct sum of $G$-equivariant simple perverse sheaves on $\cn$.
In the case where $L = T$, the corresponding Lusztig sheaf is the Springer sheaf.  Each cuspidal datum, or equivalently Lusztig sheaf,
is associated to some character $\chi$ of the center $Z$ of $G$.  Conversely, given $\chi$, there
is a unique cuspidal datum $\vc_{\chi}$ with central character $\chi$ in which
the nilpotent $L$-orbit is regular.  The minimal cuspidal data are these $\vc_{\chi}$, and the 
$\A_{\vc_{\chi}}$ are the minimal cuspidal Lusztig sheaves.

The group $Z$ acts trivially on 
$\tilde{\cn}$, and the Springer sheaf $\mu_* \shc_{\tilde{\cn}}[\dim \cn]$
 is the minimal cuspidal sheaf associated to the trivial character of $Z$.  
On the other hand, $Z$ acts faithfully on $\tilde{\cm}$.
Our first main result is the following, which generalizes \cite[Theorem~5.3]{GPR23} to all types.

\begin{Thm} \label{t:main1} Let $\psi: \tilde{\cm} \to \cn$ be the extended Springer resolution.  Then 
\[
\psi_* \shc_{\tilde{\cm}}[\dim \cn] = \bigoplus_{\chi \in \hat{Z}} \A_{\vc_{\chi}}.
\]
Here the sum is taken over all characters $\chi$ of $Z$, and $\A_{\vc_{\chi}}$ is the Lusztig sheaf associated to $\vc_{\chi}$, the minimal cuspidal datum with central character $\chi$.
\end{Thm}

The proof of Theorem~\ref{t:main1} relies on several geometric facts which were verified in type $A$ in \cite{GPR23}.  
In particular, the proof relies on our description of a particular subset $\toric_{\chi}$ of the affine toric variety $\toric$ for $T$ corresponding to the cone in $\hat{T}\otimes \R$ generated by the simple roots. These computations appear in Section~\ref{sec.toric_computation} below.  The sets $\toric_{\chi}$ lead, in a natural way, to a subset $\Delta_\chi$ of simple roots associated to each character $\chi\in Z$, and thereby to a standard Levi subgroup $L_{\chi}$.  The subgroups $L_{\chi}$ are recorded in Table~\ref{fig.Levi.table} below, and by comparing with \cite[Prop.~5.1]{AHJR3} and ~\cite[Table 2.17]{Bonnafe2004},
we see that $L_{\chi}$ coincides with the Levi subgroup in the minimal cuspidal datum $\vc_{\chi}$.  
The toric variety $\toric$ plays an essential role in the definition of $\tilde{\cm}$, which inherits the $Z$-action from $\toric$. We use the sets $\toric_\chi$ to identify an open subset $\tilde{\cn}_\chi$ of $\tilde{\cn}$ on which the natural restriction maps used to define the Lusztig sheaf $\A_{\vc_{\chi}}$ are particularly well-behaved, leading to the proof of Theorem~\ref{t:main1}. 

Our second main result describes the intersection cohomology sheaves $IC(\overline{\co}, \cl)$ appearing as summands in $\psi_* \shc_{\tilde{\cm}}[\dim \cn] $.  
Each Lusztig sheaf $\A_\vc$ is a direct sum of intersection cohomology sheaves $IC(\overline{\co}, \cl)$, where 
$\co$ is a nilpotent orbit and $\cl$ is a $G$-equivariant local system on $\co$.  
Moreover, each $IC(\overline{\co}, \cl)$ appears as a summand in exactly one Lusztig sheaf.  In type $A$, all of the Lusztig sheaves are minimal
cuspidal sheaves $\A_{\vc_{\chi}}$ as in Theorem~\ref{t:main1}.  Therefore,  $\psi_* \shc_{\tilde{\cm}}[\dim \cn] $ contains every $IC(\overline{\co}, \cl)$ as a summand, 
with multiplicity as determined by the generalized Springer correspondence.  In other types, the situation is more complex:  there are Lusztig sheaves which are not minimal cuspidal, so
not every $IC(\overline{\co}, \cl)$ is contained in $\psi_* \shc_{\tilde{\cm}}[\dim \cn] $.  

The Levi-regular orbits are the nilpotent orbits $\co$ such that $\co \cap \fl$ is the regular nilpotent $L$-orbit in some Levi subalgebra $\fl$.
These orbits can be identified using work of Russell (see \cite{Russell}),
in the following way.  Let $\nu$ be a regular nilpotent element of $\fg$, and choose an $\mathfrak{sl}_2$-triple $\{h, \nu, \nu^-\}$ in $\fg$ containing $\nu$ as nil-positive element.  Write $\fg = \bigoplus \fg_i$ for the decomposition of $\fg$ into eigenspaces for $\ad h$.  The Levi-regular orbits are those which intersect $\fg_2$.  In type $A$, every nilpotent orbit intersects $\fg_2$;
for general types, the partition classification of orbits intersecting 
$\fg_2$ is described in~\cite{Russell}.  Our second main theorem is the following, which directly generalizes the type $A$ case.

\begin{Thm} \label{t:main2}
If $\co$ is Levi-regular, then every intersection cohomology sheaf $IC(\overline{\co}, \cl)$ for $\cl$ any $G$-equivariant local system on $\co$ is obtained from the extended Springer resolution as a summand of  $\psi_* \shc_{\tilde{\cm}}[\dim \cn] $, with multiplicity as determined by the generalized Springer correspondence.
\end{Thm}

In classical types outside of type $D$, a converse statement holds: If $\co$ is not Levi-regular, then there exist intersection cohomology sheaves $IC(\overline{\co}, \cl)$ which
are not obtained from the extended Springer resolution.  In type $D$, the analogous statement holds if $\co$ is not  ``almost Levi-regular"  (see Theorem ~\ref{thm.support}).

These results provide new insight into the structure of Lusztig sheaves for the ``minimal'' Levi subgroups $L_\chi$ associated to each central character $\chi$. Indeed, these Levi subgroups are all of type $A$. Theorem~\ref{t:main2} shows that the minimal cuspidal Lusztig sheaves completely capture the data of the intersection cohomology sheaves supported on Levi-regular orbits. As far as we are aware, this property of the generalized Springer correspondence has not previously appeared in the literature.

The contents of the paper are as follows. Section~\ref{sec.notation} recalls necessary background information and definitions, including the definition of the extended Springer resolution. Our key technical computations appear in Section~\ref{sec.toric_computation}, where we identify the sets $\toric_\chi$ for each central character $\chi$. Section~\ref{s.Lusztig} completes the proof of Theorem~\ref{t:main1} by
building on the proof techniques in type $A$ from~\cite{GPR23}.  

The proof of Theorem~\ref{t:main2} appears in Section~\ref{sec.g2.intersection}. This proof requires the use of Lusztig's construction of symbols in the classical types to determine the generalized Springer correspondence~\cite{Lu84} and also the work of Lusztig--Spaltenstein~\cite{LuSpal85} for the correspondence in the case of the spin groups. The arguments of Section~\ref{sec.toric_computation} require two technical lemmas. While neither is difficult, we include short proofs in Appendix~\ref{appendix.ZintypeD} and~\ref{appendix.Bonnafe} for the interested reader. We also give an overview of the generalized Springer correspondence for the classical groups in Appendix~\ref{s.overview} to aid understanding of the techniques applied in Section~\ref{sec.g2.intersection}.

\noindent\textbf{Acknowledgments.} The authors are grateful to Anthony Henderson, who suggested the connection between the extended Springer resolution and minimal cuspidal data that appears in Theorem~\ref{t:main1}.  We would like to thank David Kazhdan and Antonio Nigro for informing us of the error in the original Theorem~\ref{thm.support} (2).
The second author is partially supported by NSF CAREER grant DMS-2237057. This material is based upon work supported by a grant from the Institute for Advanced Study School of Mathematics and by the National Science Foundation under Grants No.~DMS-2302116 and DMS-2302117.

\medskip

\noindent\textbf{AI Disclosure.} The first version of this paper was written without AI use (beyond typographical corrections), but
it contained an error in the statement of Theorem~\ref{thm.support} (2), which applied to the exceptional groups.
This error was discovered by ChatGPT and confirmed by Claude.  We have removed this statement and given a different
argument for the proof of Theorem \ref{t:main2} for the exceptional groups using results of Lusztig.  
Claude also pointed out some smaller errors after a review of the first version, which do not affect the statements of the results.
We have made the necessary minor corrections.


\section{Background and notation}\label{sec.notation}
In this paper, $G$ will denote a semisimple, simply connected complex algebraic group.  Let $B \supset T$ be a Borel subgroup and maximal torus of $G$, respectively, and let $U$ denote the unipotent radical of $B$, so $B = TU$.  We denote the Lie algebra of an algebraic group by the corresponding fraktur letter, so the Lie algebras of $G$, $B$, $T$ and $U$ are, respectively, $\fg, \fb, \ft$, and $\fu$.  The Weyl group of
$G$ with respect to $T$ is denoted $W$.  

Let $Z$ denote the center of $G$, so $G_{ad} = G/Z$ is the adjoint group corresponding to $G$, with maximal torus $T_{ad} = T/Z$.  Let $\cp = \hat{T}$ denote the group of characters of $T$.  Let $\Phi \subset \hat{T}$ denote the root system of $\fg$ with respect to $T$, and let $\cq = \langle \Phi \rangle$ denote the root lattice.   We take the positive roots $\Phi^+$ to be the set of roots $\ga\in \Phi$ such that the corresponding root space $\fg_{\ga}$ is contained in $\fb$.  For each root $\ga$, let $X_\ga \in \fg_\ga \setminus \{0\}$ be a fixed root vector.   A Levi subgroup $L$ is a Levi factor (i.e.~maximal reductive subgroup) of a parabolic subgroup of $G$.
Any Levi subgroup contains a subgroup conjugate to $T$.   Unless otherwise stated, we will assume that $L$ is standard, i.e., $L$ contains $T$.   By the Levi type of a Levi subgroup, we
mean the type of the semisimple subalgebra $[\fl, \fl]$.  If $L$ is standard, this corresponds to a subdiagram of the Dynkin diagram of $G$ whose nodes we color black.

Let $\ga_1, \ldots, \ga_n\in \Phi^+$ denote the simple roots and $\gl_1, \ldots, \gl_n$ the fundamental dominant weights. Throughout this paper, we use the notational conventions of~\cite{HumphreysLieAlg} for the labelling of the simple roots and of the fundamental dominant weights. The \emph{minuscule weights} are the nonzero dominant weights which are minimal with respect to
the partial ordering $\mu \leq \lambda$ if $\lambda - \mu$ is a nonnegative sum of simple roots.  The minuscule weights are a subset of the fundamental dominant weights, with equality in type $A$, and they form a complete set of coset representatives for $\cp/\cq$.   A list of the minuscule weights for each irreducible root system can be found in  \cite[\S 13.4, Exercise 13]{HumphreysLieAlg}.

Suppose $\gl_k$ is minuscule.   We can write $\gl_k = \sum_j c_{kj} \ga_j$, where the $c_{kj}$ are in $\QQ$.  Set $d_{kj} = r(c_{kj})$, where $r(c) := c - \lfloor c \rfloor$ denotes the fractional part of $c$ (e.g.~ $r(\frac{7}{3}) = \frac{1}{3}$), and set $\mu_k := \sum_j d_{kj} \ga_j$.  We refer to $\mu_k$ as the \emph{minimal element corresponding to $\gl_k$}; it is the element of $\cp/\cq$ in the $\cq$-coset of $\gl_k$ whose coefficients in the $\ga_i$-basis are as small as possible while being nonnegative. These are discussed in type $A$ in \cite{GPR23}, and described for other types in Table~\ref{fig.minimal.table} below.

The group of characters of $Z$ is denoted by $\hat{Z}$ (with analogous notation for other finite groups). The restriction map $\cp \to \hat{Z}$ is surjective with kernel $\cq$, so there is an isomorphism $\cp/\cq \to \hat{Z}$,
which we denote by $\mu \mapsto \overline{\mu}$.  To simplify notation, we will sometimes write $\mu$ instead of $\overline{\mu}$ if it is understood that
we are working in $\cp/\cq$.  We will usually
write the group operation on $\cp/\cq$ as addition, and the group operations on $Z$ and $\hat{Z}$ as multiplication.
The group $\cp/\cq \cong \hat{Z}$  consists of $0$, together with the $\overline{\lambda}_i = \overline{\mu}_i$, where the $\lambda_i$ are the minuscule dominant weights.     By the above remarks, any character $\chi \in \hat{Z}$ is represented by some minuscule weight $\gl_i$, or equivalently by the corresponding $\mu_i$; that is, $\chi = \overline{\lambda}_i = \overline{\mu}_i$ for some~$i$.

\subsection{Characters and component groups}  We write $H_0$ for the identity component of an algebraic group $H$, so $H/H_0$ is the component group of $H$. This section contains some results about the component groups of centralizers of regular nilpotent elements, and the characters of these groups.  

Suppose that $L$ is a (standard) Levi subgroup of $G$.  Let $\Phi_L$ denote the root system of $\fl$ with respect to $T$.   Let $Z(L)$ denote the center of $L$, and define $\cz(L) := Z(L)/Z(L)_0$. Since the center $Z = Z(G)$ of $G$ is finite,
we have $\cz(G) = Z$.   Recall that an element $\nilp \in \fg$ is \emph{regular} if its centralizer $G^{\nilp}$ has dimension as small as possible.  If $\nilp$ is a regular nilpotent element of $\fg$ and $\nilp_L$ is a regular nilpotent 
element of $\fl$, then the corresponding component groups are $L^{\nilp_L} / L^{\nilp_L}_0 =\cz(L)$ and $G^{\nilp}/G^{\nilp}_0 = Z$ (see~\cite[Prop.~3.1]{Bonnafe2004}).

The inclusion $Z \subset Z(L)$ induces a surjective map $h_L: Z \to \cz(L)$ \cite{Bonnafe2004}.  Pullback via $h_L$ defines an injective map
$h_L^*: \widehat{\cz(L)} \to \hat{Z}$ of character groups.  This map can be described as follows under our assumption that $L$ contains $T$.
By  \cite[Prop.~2.1]{Bonnafe2004}, we have
\[ 
\hat{Z} \cong (\cp/\cq )_{\tors} = \cp/\cq \quad \text{ and } \quad \widehat{\cz(L)} \cong (\cp /\langle \Phi_L \rangle)_{\tors},
\]
where the subscript ``tors" indicates the torsion subgroup.
The inclusion $\Phi_L \subset \Phi_G$
induces the map
$$
h_L^*: \widehat{\cz(L)} \cong (\cp /\langle \Phi_L \rangle)_{\tors} \to \widehat{\cz(G)} = \hat{Z} \cong \cp/\cq .
$$

\begin{Lem} \label{l:image}
Let $\chi \in \cp/\cq$ be represented by $\mu = \sum a_i \alpha_i$, where $\alpha_1, \ldots, \alpha_n$ are simple roots and $a_i \in \QQ$.  Then $\chi$ is in the image of $h_L^*$ if and only if for each $\alpha_i$ which is not in $\Phi_L$, we have $a_i \in \Z$.
\end{Lem}

\begin{proof}
The element $\chi$ is in the image of $h_L^*$ if and only if $\mu - \mu'\in \cq$ for some $\mu' = \sum b_i \alpha_i \in \cp$ (where $b_i \in \QQ$)  which represents a torsion element of $\cp/\langle \Phi_L \rangle$.
The condition that $\mu'$ represents a torsion element is equivalent to the condition that $b_i = 0$ whenever $\alpha_i \not\in \Phi_L$.  

If $a_i \in \Z $ whenever $\alpha_i \not\in \Phi_L$, define $\mu' =  \sum b_i \alpha_i$ by 
\[
b_i = \begin{cases}
a_i & \mbox{ if } \alpha_i \in \Phi_L \\
0 & \mbox{ if }  \alpha_i \not \in \Phi_L.
\end{cases}
\]
By construction, $\mu-\mu' = \sum (a_i-b_i) \alpha_i = \sum_{\alpha_i\notin \Phi_L} a_i \alpha_i \in \cq$. 
Thus, if $\chi_L \in \widehat{\cz}(L)$ is the character represented by $\mu'$, then $h_L^*(\chi_L) = \chi$.  Conversely, if $\chi = h_L^*(\chi_L)$ for some $\chi_L$, let $\mu' = \sum b_i \alpha_i \in \cp$ represent $\chi_L$.
Then $\mu - \mu' \in \cq$.  Since $b_i = 0$ whenever $\alpha_i \not\in \Phi_L$, this implies that the corresponding $a_i$ are
integers, as desired.
\end{proof}


\subsection{The extended Springer resolution}
Let $\cn$ denote the nilpotent cone in $\fg$. The Springer resolution is the map 
$$
\mu: \tilde{\cn} = G \times^B \fu \to \cn
$$
given by $[g,x] \mapsto g \cdot x$, where $g \cdot x = \Ad(g) (x)$.  
We briefly recall the construction of the extended Springer resolution; see \cite{Graham2022} or \cite{GPR23} for more details.  Let $\{h, \nilp, \nilp^- \}$ be an $\mathfrak{sl}_2$-triple such that $h \in \ft$, and $\nilp$ is the regular nilpotent element $\nilp = \sum_{i=1}^n X_{\ga_i}$.  Let $\fg_i$ denote the $i$-eigenspace for $\ad h$ on $\fg$; then $\fg = \oplus \fg_i$.  Define $\toric_{ad} = \mathfrak{g}_2$; this space is the span of the root spaces $\fg_{\ga}$ where $\ga$ is simple.   There is a $B$-equivariant projection $p:\fu \to \toric_{ad}$ obtained from the identification of $\toric_{ad}$ with $\fu/[\fu, \fu]$. 
 
The space $\toric_{ad}$ is the affine toric variety for $T_{ad}$ corresponding to the cone in $\hat{T} \otimes \R$ generated by the simple roots.  The lattice associated with this toric variety is the root lattice $\cq$.  
Define $\toric$ to be the affine toric variety for $T$ corresponding to the same cone, but with lattice equal to $\cp$.  
Since $\cq\subset \cp$ there is a natural projection $\pi:\toric \to \toric_{ad}$.   
Define $\tilde{\fu}$ as the fiber product $\toric\times_{\toric_{ad}} \fu$ corresponding to the following commutative diagram:
\[
\xymatrix{
 \tilde{\fu}:= \toric\times_{\toric_{ad}} \fu \ar[r]\ar[d]& \fu \ar[d]_p \\ 
 \toric \ar[r]^\pi & \toric_{ad} }
\]
Projection from $\tilde{\fu}$ to the second factor yields a $B$-equivariant map $\eta: \tilde{\fu} \to \fu$.  This induces a map $\tilde{\eta}: \tilde{\cm} = G \times^B \tilde{\fu} \to \tilde{\cn} = G \times^B \fu$.
The action of $Z$ on $\toric$ induces an action on $\tilde{\cm}$ and the map $\tilde{\eta}$
is the quotient by $Z$.  By definition, the extended Springer resolution $\psi$ is the composition of $\tilde{\eta}$ with the Springer resolution $\mu: \tilde{\cn}$: that is,
\begin{equation} \label{e:extended}
\psi: \tilde{\cm} \stackrel{\tilde{\eta}}{\longrightarrow} \tilde{\cn} \stackrel{\mu}{\longrightarrow} \cn.
\end{equation}

The coordinate rings $\C[\toric_{ad}]$ and $\C[\toric]$ of the affine toric varieties $\toric_{ad}$ and $\toric$ have the following description (see \cite[\S 3.1]{GPR23} for more details). Given a weight $\lambda$ of $T$, we let $e^{\gl}$ denote the corresponding function on $T$. Recall that $\ga_1, \ldots, \ga_n$ are the simple roots.  Let $\gl_{i_1}, \ldots, \gl_{i_k}$ denote the minuscule weights, with corresponding minimal elements $\mu_{i_1}, \ldots, \mu_{i_k}$. We write $x_i = e^{\ga_i}$, and $v_j = e^{\mu_j}$.  Then $\C[\toric_{ad}] = \C[x_1, \ldots, x_n]$ is a polynomial ring and  $\C[\toric] =  \C[x_1, \ldots, x_n, v_{i_1}, \ldots, v_{i_k}]/I$ for some ideal $I$.  The action of $z \in Z$ on $\C[\toric]$ is given by
\begin{equation} \label{e:zaction}
z \cdot x_i = x_i, \hspace{.2in} z \cdot v_j = e^{\mu_j}(z^{-1}) v_j,
\end{equation}
where the first equation holds because $z$ acts trivially on $\C[\toric_{ad}] $.

In type $A_n$, each of the fundamental dominant weights is minuscule, but in other simple types there are at most $3$ minuscule weights.  An explanation
of how to compute the minimal weights for type $A_n$ can be found in \cite{GPR23}.  Table~\ref{fig.minimal.table} below lists the minuscule and corresponding minimal weights for each type besides $A_n$ (using the numbering of \cite[\S 13.4, Exercise 13]{HumphreysLieAlg}, and assuming that for types $B_n$ and $C_n$, $n \geq 2$, and for $D_n$, $n \geq 4$).  For readability, we use the following notation in the table.  If $c \in \QQ$, recall that $r(c) = c - \lfloor c \rfloor$ denotes the fractional part of $c$.  Let $d = n-2$ if $n$ is odd, and $d = n-3$ if $n$ is even, and write
\begin{align*}
A & =   \alpha_1 + 2\alpha_2 + \cdots + (n-2) \alpha_{n-2} \\
B & = \alpha_1 + \alpha_3 + \cdots + \alpha_{d}.
\end{align*}

\begin{table}[h] 
\small{
\begin{center}
\resizebox{\linewidth}{!}{
\begin{tabular}{|Sc||Sc|Sc|} \hline
Type of $G$ & Minuscule weight & Minimal weight\\
\hline\hline
$B_n$ & $\lambda_n =  \frac{1}{2}\left( \alpha_1+2\alpha_2 + \cdots + n\alpha_n \right) $  & $\mu_n  =  \frac{1}{2} (\alpha_1 + \alpha_3 + \cdots )$   \\
\hline
$C_n$ & $\lambda_1 =  \alpha_1 + \cdots + \alpha_{n-1}+\frac{1}{2}\alpha_n $  & $\mu_1  =  \frac{1}{2} \alpha_n$   \\
\hline
$D_n$ & $\lambda_1 =  \alpha_1 + \cdots + \alpha_{n-2}+\frac{1}{2} (\alpha_{n-1} + \alpha_n) $  & $\mu_1  = \frac{1}{2} (\alpha_{n-1} + \alpha_n) $   \\
            & $\lambda_{n-1} =  \frac{1}{2} A +  \frac{n}{4}\alpha_{n-1}+\frac{n-2}{4}\alpha_n $ 
            & $\mu_{n-1} =  \frac{1}{2} B+ r\left( \frac{n}{4} \right ) \alpha_{n-1}+ r\left( \frac{n-2}{4} \right )  \alpha_n$ \\
            & $\lambda_{n} =   \frac{1}{2} A+ \frac{n-2}{4}\alpha_{n-1}+\frac{n}{4}\alpha_n$ 
            & $\mu_{n} =  \frac{1}{2} B + r\left( \frac{n-2}{4} \right ) \alpha_{n-1}+ r\left( \frac{n}{4} \right )  \alpha_n$ \\
\hline
$E_6$ &  $\lambda_{1} =  \frac{1}{3} \left(4  \alpha_1 + 3 \alpha_2 + 5 \alpha_3  + 6 \alpha_4 + 4 \alpha_5 + 2 \alpha_6  \right)$
            &  $\mu_{1} =  \frac{1}{3} \left(  \alpha_1 + 2 \alpha_3   +  \alpha_5 + 2 \alpha_6  \right)     $ \\   
            & $\lambda_{6} =   \frac{1}{3} \left(2  \alpha_1 + 3 \alpha_2 + 4 \alpha_3  + 6 \alpha_4 + 5 \alpha_5 + 4 \alpha_6  \right)$ 
            & $\mu_{6} =   \frac{1}{3} \left(2  \alpha_1 +  \alpha_3  + 2 \alpha_5 + \alpha_6  \right)$ \\
\hline
 $E_7$ &  $\lambda_7 = \frac{1}{2} \left(2  \alpha_1 + 3 \alpha_2 + 4 \alpha_3  + 6 \alpha_4 + 5 \alpha_5 + 4 \alpha_6  + 3 \alpha_7  \right)$ 
            & $\mu_7 = \frac{1}{2} \left( \alpha_2 +  \alpha_5  + \alpha_7  \right)$ \\ \hline
\end{tabular}
}
\end{center}
}
\medskip
\caption{Minuscule and corresponding minimal weights} \label{fig.minimal.table}
\end{table}

Associated to each nontrivial character $\chi$ in $ \hat{Z}\simeq \cp/\cq$, there is a standard Levi subgroup $L_{\chi}$, defined as follows.  Suppose $\chi$ is represented by the minuscule weight $\lambda_k$
or equivalently, by the minimal element $\mu_k$. Let $\Delta_{\chi}$ denote the set of simple roots which appear with nonzero coefficient in $\mu_k$, and let $L = L_{\chi}$ be the standard Levi subgroup determined by the subset $\Delta_{\chi}$.   

The last column of Table~\ref{fig.Levi.table} lists the Dynkin diagrams corresponding to the pairs $(G, L_{\chi})$, where the simple roots in $\Delta_{\chi}$ correspond to the vertices colored in black except in type $A$, where our convention is that the boxed parts of the displayed Dynkin diagram are colored in black.  
The notation for the type $D_n$ with $n$ even rows of Table~\ref{fig.Levi.table} is as follows.  In this case, the center $Z$ is isomorphic to $\Z_2 \times \Z_2$.  For any nontrivial character $\chi \in \hat{Z}$, the group $\ker \chi$ is generated by some $z \in Z$ (necessarily of order 2).  Moreover, $\chi$ is uniquely determined by $\ker \chi = \langle z \rangle$. We write $Z = \{ 1, z_1, z_{n-1}, z_n \}$, where the elements $z_i$ are defined in Section~\ref{ss:Dn} below.  In the entries of the table, the character $\chi$ of $Z$ is specified by its kernel. Note that Table~\ref{fig.Levi.table} is very similar to \cite[Table 2.17]{Bonnafe2004}, and $z_i$ is
the element denoted in \cite{Bonnafe2004} as $\ibarT(\varpi_i^{\vee})$.   We discuss the relationship between the tables further in Lemma~\ref{lem:Bonnafe} and Remark \ref{r:bonnafe} below, with more detail presented in Appendix~\ref{appendix.Bonnafe}.

\begin{table}[h] 
\small{
\begin{center}
\resizebox{\linewidth}{!}{
\begin{tabular}{|Sc||Sc|Sc|Sc|Sc|}
\hline
Type of $G$ & $Z$ & $\chi\in \widehat{Z}, \,\chi \neq 1$ &  Levi Type & Levi Diagram \\ \hline \hline
\multirow{2}*{$A_{n-1}$} & \multirow{2}*{$\Z_n$}  & $\chi^{d} = 1,$  &  \multirow{2}*{$A_{d-1} \times \cdots \times A_{d-1}$}  &  
\multirow{2}*{
\begin{tikzpicture}[scale=1]
\node[rectangle, draw, inner sep=2pt] (1) at (1,0) {$A_{d-1}$};
\node[circle, draw, inner sep=2pt] (2) at (2,0) {};
\node[rectangle, draw, inner sep=2pt] (3) at (3,0) {$A_{d-1}$};
\node[circle, draw, inner sep=2pt] (4) at (4,0) {};
\node[circle, draw, inner sep=2pt] (5) at (5,0) {};
\node[rectangle, draw, inner sep=2pt] (6) at (6,0) {$A_{d-1}$};
\draw (1)--(2);
\draw (2)--(3);
\draw (3)--(4);
\draw[dashed] (4)--(5);
\draw (5)--(6);
\end{tikzpicture} }  \\ 
& & $d$ a divisor of $n$ &&\\ \hline \hline
$C_n$ & $\Z_2$ & $\chi^2= 1$ & $A_1$ \text{(long root)} & \begin{tikzpicture}[scale=1]
\node[circle, draw, inner sep=2pt] (1) at (1,0) {};
\node[circle, draw, inner sep=2pt] (2) at (2,0) {};
\node[circle, draw, inner sep=2pt] (3) at (3,0) {};
\node[circle, draw, inner sep=2pt] (4) at (4,0) {};
\node[circle, draw, inner sep=2pt] (5) at (5,0) {};
\node[inner sep=-0.5pt, outer sep=-0.5pt] (7) at (5.4,0) {};
\node[circle, fill=black, draw, inner sep=2pt] (6) at (6,0) {};
\draw (1)--(2);1
\draw (2)--(3);
\draw (3)--(4);
\draw[dashed] (4)--(5);
\draw[double, double distance=1.5pt] (5)--(7);
\draw[double, double distance=1.5pt, <-] (7)--(6);
\end{tikzpicture}  \\  \hline \hline
$B_n$ & \multirow{2}*{$\Z_2$} & \multirow{2}*{$\chi^2=1$} &\multirow{2}*{$A_1\times A_1\times \cdots \times A_1$}  &
\multirow{2}*{\begin{tikzpicture}[scale=1]
\node[circle, fill=black, draw, inner sep=2pt] (1) at (1,0) {};
\node[circle, draw, inner sep=2pt] (2) at (2,0) {};
\node[circle, fill=black, draw, inner sep=2pt] (3) at (3,0) {};
\node[circle, draw, inner sep=2pt] (4) at (4,0) {};
\node[circle, fill = black, draw, inner sep=2pt] (5) at (5,0) {};
\node[inner sep=-0.5pt, outer sep=-0.5pt] (7) at (5.6,0) {};
\node[circle, draw, inner sep=2pt] (6) at (6,0) {};
\draw (1)--(2);
\draw (2)--(3);
\draw (3)--(4);
\draw[dashed] (4)--(5);
\draw[double, double distance=1.5pt, ->] (5)--(7);
\draw[double, double distance=1.5pt] (7)--(6);
\end{tikzpicture} }  \\
$n$ even & & & & \\  \hline 
$B_n$ & \multirow{2}*{$\Z_2$} & \multirow{2}*{$\chi^2=1$} & \multirow{2}*{$A_1\times A_1\times \cdots \times A_1$} & \multirow{2}*{\begin{tikzpicture}[scale=1]
\node[circle, fill=black, draw, inner sep=2pt] (1) at (1,0) {};
\node[circle, draw, inner sep=2pt] (2) at (2,0) {};
\node[circle, fill=black, draw, inner sep=2pt] (3) at (3,0) {};
\node[circle, draw, inner sep=2pt] (4) at (4,0) {};
\node[circle, draw, inner sep=2pt] (5) at (5,0) {};
\node[inner sep=-0.5pt, outer sep=-0.5pt] (7) at (5.6,0) {};
\node[circle, fill=black, draw, inner sep=2pt] (6) at (6,0) {};
\draw (1)--(2);
\draw (2)--(3);
\draw (3)--(4);
\draw[dashed] (4)--(5);
\draw[double, double distance=1.5pt, ->] (5)--(7);
\draw[double, double distance=1.5pt] (7)--(6);
\end{tikzpicture}}  \\
$n$ odd & & & &   \\  \hline \hline
$D_n$ & $\Z_2\times \Z_2$ & $\ker \chi = \langle z_n \rangle$ & $A_1\times A_1\times \cdots \times A_1$ &
\begin{tikzpicture}[scale=1]
\node[circle, fill=black, draw, inner sep=2pt] (1) at (1,0) {};
\node[circle, draw, inner sep=2pt] (2) at (2,0) {};
\node[circle, fill=black, draw, inner sep=2pt] (3) at (3,0) {};
\node[circle, draw, inner sep=2pt] (4) at (4,0) {};
\node[circle, fill= black, draw, inner sep=2pt] (5) at (5,0) {};
\node[circle, draw, inner sep=2pt] (6) at (6,0) {};
\node[circle, fill=black, draw, inner sep=2pt] (7) at (7,0.3) {};
\node[circle, draw, inner sep=2pt] (8) at (7,-0.3) {};
\draw (1)--(2);
\draw (2)--(3);
\draw[dashed] (3)--(4);
\draw (4)--(5);
\draw (5)--(6);
\draw (6)--(7);
\draw (6)--(8);
\end{tikzpicture}
 \\
$n$ even & & & &   \\ 
\empty & \empty & $\ker \chi = \langle z_{n-1} \rangle$ & $A_1\times A_1\times \cdots \times A_1$ & 
\begin{tikzpicture}[scale=1]
\node[circle, fill=black, draw, inner sep=2pt] (1) at (1,0) {};
\node[circle, draw, inner sep=2pt] (2) at (2,0) {};
\node[circle, fill=black, draw, inner sep=2pt] (3) at (3,0) {};
\node[circle, draw, inner sep=2pt] (4) at (4,0) {};
\node[circle, fill= black, draw, inner sep=2pt] (5) at (5,0) {};
\node[circle, draw, inner sep=2pt] (6) at (6,0) {};
\node[circle, draw, inner sep=2pt] (7) at (7,0.3) {};
\node[circle, fill=black, draw, inner sep=2pt] (8) at (7,-0.3) {};
\draw (1)--(2);
\draw (2)--(3);
\draw[dashed] (3)--(4);
\draw (4)--(5);
\draw (5)--(6);
\draw (6)--(7);
\draw (6)--(8);
\end{tikzpicture}\\ 
 & & & &   \\ 
\empty & \empty & $\ker \chi = \langle z_1 \rangle$ & $A_1\times A_1$ &  
\begin{tikzpicture}[scale=1]
\node[circle, draw, inner sep=2pt] (1) at (1,0) {};
\node[circle, draw, inner sep=2pt] (2) at (2,0) {};
\node[circle, draw, inner sep=2pt] (3) at (3,0) {};
\node[circle, draw, inner sep=2pt] (4) at (4,0) {};
\node[circle, fill=black, draw, inner sep=2pt] (5) at (5,0.3) {};
\node[circle, fill=black, draw, inner sep=2pt] (6) at (5,-0.3) {};
\draw (1)--(2);
\draw[dashed] (2)--(3);
\draw (3)--(4);
\draw (4)--(5);
\draw (4)--(6);
\end{tikzpicture}\\ 
 & & & &   \\ 
\hline 
$D_n$ & $\Z_4$ & $\chi^2 = 1$ & $A_1\times A_1$ &
\begin{tikzpicture}[scale=1]
\node[circle, draw, inner sep=2pt] (1) at (1,0) {};
\node[circle, draw, inner sep=2pt] (2) at (2,0) {};
\node[circle, draw, inner sep=2pt] (3) at (3,0) {};
\node[circle, draw, inner sep=2pt] (4) at (4,0) {};
\node[circle, fill=black, draw, inner sep=2pt] (5) at (5,0.3) {};
\node[circle, fill=black, draw, inner sep=2pt] (6) at (5,-0.3) {};
\draw (1)--(2);
\draw[dashed] (2)--(3);
\draw (3)--(4);
\draw (4)--(5);
\draw (4)--(6);
\end{tikzpicture} \\
$n$ odd & & & &  \\ 
\empty & \empty &  & & \\ 
 \empty & \empty & $\chi^4 = 1, \, \chi^2\neq 1$ & $A_1\times \cdots \times A_1 \times A_3$ &
 \begin{tikzpicture}[scale=1]
\node[circle, fill=black, draw, inner sep=2pt] (1) at (1,0) {};
\node[circle, draw, inner sep=2pt] (2) at (2,0) {};
\node[circle, fill=black, draw, inner sep=2pt] (3) at (3,0) {};
\node[circle, fill=black, draw, inner sep=2pt] (4) at (4,0) {};
\node[circle, draw, inner sep=2pt] (5) at (5,0) {};
\node[circle, fill=black, draw, inner sep=2pt] (6) at (6,0) {};
\node[circle, fill=black, draw, inner sep=2pt] (7) at (7,0.3) {};
\node[circle, fill=black, draw, inner sep=2pt] (8) at (7,-0.3) {};
\draw (1)--(2);
\draw (2)--(3);
\draw[dashed] (3)--(4);
\draw (4)--(5);
\draw (5)--(6);
\draw (6)--(7);
\draw (6)--(8);
\end{tikzpicture} \\
 & & & &   \\ 
\hline \hline
$E_6$ & $Z_3$ & $\chi^3 = 1$ & $A_2 \times A_2$ & 
\begin{tikzpicture}[scale=1]
\node[circle, fill=black, draw, inner sep=2pt] (1) at (1,0) {};
\node[circle, fill=black, draw, inner sep=2pt] (2) at (2,0) {};
\node[circle, draw, inner sep=2pt] (3) at (3,0) {};
\node[circle, fill=black, draw, inner sep=2pt] (4) at (4,0) {};
\node[circle, fill=black, draw, inner sep=2pt] (5) at (5,0) {};
\node[circle, draw, inner sep=2pt] (6) at (3,1) {};
\draw (1)--(2);
\draw (2)--(3);
\draw (3)--(4);
\draw (4)--(5);
\draw (3)--(6);
\end{tikzpicture}  \\ \hline \hline
$E_7$ & $\Z_2$ & $\chi^2 =1$ & $A_1\times A_1\times A_1$ & 
\begin{tikzpicture}[scale=1]
\node[circle,  draw, inner sep=2pt] (1) at (1,0) {};
\node[circle, draw, inner sep=2pt] (2) at (2,0) {};
\node[circle, draw, inner sep=2pt] (3) at (3,0) {};
\node[circle,  fill=black, draw, inner sep=2pt] (4) at (4,0) {};
\node[circle,  draw, inner sep=2pt] (5) at (5,0) {};
\node[circle, fill=black, draw, inner sep=2pt] (6) at (6,0) {};
\node[circle, fill=black, draw, inner sep=2pt] (7) at (3,1) {};
\draw (1)--(2);
\draw (2)--(3);
\draw (3)--(4);
\draw (4)--(5);
\draw (5)--(6);
\draw (3)--(7);
\end{tikzpicture}
\\ \hline
\end{tabular}
}
\end{center}
\medskip
}
\caption{The Levi subgroups $L_\chi$ for nontrivial $\chi\in \hat{Z}$.} \label{fig.Levi.table}
\end{table}


\section{The subsets $\toric_\chi$}  \label{sec.toric_computation}
The center $Z$ of $G$ acts on the affine toric variety $\toric$. Given any character $\chi \in \widehat{Z}$, we recall from \cite{GPR23} the definition
\[
\toric_\chi := \{ x\in \toric\mid  Z^x \subseteq \ker \chi  \}. 
\]
Let $\toric_{ad, \chi}$ denote the image of $\toric{\chi}$ under the surjective map $\pi:\toric \to \toric_{ad}$. The subsets $\toric_\chi$ play a key role in understanding how the extended Springer resolution relates to the Lusztig sheaves.  These sets were described in \cite{GPR23} in type $A$.  In this section we prove the following proposition, which describes $\toric_{\chi}$ in general.  

\begin{Prop} \label{p:Levi-chi}
Let $\chi$ be a nontrivial character in $ \hat{Z}$.  Suppose $\chi$ is represented by the minuscule weight $\lambda_k$ (or equivalently, by the minimal element $\mu_k$). Let $\Delta_{\chi}$ denote the set of simple roots which appear with nonzero coefficient in the minimal weight $\mu_k$ and $L = L_{\chi}$ be the standard Levi subgroup whose simple roots are $\Delta_{\chi}$.  Then:
\begin{enumerate}
\item $\chi = h_L^*(\chi_L)$ for a (necessarily unique) $\chi_L \in \widehat{\cz(L)}$.
\item $\toric_{ad, \chi}$ is defined by the conditions $x_i \neq 0$ for $\alpha_i \in \Delta_{\chi}$.
\end{enumerate}
\end{Prop}

Part (1) of the proposition is an immediate consequence of Lemma \ref{l:image}, since $\chi$ is represented by $\mu_k$, and the simple root $\alpha_i$ appears in the formula for $\mu_k$ if and only if $\alpha_i \in \Phi_L$.  In type $A_n$, the simple roots appearing in $\mu_k$ were computed and part (2) was proved in \cite[\S 3.2]{GPR23}. For types besides type $A_n$, the simple roots appearing in $\mu_k$ can be determined using Table~\ref{fig.minimal.table}. In all cases, the sets $\Delta_\chi$ for $\chi\in \hat{Z}$ nontrivial are those corresponding to the black nodes of the Dynkin diagrams in the last column of Table~\ref{fig.Levi.table} (in type $A_n$, recall that the nodes corresponding to the boxed part of the Dynkin diagram are viewed as colored in black). It remains to prove part (2) in types besides type $A_n$; we do this below.   As above,
we use the notational conventions of~\cite{HumphreysLieAlg} in numbering the simple roots and fundamental dominant weights.

\subsection{Types $B_n$, $C_n$, and $E_7$.}
In types $B_n$, $C_n$, and $E_7$, we have $Z \cong \Z_2$.
There is only one nontrivial element $z$ in $Z$ and one nontrivial character
$\chi$ in $\hat{Z}$, and $\chi(z) = -1$.  There is only one minuscule weight $\lambda_k$, and hence only one minimal element $\mu_k$, to consider. We have $\C[\toric] = \C[x_1, \ldots, x_n, v_k]/I$, where $\ga_1, \ldots, \ga_n$ are the simple roots, $x_i = e^{\alpha_i}$, and $v_k = e^{\mu_k}$.  By \eqref{e:zaction}, we have $z \cdot x_i = x_i$, and since
$\mu_k$ represents $\chi$, we have $v_k (z) = e^{\mu_k}(z) = \chi(z) = -1$.  This implies that the 
action of $z$ on $v_k$ is given by $z \cdot v_k = - v_k$.

The element $\mu_k$ is a generator of $\cp/\cq \cong \Z_2\cong \hat{Z}$, so $2 \mu_k \in \cq$.  If we write 
$$
\mu_k = \sum_{\alpha_i  \in \Delta_{\chi}}  a_i \ga_i,
$$
then $0 < a_i < 1$ and $2 a_i \in \Z$, so each $a_i = \frac{1}{2}$.  Hence
$$
\mu_k = \frac{1}{2} \sum_{\ga_i \in \Delta_{\chi}} \ga_i,
$$
which the reader can confirm using Table~\ref{fig.minimal.table}. This implies 
that in $\C[\toric]$, we have 
\begin{equation} \label{eqn.vtoxBCE}
v_k^2 = \prod_{ \alpha_i \in \Delta_{\chi}} x_i.
\end{equation}

We claim that $\toric_{\chi}$ consists of the elements $x$ in $\toric$ such that $v_k(x) \neq 0$.  Indeed, by
definition, $\toric_{\chi}$ consists of the $x\in \toric$
such that the stabilizer group $Z^x$ is a subset of $\ker \chi$.  
Since $\ker \chi$ is trivial, $\toric_{\chi}$ consists of the elements $x$ with trivial stabilizer.  This means $z \cdot x \neq x$, which is equivalent to the condition that $v_k(x) \neq 0$, proving the claim.

By definition,  $\toric_{ad, \chi}$ is the image of $\toric_{\chi}$ under the quotient map $\toric \to \toric_{ad}$.
Since $\toric_{\chi}$ is defined by the relation $v_k \neq 0$, equation~\eqref{eqn.vtoxBCE} implies that $\toric_{ad, \chi}$ is defined by $\prod_{ \alpha_i \in \Delta_{\chi}} x_i \neq 0$.    This proves part (2) of Proposition~\ref{p:Levi-chi} for types $B_n$, $C_n$, and $E_7$.


\subsection{Type $E_6$}
In type $E_6$, we have $Z \cong \Z_3$.  The minuscule weights $\lambda_{1},\lambda_{6}$, and the corresponding minimal
elements $\mu_{1}, \mu_{6}$, are recorded in Table~\ref{fig.minimal.table}. 
We have 
\[
\C[\toric] = \C[x_1, \ldots, x_6, v_1, v_6]/I,
\] 
where  $v_1^3 =  x_1 x_3^2 x_5 x_6^2$ and $v_6^3 = x_1^2 x_3 x_5^2 x_6$.
The group $\cp/\cq\simeq \hat{Z}$ consists of $\{0, \overline{\mu}_1, \overline{\mu}_6 \}$, subject to the relations $2 \overline{\mu}_1 = \overline{\mu}_6$ and $3 \overline{\mu}_1 = 3\overline{\mu}_6=0$.
Let $\chi_j$ (for $j = 1,6$) denote the character of $Z$ represented by $\mu_j$, so $\hat{Z} = \{ 1, \chi_1, \chi_6 \}$, with
$\chi_1^2 = \chi_6$ and $\chi_1^3 = 1$. Note that $\Delta_{\chi_1} = \Delta_{\chi_6} = \{ \ga_1, \ga_3, \ga_5, \ga_6 \}$.
The center $Z$ consists of $\{1,z, z^2 \}$, where $z^3 = 1$, and we can choose $z$ so that $\chi_1(z) = \omega = e^{2 \pi i/3}$.  Then by \eqref{e:zaction},
\begin{align*}
z \cdot v_1 &= e^{\mu_1}(z^{-1}) v_1 = \chi_1(z^{-1})v_1 = \omega^{-1}v_1 = \omega^2 v_1 \\
z \cdot v_6 &= e^{\mu_6}(z^{-1})v_6 = \chi_6(z^{-1})v_6 = \omega^{-2}v_6 =  \omega v_6.
\end{align*} 
If $\chi$ is either $\chi_1$ or $\chi_6$, then $\ker \chi$ is trivial.  Hence, as in the previous case, we get $x\in \toric_\chi$ if and only if $x$ has a trivial stabilizer subgroup in $Z$, which in turn holds if and only if either $v_1(x)$ or $v_6(x)$ is nonzero.
From the expressions for $v_1^3$ and $v_6^3$ in terms of the $x_i$, we see that these conditions are equivalent to $x_k \neq 0$ for
$k = 1,3,5,6$.  Hence, $\toric_{ad, \chi}$ is defined by $x_k \neq 0$ for $k=1,3,5,6$, proving part (2) of Proposition~\ref{p:Levi-chi} for type~$E_6$.


\subsection{Type $D_n$} \label{ss:Dn}
We now consider type $D_n$ with $n\geq 4$.  Recall that $T_{ad} =  T/Z$.  Let $y_i=e^{\lambda_i}$ and $x_i=e^{\alpha_i}$.
The $y_i$ may be viewed as coordinates on $T$, and the $x_i$ as coordinates on $T_{ad}$.  The following lemma describes $Z$ in terms of the coordinates $y_i$.  For lack of a simple reference, we give a proof in the appendix.  One can avoid the use of this lemma, using only the existence of generators $z_i$ of $Z$, but we include it to make the arguments more concrete. 

\begin{Lem}\label{lem.ZintypeD} Let $y_i=e^{\lambda_i}$ and $x_i=e^{\alpha_i}$ for each $i$, and recall that $Z$ is the kernel of the natural projection $T \to T_{ad}$ defined by $(y_1, \ldots, y_n)\mapsto (x_1, \ldots, x_n)$. 
\begin{enumerate} 
\item If $n$ is even, then $Z\simeq \Z_2\times \Z_2$, with generators given by
\[
w_1 := (-1,1,-1, \ldots, -1,1,1,-1) \; \text{ and } \; w_2 := (-1,1, -1,\ldots, -1,1,-1,1).
\]
Here both $w_1$ and $w_2$ have entry $-1$ in positions $1,3, \ldots, n-3$; they end with $(1,1,-1)$ and $(1,-1,1)$, respectively.
\item If $n$ is odd, then $Z\simeq \Z_4$ with generator 
\[
z_0 := (-1,1,-1, \ldots,1,-1,\omega,\omega^3) 
\]
where $\omega$ is a primitive fourth root of unity.
\end{enumerate}
\end{Lem}

When $n$ is even, we have $Z\simeq \Z_2\times \Z_2$.  Since we make use of results of Bonnaf\'e in~\cite{Bonnafe2004}, we need
to align the generators from statement (1) of Lemma~\ref{lem.ZintypeD} with those used in that reference. In order to do so, we set some notation.
We write $Z=\{1, z_1, z_{n-1}, z_n\}$, where we set $z_1=w_1w_2$ and 
\begin{align*}
z_{n-1} &= w_1 \; \text{ and }\; z_{n}=w_2 \; \quad \text{if $n \equiv 0 \textup{ mod } 4$, and  }\\
z_{n-1} &= w_2 \; \text{ and }\; z_{n}=w_1 \; \quad \text{if $n \equiv 2 \textup{ mod } 4$.  }
\end{align*} 
The following is proved in Appendix~\ref{appendix.Bonnafe}.

\begin{Lem} \label{lem:Bonnafe}
Let $G$ be of type $D_n$ for $n$ even, and let $z_1, z_{n-1}, z_n$ be as defined above.  Then $z_1 = \ibarT(\varpi_1^{\vee})$,  $z_{n-1}=\ibarT(\varpi_{n-1}^{\vee})$, and  $z_{n}=\ibarT(\varpi_{n}^{\vee})$, where $\ibarT(\varpi_1^{\vee})$, $\ibarT(\varpi_{n-1}^{\vee})$, and  $\ibarT(\varpi_{n}^{\vee})$ are the nontrivial elements of $Z$ in the notation of \cite{Bonnafe2004}.
\end{Lem}

In order to compute $\toric_{ad, \chi}$, we introduce some additional notation. We have
$$
\cp/\cq = \{ 0, \overline{\mu}_1, \overline{\mu}_{n-1}, \overline{\mu}_n \}.
$$
Let $\chi_j$ (for $j = 1, n-1, n$) denote the character of $Z$ represented by $\mu_j$.  If $z \in Z$,
$$
\chi_j(z) = e^{\mu_j}(z) = e^{\gl_j}(z) = y_j(z).
$$

\begin{Lem} \label{lem.kernel}
Let $Z_0 \subset Z $ denote the kernel of the covering map $\Spin(2n) \to \SO(2n)$.  
Then
$$
Z_0 = \ker \chi_1 = \begin{cases} \{ 1, z_1 \} & \mbox{if} \  \  n  \ \ \mbox{is even} \\
\{1, z_0^2 \} &\mbox{if} \  \  n  \ \ \mbox{is odd}.
\end{cases}
$$
The characters $\chi_{n-1}$ and $\chi_n$ take the value $-1$ on the nontrivial element of $Z_0$.
\end{Lem}

\begin{proof}
The standard representation $V$ of $\SO(2n)$ has highest weight $\gl_1$.  The subgroup $Z_0$ consists of the elements of $Z$ which
act trivially on the standard representation.  An element $z \in Z$ acts on $V$ as multiplication by the scalar $e^{\gl_1}(z) =  \chi_1(z)=y_1(z)$.
It follows that $Z_0 = \ker \chi_1$.  The coordinates of the elements in the center $Z$ are determined in Lemma \ref{lem.ZintypeD}.  By inspection, we see that the nontrivial  elements $z \in Z$ such that $y_1 (z) = 1$ are precisely $z=w_1w_2$ when $n$ is even and $z=z_0^2$ when $n$ is odd.   The fact that $\chi_{n-1}$ and $\chi_n$ take the value $-1$ on this nontrivial element of $Z_0$ follows since the $y_{n}$ and $y_{n-1}$ coordinates of both $w_1w_2$ and $z_0^2$ are equal to $-1$. 
\end{proof}

We recall from  Table~\ref{fig.minimal.table} the formulas for minimal elements $\mu_1, \mu_{n-1},\mu_n$.  We have $\mu_1 =  \frac{1}{2}(\alpha_{n-1}+\alpha_n)$, and, writing 
$B = \ga_1 + \ga_3 + \cdots + \ga_d$, where $d = n-2$ if $n$ is odd, and $d = n-3$ if $n$ is even, we have:
\begin{itemize}
\item If $n \equiv 1 \textup{ mod }4$ then $\mu_{n-1} =  \frac{1}{2} B +  \frac{1}{4} (\alpha_{n-1}+3 \alpha_n)$ and $\mu_n =  \frac{1}{2} B + \frac{1}{4}( 3 \alpha_{n-1}+\alpha_n)$.

\item If $n \equiv 3 \textup{ mod }4$ then 
$\mu_{n-1} =  \frac{1}{2}B  +   \frac{1}{4}( 3 \alpha_{n-1}+\alpha_n)$ and $\mu_n =  \frac{1}{2} B + \frac{1}{4} (\alpha_{n-1}+3 \alpha_n)$.

\item If $n \equiv 2 \textup{ mod }4$ then $\mu_{n-1} =  \frac{1}{2} (B +  \alpha_{n-1})$ and $\mu_{n} =  \frac{1}{2} (B +  \alpha_n)$.

\item If $n \equiv 0 \textup{ mod }4$ then $\mu_{n-1} =  \frac{1}{2} (B +  \alpha_n)$ and $\mu_{n} =  \frac{1}{2} (B +  \alpha_{n-1})$.
\end{itemize}

We are now ready to verify part (2) of Proposition \ref{p:Levi-chi} for type $D_n$.  We consider the cases of $n$ odd and $n$ even separately. 

\textbf{Case 1: $n$ is odd.} We have $\cp/\cq \cong \Z_4$. Both ${\mu}_{n-1}$ and ${\mu}_n$ are generators for $\cp/\cq$, and, correspondingly, $\hat{Z}$ is generated by $\chi_{n-1}$ and $\chi_n$.  
Let $z_0$ be the generator of $Z$ from Lemma \ref{lem.ZintypeD}(2). We have $\C[\toric] = \C[x_1, \ldots, x_n, v_1, v_{n-1}, v_{n}]/I$, and the following relations hold:
\begin{align*} 
v_1^2 & =  x_{n-1} x_n \\
v_{n-1}^4 & =  (x_1 x_3 \cdots x_{n-2})^2 x_{n-1}^a x_n^b \\
v_{n}^4 & =  (x_1 x_3 \cdots x_{n-2})^2 x_{n-1}^b x_n^a
\end{align*}
where $\{a,b \} = \{ 1,3 \}$.
Since the action of $Z$ on $f \in \C[\toric]$ is given by $(z \cdot f)(x) = f(z^{-1} x)$ for $z \in Z$, $x \in \toric$, Lemma \ref{lem.ZintypeD}(2) implies that $z_0 \cdot v_1 = e^{\mu_1}(z_0^{-1}) v_1 = - v_1$.  Similarly, $z_0 \cdot v_{n-1} = \omega^3 v_{n-1}$, and $z_0 \cdot v_n = \omega v_n$. 

If $\chi$ has order $4$ (i.e., if $\chi$ equals $\chi_{n-1}$ or $\chi_n$), then $\ker \chi = \{ 1 \}$.  Thus, an element $x \in \toric$ is in $\toric_{\chi}$ if and only if $x$ has trivial stabilizer subgroup in $Z$, which in turn holds if and only if $v_{n-1}(x)\neq 0$ or $v_{n}(x)\neq 0$.
From the expressions for $v_{n-1}^4$ and $v_n^4$ above, we see that $\toric_{ad, \chi}$ is defined by the conditions $x_k \neq 0$ if $k$ is odd and $1\leq k \leq n-2$, or if $k \in  \{ n-1, n \}$.  Using our formulas for $\mu_{n-1}$ and $\mu_n$ above, we see that these $k$ are precisely the indices such that $\ga_k \in \Delta_{\chi}$.  Hence part (2) of Proposition~\ref{p:Levi-chi} holds when $\chi\in \hat{Z}$ has order $4$.

If $\chi$ has order $2$ (i.e., if $\chi = \chi_1$) then $\ker \chi = \{1, z_0^2 \}$, so $x$ is in $\toric_{\chi}$ if and only if $z_0 x \neq x$ and $z_0^3 x \neq x$.  This holds if and only if $v_1(x) \neq 0$.   From the equation $v_1^2  =  x_{n-1} x_n$ we see that $\toric_{ad, \chi}$ is defined by the conditions $x_{n-1}\neq 0$ and  $x_n \neq 0$.  Since $\Delta_{\chi} = \{ \alpha_{n-1}, \alpha_n \}$, we again get that part (2) of Proposition~\ref{p:Levi-chi} holds.

\textbf{Case 2: $n$ is even.} We have $\cp/\cq \cong \Z_2 \times \Z_2$. Each of the elements $\overline{\mu}_k$ has order $2$, for $k\in \{1,n-1,n\}$, and the sum of any two distinct elements equals the third.  The relations  in the ideal $I$ defining $\C[\toric] = \C[x_1, \ldots, x_n, v_1, v_{n-1}, v_{n}]/I$ depend on the value of $n$ modulo $4$.  We have
\begin{itemize}
\item $v_1^2  =  x_{n-1} x_n$ in either case.
\item If $n \equiv 0 \textup{ mod }4$, then $v_{n-1}^2  = x_1 x_3 \cdots x_{n-3} x_n$ and $v_n^2  = x_1 x_3 \cdots x_{n-3} x_{n-1}$.
\item If $n \equiv 2 \textup{ mod }4$, then $v_{n-1}^2  = x_1 x_3 \cdots x_{n-3} x_{n-1}$ and $v_n^2  = x_1 x_3 \cdots x_{n-3} x_n$. 
\end{itemize}
We have $\hat{Z} = \{ 1, \chi_1, \chi_{n-1}, \chi_n \}$.  Also, for $j \in \{1, n-1, n \}$, the character $\chi_j$ is represented by $e^{\mu_j} = v_j$, and also by
$e^{\lambda_j} = y_j$.  In other words, $\chi_j(z) = v_j(z) = y_j(z)$.  By Lemma \ref{lem.ZintypeD}(1), 
the elements of $\hat{Z}$ take the following values on the generators $w_1, w_2$ of $Z$:
\begin{equation} \label{e:gen-z-group}
\begin{aligned}
  \chi_1(w_1) &= \chi_1(w_2) = -1 \\
  \chi_{n-1}(w_1) &= 1, \quad \chi_{n-1}(w_2) = -1 \\
  \chi_n(w_1) &= -1, \quad \chi_n(w_2)= 1.
\end{aligned}
\end{equation}
Note that the center $Z$ acts on the coordinate function $v_j$ by the character $\chi_j^{-1} = \chi_j$.

To verify part (2) of the proposition, we consider each of the possibilities for $\ker \chi$.  
If $\ker \chi =  \langle z_1 \rangle = \{1, w_1 w_2 \} $, then $\chi = \chi_1$ by~\eqref{e:gen-z-group}.  By definition, $x \in \toric_{\chi_1}$ if and only if $Z^x \subset \ker \chi_1 = \{1, w_1 w_2 \}$.  In particular, we must have both $w_1 x\neq x$ and $w_2x\neq x$. The condition $w_1 x \neq x$ implies that either $v_1(x) \neq 0$ or $v_n(x) \neq 0$; the condition $w_2 x \neq x$ implies that either $v_1(x) \neq 0$ or $v_{n-1}(x) \neq 0$.  Thus, either $v_1(x) \neq 0$ or both $v_{n-1}(x)$ and $v_n(x)$ are nonzero.
In the latter case, since $v_1^2$ divides $v_{n-1}^2 v_n^2$, we see that $v_1(x) \neq 0$ as well.  Thus, $\toric_{\chi_1}$ is defined by
the relation $v_1(x) \neq 0$, which is equivalent to $x_{n-1} \neq 0$ and $x_n \neq 0$.  Since $\Delta_{\chi_1} = \{ \ga_{n-1}, \ga_n \}$, we see
part (2) of Proposition~\ref{p:Levi-chi} holds in this case.

Now suppose $\ker \chi = \langle z_{n-1} \rangle= \{1, z_{n-1}\}$.  By~\eqref{e:gen-z-group} and the definition of $z_{n-1}$, we have $\chi=\chi_{n-1}$ when $n \equiv 0 \textup{ mod }4$ and $\chi= \chi_n$ when $n \equiv 2 \textup{ mod }4$. Our computations for $\mu_n$ and $\mu_{n-1}$ above show that $\Delta_\chi = \{\alpha_1, \alpha_3, \ldots, \alpha_{n-3}, \alpha_n\}$, and thus we must check that $\toric_{ad, \chi}$ is defined by the relations $x_1 \neq 0$, $x_3 \neq 0$, \ldots, $x_{n-3} \neq 0$, $x_n \neq 0$.

First consider the case $n \equiv 0 \textup{ mod }4$, so $z_{n-1} = w_1$ and $\chi = \chi_{n-1}$. By definition, $x \in \toric_{\chi_{n-1}}$ if and only if $Z^x \subset \ker \chi_{n-1} = \{1, w_1 \}$.  In particular, we must have $w_2x\neq x$ and $w_1w_2x\neq x$.  The condition $w_2 x \neq x$ implies that either $v_1(x) \neq 0$ or $v_{n-1}(x) \neq 0$;
the condition $w_1 w_2 x \neq x$ implies that either $v_{n-1}(x) \neq 0$ or $v_n(x) \neq 0$.  Hence, either $v_{n-1}(x) \neq 0$, or both
$v_1(x)$ and $v_n(x)$ are nonzero.  In the latter case, since $v_{n-1}^2$ divides $v_1^2 v_n^2$, we see that $v_{n-1}(x) \neq 0$ as well.
Thus $\toric_{\chi_{n-1}}$ is defined by the relation $v_{n-1}(x) \neq 0$.  Since $n \equiv 0 \textup{ mod }4$, this is equivalent to the relations
$x_1 \neq 0$, $x_3 \neq 0$, \ldots, $x_{n-3} \neq 0$, $x_n \neq 0$, which is what we wanted to show.

Now consider the case $n \equiv 2 \textup{ mod }4$, so $z_{n-1} = w_2$ and $\chi = \chi_n$.  By definition, $x \in \toric_{\chi_n}$ if and only if $Z^x \subset \ker \chi_n = \{1, w_2 \}$.   In particular $w_1x\neq x$ and $w_1w_2x\neq x$.  The condition $w_1 x \neq x$ implies that either $v_1(x) \neq 0$ or $v_n(x) \neq 0$; the condition $w_1 w_2 x \neq x$ implies that either $v_{n-1}(x) \neq 0$ or $v_n(x) \neq 0$.
Hence, either $v_n(x) \neq 0$, or both $v_1(x)$ and $v_{n-1}(x)$ are nonzero.  In the latter case, since $v_n^2$ divides $v_1^2 v_{n-1}^2$, we see that $v_n(x) \neq 0$ as well.
Thus,  $\toric_{\chi_n}$ is defined by the relation $v_n(x) \neq 0$.  Since $n \equiv 2 \textup{ mod }4$, this is equivalent to the relations
$x_1 \neq 0$, $x_3 \neq 0$, \ldots, $x_{n-3} \neq 0$, $x_n \neq 0$, which is what we wanted to show.  This proves part (2) of Proposition~\ref{p:Levi-chi} in the case $\ker \chi = \langle z_{n-1} \rangle$.

The final possibility is when $\ker \chi = \langle z_n \rangle = \{1,z_n\}$.  The argument is almost the same as when $\ker \chi = \langle z_{n-1} \rangle$, and we omit it. This completes the proof of Proposition~\ref{p:Levi-chi}.

\section{Lusztig sheaves and the extended Springer resolution} \label{s.Lusztig}
In this section we prove Theorem \ref{t:main1}, which describes the pushforward of the constant sheaf on the extended Springer resolution as a direct sum of Lusztig sheaves.
This extends results of type $A$ from \cite{GPR23}.  The proof is similar to the proofs in \cite{GPR23}, and we will refer to that paper for some details.

\subsection{Minimal cuspidal data} \label{ss:minimal}
Following~\cite{AHJR3}, associated to each character $\chi$ of the center $Z$ there is a triple
$$
\vc_{\chi} = (L_{\chi},\co^{\reg}_{L_{\chi}}, \cl^{L_\chi}),
$$
called the \emph{minimal cuspidal datum} for $\chi$ (see Proposition 5.1 in loc.~cit.). Here $L_\chi$ is a standard Levi subgroup and $\cl^{L_\chi}$ is an $L_{\chi}$-equivariant local system on the regular nilpotent orbit $\co^{\reg}_{L_{\chi}}$ of $L_{\chi}$.  This local system is defined more precisely in Remark~\ref{r:local} below.  

The definition of $L_{\chi}$ in \cite[Prop.~5.1]{AHJR3} differs from the definition of $L_\chi$ in Proposition ~\ref{p:Levi-chi} of this paper.
However, the next proposition implies that our notation is consistent with that of~\cite{AHJR3}.
More precisely, they note that the subgroup they define is characterized up to $G$-conjugacy by a certain property, 
and the next proposition shows that the subgroup defined in Proposition~\ref{p:Levi-chi} satisfies that property. 

Before stating this proposition, we make some definitions. Given a subgroup $K$ of $Z = \cz(G)$, let $\sfl(K)$ denote the set of Levi subgroups $L$ of $G$ 
such that  $\ker h_L \subseteq K$.
For a fixed $K$, the elements of $\sfl(K)$ are ordered by inclusion; let $\mathsf{L}_{min}(K)$ denote the minimal elements of 
$\sfl(K)$.   
The elements of $\sfl_{min}(\ker \chi)$ are all $G$-conjugate by \cite[Lemma 2.16]{Bonnafe2004}, so $\sfl_{min}(\ker \chi)$ 
contains a standard Levi subgroup.  These standard Levi subgroups are described in \cite[Table 2.17]{Bonnafe2004}.
From this description, we can see that each of these Levi subgroups
is the unique standard Levi subgroup in its conjugacy class.  The reason is that
two standard Levi subgroups are $G$-conjugate exactly when they correspond to subsets of simple roots which are
conjugate under the Weyl group (see \cite{Solleveld2020}).  By inspection,
each of the subsets of simple roots in the table is $W$-conjugate to no other set of simple roots.

\begin{Prop} \label{p:L-chi}
Let $\chi$ be a character of $Z$.  The Levi subgroup $L_{\chi}$ defined in Proposition~\ref{p:Levi-chi} has the property that for any Levi subgroup $M$ of $G$,  the character $\chi$ factors through the surjection $h_M: Z \to \cz(M)$ if and only if $M$ contains a conjugate of $L_{\chi}$.   Equivalently, $L_{\chi}$ is the standard Levi subgroup in $\sfl_{min}(\ker \chi)$.  Hence $L_{\chi}$ coincides with the Levi subgroup in the minimal cuspidal datum~$\vc_{\chi}$ for $\chi$.
\end{Prop}

\begin{proof}
The character $\chi$ factors through $h_M$ if and only if $\ker h_M \subset \ker \chi$.  This means
$M \in \sfl(\ker \chi)$, so $M$ contains an element of $\sfl_{min}(\ker \chi)$.  
For any subgroup $K$ of $Z$, \cite[Table 2.17]{Bonnafe2004} lists the standard Levi subgroup in $\sfl_{min}(K)$.  The Levi subgroups $L_{\chi}$ are recorded in Table~\ref{fig.Levi.table}, and comparing with \cite[Table 2.17]{Bonnafe2004} using Lemma~\ref{lem:Bonnafe} shows that each $L_{\chi}$ is in $\sfl_{min}(\ker \chi)$.
By \cite[Lemma 2.16]{Bonnafe2004}, any two elements of $\sfl_{min}(\ker \chi)$ are $G$-conjugate.  Hence the character $\chi$ factors through $h_M$ if and only if $M$ contains a conjugate of $L_{\chi}$. 
The discussion above shows that this condition is equivalent to the statement that $L_{\chi}$ is the standard Levi subgroup in $\sfl_{min}(\ker \chi)$.
The last statement of the proposition follows from \cite[Prop.~5.1]{AHJR3}.
\end{proof}

\begin{Rem} \label{r:local}
Proposition~\ref{p:Levi-chi}(1) implies that $\chi = h_{L_{\chi}}^*(\chi_{L_\chi})$, where $\chi_{L_\chi}$ is a uniquely determined character of $\cz(L_{\chi})$.  Characters of $\cz(L_{\chi})$ correspond to $L_{\chi}$-equivariant local systems on the regular
nilpotent orbit $\co^{\reg}_{L_{\chi}}$.  The local system on $\co^{\reg}_{L_{\chi}}$ corresponding to
$\chi_{L_\chi}$ is denoted  $\cl^{L_\chi}$.
\end{Rem}

\begin{Rem} \label{r:bonnafe}
The rows of Table~\ref{fig.Levi.table} omit one row of \cite[Table 2.17]{Bonnafe2004}, but are otherwise the same.
The reason is as follows.  For each subgroup $K$ of $Z$, Bonnaf\'e records the standard Levi subgroup in $\sfl_{min}(K)$.
If $K = \ker \chi$, then Proposition~\ref{p:L-chi} shows that $L_{\chi}$ coincides with the subgroup in Bonnaf\'e's table.
There is only one case where the subgroup $K$ is not of the form $\ker \chi$.  This occurs when $G$ is of type $D_{2r}$,
and $K = \{ 1 \}$ is the trivial subgroup.  In this case, $Z=\cz(G) \cong \Z_2 \times \Z_2$, so there is no character $\chi\in \hat{Z}$ whose kernel is
the trivial group $K$.  Thus, the row corresponding to this case appears in \cite[Table 2.17]{Bonnafe2004} but not in Table \ref{fig.Levi.table}.
\end{Rem}


\subsection{Cuspidal data and Lusztig sheaves} \label{ss.cuspidal}
We briefly describe the construction of the Lusztig sheaf associated to a cuspidal datum; more details can be found in \cite{GPR23}. 

Let $L$ be a standard Levi subgroup of $G$. Let $P = L U_P$ denote the standard parabolic subgroup containing $B$ with Levi factor $L$, and let $\fp = \fl + \fu_P$ denote the associated Levi decomposition of its Lie algebra.  Following Borho--MacPherson~\cite{Borho1981}, set $\tilde{\cn}^P = G \times^P (\cn_L + \fu_P)$ where $\cn_L$ is the nilpotent cone in $\fl$.  We define $\tilde{\cn}^P_{\reg} = G \times^P (\co_L^{\reg} + \fu_P)$, where $\co_L^{\reg}$ is the regular nilpotent orbit of $L$.  There are natural maps which factor the Springer resolution
\begin{equation}\label{eqn.Springer-factor}
\tilde{\cn} \xrightarrow{\;\;\rho_P\;\;} \tilde{\cn}^P \xrightarrow{\;\;\mu_P\;\;} \cn,
\end{equation}
defined by
$$
[g,x]_B \mapsto [g,x]_P \mapsto \Ad(g) x.
$$
These maps are isomorphisms over the regular nilpotent orbit $\co^{\reg}$ in $\cn$.

Recall that a cuspidal datum is a triple $\vc = (L,\co_L, \cl^L)$ satisfying certain conditions.  Here, $L$ is a Levi subgroup of $G$ and $\cl^L$ is a local system on a nilpotent orbit $\co_L$ of $L$. 
From a cuspidal datum $\vc$, Lusztig constructs a perverse sheaf $A_{\vc}$ on $\cn$ as follows.  Consider the diagram
\[
\cn\xleftarrow{\;\mu_P\;} \widetilde{\cn}^P:=G\times^P (\cn_L+\fu_P)\xrightarrow{\;\pi_P\;} G\times^P \cn_L{\longrightarrow} \cn_L.
\]
There is an induced perverse sheaf $\Ind_P^G IC(\overline{\co}_L, \cl^L)$ on $G\times^P \cn_L{\longrightarrow} \cn_L$.  Therefore,
\[
A_{\vc} := (\mu_P)_* (\pi_P)^* \Ind_P^G IC(\overline{\co}_L, \cl^L)
\] 
is a perverse sheaf on $\cn$, which we call the \emph{Lusztig sheaf} associated to $\vc$.
Of greatest importance in this paper are the Lusztig sheaves $A_{\vc_{\chi}}$ associated to minimal cuspidal data.

\subsection{Proof of Theorem \ref{t:main1}} Recall that this theorem states that 
$$
\psi_* \shc_{\tilde{\cm}}[\dim \cn] = \bigoplus_{\chi \in \hat{Z}} \A_{\vc_{\chi}},
$$
where $\psi$ is the extended Springer resolution from \eqref{e:extended}:
$$
\psi: \tilde{\cm} \stackrel{\tilde{\eta}}{\longrightarrow} \tilde{\cn} \stackrel{\mu}{\longrightarrow} \cn.
$$
The map $\tilde{\eta}$ is the quotient by the action of $Z$.
The proof below uses
the factorization of $\mu$ from~\eqref{eqn.Springer-factor}, as well as the projection $\pi_P: \tilde{\cn}^P \to G\times ^P \cn_L$.  
We divide the proof into steps. 

\medskip

{\em Step 1.}  Let $\cl_{\chi}$ be the local system on $\co^{\reg}\subset \cn$ corresponding to
$\chi \in \hat{Z}$. We have 
$$
\tilde{\eta}_* \shc_{\tilde{\cm}}[\dim \cn] = \bigoplus_{\chi \in \hat{Z}} IC(\tilde{\cn}, \cl_{\chi})
$$
by \cite[Prop.~3.19]{GPR23} (which was stated in loc.~cit.~only for $G$ of type $A_n$, but is true in all types by the same proof). 

\medskip

{\em Step 2.} The regular nilpotent orbit $\co^{\reg}$ can be identified with an open subset of $\tilde{\cn}$.   
There is a larger open subset  $\tilde{\cn}_{\chi}$ of $\tilde{\cn}$ such that the local system $\cl_{\chi}$ on $\co^{\reg}$ extends to a local system (still denoted $\cl_{\chi}$) on $\tilde{\cn}_{\chi}$.

The subset $\tilde{\cn}_{\chi}$ of $\tilde{\cn}$ is defined as follows (cf.~\cite[Prop.~3.17]{GPR23}).  We have open embeddings $\toric_{ad, \chi} \subset \toric_{ad}$ and $\toric_{\chi} \subset \toric$.  We define $\tilde{\fu}_{\chi} = {\toric_{\chi}} \times_{\toric_{ad}} \fu$, an open subset of $\tilde{\fu}$.  Similarly, $\fu_{\chi} = \toric_{ad, \chi} \times_{\toric_{ad}} \fu$.  This is identified with the open subset of $\fu$ consisting of elements of the form $\sum_{\ga \in \Phi^+} a_{\alpha} X_{\alpha} $ where $ a_{\ga} \neq 0$ for $\ga \in \Delta_{\chi}$.
We define $\tilde{\cm}_{\chi} = G \times^B \tilde{\fu}_{\chi}$ and $\tilde{\cn}_{\chi} = G \times^B \fu_{\chi}$.  These are open subsets of $\tilde{\cm}$ and $\tilde{\cn}$, respectively.
Since $\co^{\reg} \cap \fu$ consists of the elements  $\sum_{\ga \in \Phi^+} a_{\alpha} X_{\alpha} $ where $ a_{\ga} \neq 0$ for all $\ga \in \Delta$, we see that $\co^{\reg} \cap \fu$ is an open subset of $\fu_{\chi}$.  Therefore, $\co^{\reg}$ is an open subset of $\tilde{\cn}_{\chi}$.

The universal cover $\tilde{\co}^{\reg}$ of the regular nilpotent orbit is identified with the inverse image of $\co^{\reg}$ under the quotient by $Z$, $\tilde{\eta}: \tilde{\cm} \to \tilde{\cn}$. The action of $Z$ on $\toric$ induces an action of $Z$ on $\tilde{\cm}$, which restricts to a free action on the open set $\tilde{\co}^{\reg}$. Moreover, $\tilde{\cm}_{\chi}$ is $Z$-invariant, and $\tilde{\cm}_{\chi}/Z = \tilde{\cn}_{\chi} $. The local system $\cl_{\chi}$ on $\co^{\reg}$ corresponds to the line bundle $\tilde{\co}^{\reg} \times^Z \C_{\chi}$. By definition, $\toric_{\chi}$ is the set of all $x \in \toric$ such that $Z^x \subset \ker \chi$.  Thus, $\tilde{\cm}_{\chi}$ is the set of $x \in \tilde{\cm}$ such that $Z^x \subset \ker \chi$.  Hence, by~\cite[Prop.~3.5]{GPR23},  the local system $\cl_{\chi}$ on $\co^{\reg}$ extends to a local system on $\tilde{\cn}_{\chi} $.

\medskip

{\em Step 3.} Let $L = L_{\chi}$ and let $P = L U_P$ be the corresponding parabolic subgroup.  We have a commutative
diagram
$$
\begin{tikzcd}
 \tilde{\cn}_{\chi}  \arrow[hook, r] \arrow[d,  "\rho_P|_{\tilde{\cn}_\chi}"] & 
\tilde{\cn} \arrow[d,  "\rho_P"] \\
 \tilde{\cn}^P_{\reg} \arrow [hook, r] &
\tilde{\cn}^P 
\end{tikzcd}
$$
and we will prove below that the left vertical arrow is an isomorphism.  Hence the local system $\cl_{\chi}$ can be viewed as a local
system on $ \tilde{\cn}^P_{\reg}= G\times^P \left( \co_L^{\reg} +\fu_P \right)$.  

The fact that the left vertical arrow is an isomorphism is a consequence of the next two lemmas. The first extends \cite[Lemma 4.4]{GPR23} to all types.

\begin{Lem} \label{l:intersection}
With notation as above, $(\co_L^{reg} + \fu_P) \cap \fu = \fu_{\chi}$.
\end{Lem}

\begin{proof}
Let $\fu_L$ denote the Lie algebra of $U\cap L$. 
The intersection $\co_L^{\reg} \cap \fu_L$ is the set of elements in $\fl$ of the form $ \sum_{\ga \in \Phi^+_L} a_{\alpha} X_{\alpha}$ where $a_{\ga} \neq 0$ for $\ga \in \Delta_L$. Because $\fu = \fu_L \oplus \fu_P$, we have
$$
(\co_L^{\reg} + \fu_P) \cap \fu = (\co_L^{\reg} \cap \fu_L) + \fu_P = \left\{  \sum_{\ga \in \Phi^+} a_{\alpha} X_{\alpha} \mid a_{\ga} \neq 0 \mbox{ for }\ga \in \Delta_L \right\}.
$$
Since $\Delta_L = \Delta_{\chi}$, we see that $(\co_L^{\reg} \cap \fu_L) + \fu_P = \fu_{\chi}$, as desired.
\end{proof}

The next lemma generalizes \cite[Prop.~4.5]{GPR23} to all types.
\begin{Lem} \label{l:openisom}
The map $\rho_P: \tilde{\cn} \to \tilde{\cn}^P$ restricts to an isomorphism of $\tilde{\cn}_{\chi} \to \tilde{\cn}^P_{\reg}$.  Moreover,
$\rho_P^{-1}(\tilde{\cn}^P_{\reg}) = \tilde{\cn}_{\chi}$.
\end{Lem}

\begin{proof}
The proof of \cite[Prop.~4.5]{GPR23} applies, using Lemma~\ref{l:intersection} above in place of \cite[Lemma 4.4]{GPR23}.
\end{proof}

\medskip

{\em Step 4.} 
We have $(\rho_P)_* IC(\tilde{\cn}, \cl_{\chi}) = IC(\tilde{\cn}^P, \cl_{\chi})$.  The analogous statement in type $A$ is \cite[Theorem~4.6]{GPR23}; the proof there extends to all types.

\medskip

{\em Step 5.}  
By Proposition \ref{p:Levi-chi}(1), $\chi = h_L^*(\chi_L)$ for $\chi_L \in \widehat{\cz(L)}$. Let $\cl^L$ denote the local system on $\co^{\reg}_L$ corresponding to the character $\chi_L$ (c.f.~Remark~\ref{r:local}). Let $\nilp_L$ denote a regular nilpotent element in $\fl$ and  recall that $\nilp$ denotes a regular nilpotent element in~$\fg$.

We claim that $\pi_P^* \Ind_P^G IC(\cn_L, \cl^L) = IC(\tilde{\cn}^P, \cl_{\chi})$.  
This statement is the analogue of equation (5.1) of \cite{GPR23}; we give the argument, which is similar.  
As in \cite{GPR23}, the claim reduces to showing the following statement: If $q: \co_L^{reg} + \fu_P \to \co_L^{reg}$ denotes the projection, then the local systems $q^* (\cl^L)$  and $\cl_{\chi}$ have isomorphic restrictions to the $P$-orbit $P \cdot \nilp$.  (By the argument of \cite[Remark 5.2]{GPR23}, $P \cdot \nilp$ is open in $\co^{\reg}_L + \fu_P$.)

The map $q$ induces (by restriction) a map of orbits
\begin{equation} \label{e.orbitmap}
q: P \cdot \nilp \to L \cdot \nilp_L = \co^{\reg}_L.
\end{equation}
We need to check that 
\begin{equation} \label{e.pullbackchi}
q^* (\cl^L) \cong \cl_{\chi}|_{P \cdot \nilp}.
\end{equation}    
By definition, the local system $\cl^L$ corresponds to the representation $\chi_L$ of the component group  $L^{\nilp_L}/L^{\nilp_L}_0 = \cz(L)$.  
The pullback local system $q^* (\cl^L) $ corresponds to the character of the component group $P^{\nilp}/P^{\nilp}_0 = Z$
which is the pullback of $\chi_L$ under the map of component groups induced by $q$.  Since this is the map $h_L$, the corresponding character of $Z$ is indeed $h_L^*(\chi_L) = \chi$.  Therefore, the pullback local system is indeed $\cl_{\chi}|_{P \cdot \nilp}$ and~\eqref{e.pullbackchi} follows.

\medskip

{\em Step 6.}
By Steps 4 and 5, as perverse sheaves on $\cn$, we have
\begin{align*}
\mu_*  IC(\tilde{\cn}, \cl_{\chi}) & = (\mu_P)_*(\rho_P)_* IC(\tilde{\cn}, \cl_{\chi}) = (\mu_P)_* IC(\tilde{\cn}^P, \cl_{\chi}) \\
& = (\mu_P)_* \pi_P^* \Ind_P^G IC(\cn_L, \cl^L) = \A_{\vc_{\chi}}.
\end{align*}

\medskip

{\em Step 7.}  By Steps 1 and 6, 
$$
\psi_* \shc_{\tilde{\cm}}[\dim \cn]  = \mu_* \tilde{\eta}_*  \shc_{\tilde{\cm}}[\dim \cn] = \bigoplus_{\chi \in \hat{Z}} \mu_* IC(\tilde{\cn}, \cl_{\chi}) \\
= \bigoplus_{\chi \in \hat{Z}} \A_{\vc_{\chi}},
$$
completing the proof.

\section{A support theorem for the extended Springer resolution} \label{sec.g2.intersection}

Recall that each pair $(\co, \cl)$ where $\co$ is a nilpotent orbit and $\cl$ an equivariant local system on $\co$ defines a simple perverse sheaf $IC(\overline{\co}, \cl)$ on $\cn$. In this section, we classify the collection of all such sheaves that appear as summands of $\psi_*\shc_{\tilde{\cm}}[\dim \cn]$. Recall that $\fg_2 = \toric_{ad}$ is the sum of root spaces corresponding to the simple roots. 

\begin{Def} A nilpotent orbit $\co\subset \cn$ is \emph{Levi-regular} if $\co\cap \fg_2 \neq \emptyset$, or equivalently, if there exists a Levi subalgebra $\fm$ of $\fg$ such that $\co\cap \fm = \co^{\reg}_M$ is the regular nilpotent orbit of $M$ on $\cn_M$.
\end{Def}

 In type $A_n$, every nilpotent orbit is Levi-regular.  The Levi-regular nilpotent orbits in other classical types were described by Russell
  in terms of the partition classification (\cite[Prop.~3.1]{Russell}; see Lemma~\ref{lem.g2orbits} below).  As the next result shows, the Levi-regular
  orbits are closely related to the extended Springer resolution.

\begin{Thm} \label{thm.g2.intersection} 
If the nilpotent orbit $\co$ is Levi-regular, then every simple perverse sheaf $IC(\overline{\co}, \cl)$, for $\cl$ an equivariant local system on $\co$, is a summand of $\psi_*\shc_{\tilde{\cm}}[\dim \cn]$. 
\end{Thm}

Theorem~\ref{t:main1} above implies that $\psi_*\shc_{\tilde{\cm}}[\dim \cn]$  is the direct sum of all the minimal cuspidal Lusztig sheaves.
Therefore, this theorem is equivalent to the statement that if $\co$ is Levi-regular, then each $IC(\overline{\co}, \cl)$ appears in a minimal cuspidal Lusztig sheaf.  

For exceptional groups, Theorem~\ref{thm.g2.intersection} follows from \cite[Section 15]{Lu84}.  Here, Lusztig determines the cuspidal
data in exceptional types.  For each exceptional group, all cuspidal data with $L$ a proper Levi subgroup have $L$ of type $A$, thus
are minimal.  The other cuspidal data all have $L = G$.  Lusztig determines the relevant orbit in each case.  By inspecting the
tables in \cite[Section 8.4]{NilOrbBk}, one can see that all of these orbits are distinguished but not regular, so they are
not Levi-regular.  Hence the theorem holds for exceptional groups.

For classical groups, Theorem~\ref{thm.g2.intersection} is a consequence of a more precise statement, which in type $D_n$ involves a broader
notion called {\em almost Levi-regular} (see Definition \ref{def:almost}).  This notion is specific to type $D_n$.

\begin{Thm} \label{thm.support} 
If $G$ is of classical type, every simple perverse sheaf $IC(\overline{\co}, \cl)$, for $\cl$ an equivariant local system on $\co$, is a summand of $\psi_*\shc_{\tilde{\cm}}[\dim \cn]$ if and only if $\co$ is Levi-regular or almost Levi-regular. 
\end{Thm}

\begin{Rem}
Note that not all nilpotent orbits in the simple Lie algebras of types $E_6$ and $E_7$ are Levi-regular
(see \cite[Tables 2-4]{Russell}). 
\end{Rem}


\subsection{Strategy of the proof} \label{ss.strategy}
If $G$ is of type $A_n$, the theorem follows immediately from the main result of~\cite{GPR23}.  In this case, all nilpotent orbits intersect $\fg_2$ and all simple perverse sheaves 
appear as summands of $\psi_*\shc_{\tilde{\cm}}[\dim \cn]$.  
The remainder of Section \ref{sec.g2.intersection} is devoted to the proof of Theorem~\ref{thm.support} for $G$ of classical type $B_n$, $C_n$, or $D_n$.  In this section, we discuss the strategy of the proof.  An overview of the generalized Springer correspondence with
some details relevant to this section can be found in Appendix \ref{s.overview}.

As discussed in Section \ref{ss.cuspidal}, from a  cuspidal datum $\vc = (L,\co_L, \cl^L)$, Lusztig constructs a perverse sheaf $A_{\vc}$ on $\cn$.
Each such Lusztig sheaf decomposes as a sum of simple perverse sheaves $IC(\overline{\co}, \cl)$.  
Conversely, each simple perverse sheaf $IC(\overline{\co}, \cl)$ on $\cn$ occurs as a summand in exactly one Lusztig sheaf $A_{\vc}$.
We refer to $\vc$ as the cuspidal datum of (or associated to) $IC(\overline{\co}, \cl)$.  
To prove Theorem~\ref{thm.support}, we use Lusztig's combinatorial procedure~\cite{Lu84} to identify the cuspidal datum $\vc$ of $IC(\overline{\co}, \cl)$.  Using this, we can identify the orbits $\co$ for which the associated cuspidal datum is minimal, and thus prove Theorem~\ref{thm.support}.

Our arguments use the fact that to find the cuspidal datum $\vc$ of $IC(\overline{\co}, \cl)$, it is enough to find the Levi subgroup $L$ of $\vc$.
The reason is as follows.  There is a central character (that is, an element of  $\hat{Z}$) associated to any $IC(\overline{\co}, \cl)$ and also to any Lusztig sheaf $A_{\vc}$.  If $A_{\vc}$ contains $IC(\overline{\co}, \cl)$, then these central characters must coincide.  Moreover, any Levi subgroup occurs as part of at most one cuspidal datum with a given central character.
Since the central character associated to $\vc$ is determined by that associated to $IC(\overline{\co}, \cl)$, knowing the Levi subgroup is
enough to specify $\vc$.  It turns out that the cuspidal datum $\vc$ is minimal  if and only if its Levi subgroup $L$ is of type $A$, i.e., all 
simple factors of $[\fl, \fl]$ are of type $A$ (see Theorem \ref{lu_summary}).  
This includes the case where $L = T$, which implies $\vc = (T, \{ 0 \}, \cl_{triv})$,
the central character is trivial, and the Lusztig sheaf is the Springer sheaf.

We now outline the proof of Theorem~\ref{thm.support} for type $C_n$.
The group here is $G = \Sp(2n)$, which is simply connected. Lusztig defines a set $\Psi_{2n}$ of combinatorial objects called symbols, such that
$\Psi_{2n} = \sqcup \Psi_{2n, \gl}$, where $\gl$ runs over the set of partitions parametrizing nilpotent orbits in $\fg$.  Equivariant local systems on the nilpotent orbit $\co_{\gl}$ are in bijection with elements of $\Psi_{2n, \gl}$.  Thus, each symbol has an associated cuspidal datum,
which is minimal if and only if its Levi subgroup is of type $A$. Associated to each symbol is an integer called its defect, which determines the Levi factor of the corresponding cuspidal datum.  The proof of Theorem~\ref{thm.support} in type $C_n$ is thus a consequence of the following two results.
\begin{enumerate}
\item  The cuspidal datum associated to a symbol has Levi subgroup of type $A$ if and only if the defect of the symbol is $\pm 1$ (Corollary \ref{cor.defect} below).
\item The orbit $\co_{\gl}$ is Levi-regular if and only if the defect of every symbol in $\Psi_{\gl}$ is $\pm 1$  (Proposition \ref{prop:typecdefect} below).
\end{enumerate}
Assertion (1) is deduced from Lusztig's determination of the Levi factor from the defect. 
Assertion (2) is proved combinatorially by analyzing the defects of the symbols in $\Psi_{2n, \gl}$ for the partitions
$\gl$ corresponding to Levi-regular orbits.  The details can be found in Section~\ref{sec.Csymbols} below. 

The proofs in types $B_n$ and $D_n$ are similar, but with an additional complication arising from the fact that the orthogonal
groups $\SO(2n+1)$ and $\SO(2n)$ have simply connected double covers $\Spin(2n+1)$ and $\Spin(2n)$.  For the $IC(\overline{\co}, \cl)$ which 
are equivariant with respect to one of the orthogonal groups, our arguments are very similar to those in type $C_n$ and use the combinatorics of symbols; see Section~\ref{sec.BDsymbols}.  However, it remains to consider the $IC(\overline{\co}, \cl)$ which are equivariant for a spin group but not for the corresponding orthogonal group.  Here we apply combinatorics developed by  Lusztig--Spaltenstein~\cite{LuSpal85}. Our arguments are similar in spirit to the previous cases; we refer to Section~\ref{ss.spin} for the details.

\subsection{Levi-regular and almost Levi-regular orbits} \label{ss.g2}
Consider the complex simple Lie algebras of type $B_n$, $C_n$, and $D_n$. In each of these cases, nilpotent orbits are indexed by a particular subset of partitions~\cite[\S 5.1]{NilOrbBk}. 
\begin{itemize}
\item Type $B_n$: Nilpotent orbits in the orthogonal Lie algebra $\mathfrak{so}(2n+1)$ are in bijection with partitions of $2n+1$ in which all even parts occur with even multiplicity.
\item Type $C_n$: Nilpotent orbits in the symplectic Lie algebra $\mathfrak{sp}(2n)$ are in bijection with partitions of $2n$ in which all odd parts occur with even multiplicity.
\item Type $D_n$: Nilpotent orbits in the orthogonal Lie algebra $\mathfrak{so}(2n)$ are in bijection with partitions of $2n$ in which all even parts occur with even multiplicity, except for the \emph{very even} partitions, that is, partitions with only even parts having even multiplicity, which each correspond to two orbits.
\end{itemize}
We denote the set of partitions of $2n$ in which all odd parts occur with even multiplicity by $X_{2n}^C$ and the set of all partitions of $N$ in which all even parts occur with even multiplicity by $X_N^{BD}$.
The next lemma tells us which of these orbits are Levi-regular; see~\cite[Prop.~3.1]{Russell}.

\begin{Lem} \label{lem.g2orbits} Suppose $\lambda$ is a partition indexing a nilpotent orbit in a simple Lie algebra of type $B_n$, $C_n$ or $D_n$.  Then 
$\co_{\lambda}$ is Levi-regular (equivalently, $\co_\lambda \cap \fg_2\neq \emptyset$) if and only if
\begin{itemize}
\item \textup{Type $B_n$:} the partition $\lambda \in X_{2n+1}^{BD}$ has exactly one (necessarily odd) part with odd multiplicity,
\item \textup{Type $C_n$:} the partition $\lambda\in X_{2n}^C$ has at most one (necessarily even) part that occurs with odd multiplicity, or
\item \textup{Type $D_n$:} the partition $\lambda \in X_{2n}^{BD}$ has at most two (necessarily odd) parts with odd multiplicity, and if there are parts with odd multiplicity,
then one of those parts is equal to $1$.
\end{itemize} 
\end{Lem}

\begin{proof}  Suppose $\fg$ is of type $B_n$, $C_n$ or $D_n$. The $T$-orbits on $\fg_2$ correspond to subsets $J$ of $[n]$; an orbit representative
corresponding to $J$ is the sum of simple root vectors for $\ga_i$, where $i\notin J$.  Each $T$-orbit is contained in a single nilpotent orbit (although a single nilpotent
orbit may contain more than one $T$-orbit), and therefore determines a partition $\lambda$ of $2n+1$ (in type $B_n$) or $2n$ (in types $C_n$ and $D_n$).  The Levi-regular orbits correspond to the partitions obtained
in this way.  Russell gives a formula to compute the partition corresponding to each subset $J$ (see \cite[Prop.~3.1]{Russell}), and
by examining this formula, one can see that only the partitions in the statement of the lemma occur.  Conversely, one can reverse the formula to show
that each partition listed in the lemma is obtained from some subset $J$.

For illustration, we will construct such a subset $J$ from a partition $\lambda$ in type $B_n$.  Suppose $\lambda$ is a partition of $2n+1$ with exactly one part $p$ with odd multiplicity. In particular, suppose $\lambda = [p^{m}\ p_1^{\ell_1}\ \dots\ p_k^{\ell_k}]$ for some odd $m$, where each $\ell_i$ is even (we do not assume $p$ is the largest
part of $\lambda$).  Since $p$ is odd, $2n+1-p$ is even, and $\widehat{\lambda} = [p^{m-1}\ p_1^{\ell_1}\ \dots\ p_k^{\ell_k}]$ is a partition of $2n+1-p$ where all parts have even multiplicity.  Thus, we can reduce the multiplicity by half for each part and obtain a partition $\widetilde{\lambda}$ of $d = \frac{1}{2}(2n+1-p)  \leq n$.  Writing $\widetilde{\lambda}$ as $[r_1, \dots r_t]$, we define 
$J = \{ d, d-r_1, d-(r_1+r_2), \dots, d-(r_1+\dots +r_{t-1})=r_t\}$.  A similar process builds a set $J$ in all types.
\end{proof}

We now define the notion of almost Levi-regular in type $D_n$, which is used in the statement of Theorem \ref{thm.support}.
 Recall from \cite[Section 6.1]{NilOrbBk} that a partition is rather odd if all odd parts have multiplicity at most $1$.

 \begin{Def} \label{def:almost}
 In type $D_n$, we will say
 that $\co_{\lambda}$ is {\em almost Levi-regular} if both of the following conditions hold:
 \begin{enumerate} 
 \item $\lambda$ contains at most two (necessarily odd) parts with odd multiplicity.
\item If $\lambda$ is rather odd and has odd parts, then one of its parts is congruent to $1$ mod $4$.
\end{enumerate}
 \end{Def}
 
As noted above, the notion of almost Levi-regular is limited to type $D_n$.  A Levi-regular orbit is also almost Levi-regular.


\subsection{The symplectic group $G=\Sp(2n)$} \label{sec.Csymbols} Assume $G=\Sp(2n)$ and $n\geq 1$.  Any result stated without proof below can be found in Sections 10.4 and 11 of~\cite{Lu84}.

\begin{Def}\label{def.typeCsymbols}
Fix $n\geq 1$. Define $\widetilde \Psi_{2n}$ to be the set of ordered pairs $(A,B)$ where $A$ is a finite subset of $\{0,1,2,\dots\}$ and $B$ is a finite subset of $\{1,2,3,\dots\}$ such that
\begin{enumerate}
\item If $i$ is any integer, then $\{i, i+1\}$ is neither contained in $A$ nor in $B$, 
\item $|A|+|B|$ is odd, and 
\item $\sum_{a\in A} a + \sum_{b\in B} b = n + \frac{1}{2}\left( |A|+|B| \right) \left(|A|+|B| -1 \right)$.
\end{enumerate}
Define $\Psi_{2n}$ to be the set of equivalence classes on $\widetilde\Psi_{2n}$ with respect to the relation that 
\begin{eqnarray}\label{eqn.symbol.rel}
(A,B) \sim (\{0\}\cup (A+2), \{1\}\cup(B+2)).
\end{eqnarray}
The elements of $\Psi_{2n}$ are called \emph{type C symbols}.
\end{Def}

\begin{example} The pair $(\{0,2,5\}, \{1,4\})$ clearly satisfies conditions (1) and (2) of Definition~\ref{def.typeCsymbols}. It is an element of $\Psi_{4}$ since
\[
\sum_{a\in A} a +\sum_{b\in B} b = (0+2+5)+(1+4) = 12
\]
while 
\[
\frac{1}{2}\left( |A|+|B| \right) \left(|A|+|B| -1 \right) = \frac{1}{2}(5)(4) = 10
\]
so $2=n = 12 - 10$.  We have $(\{0,2,5\},\{1,4\}) \sim ( \{0,3\}, \{2\} )  \sim (\{0,2,4,7\}, \{1,3,6\})$.
For another example, the reader can find a list of all symbols in the set $\Psi_{4}$ in the last column of the second displayed table in Example~\ref{ex.C2} below. 
\end{example}

Two elements of $\Psi_{2n}$ are \emph{similar} if they can be represented by pairs $(A,B)$ and $(A',B')$ such that
$A\cup B=A'\cup B'$ and $A\cap B=A'\cap B'$. We call the set of all elements of $\Psi_{2n}$ similar to $(A,B)$ the \emph{similarity class} of $(A,B)$. Within each similarity class, there is a unique element 
$(A,B)$ such that 
\begin{equation} \label{e:distinguished}
A=\{a_1<a_2< \cdots<a_m<a_{m+1}\} \; \text{ and } \; B=\{b_1<b_2<\cdots < b_m\}
\end{equation}
for some $m$ with $a_1\leq b_1\leq a_2\leq b_2\leq \cdots \leq a_m \leq b_m \leq a_{m+1}$.  
We call this element \emph{distinguished}.  Note that the condition \eqref{e:distinguished}
is preserved by the equivalence relation used to define $\Psi_{2n}$, so it can be verified using any
representative $(A,B)$.

Following Lusztig~\cite[\S 11.6]{Lu84}, we define a distinguished element $(A_\lambda, B_\lambda)$ of $\Psi_{2n}$, for any partition $\lambda\in X_{2n}^C$.  

\begin{Def} Let $\lambda\in X_{2n}^C$ and write $2n = r_1+2r_2+\cdots + \ell r_\ell$, where $r_k$ is the number of parts of $\lambda$ equal to $k$.
\begin{enumerate}
\item Consider the ordered list
\[
(1,\ldots, 1, 2, \ldots, 2, \ldots, \ell, \ldots, \ell),
\]
where $1$ appears $r_1$ times, $2$ appears $r_2$ times, and so on.  
\item By appending $0$ if need be, we may assume the ordered list has an even number of terms.  Add $i-1$ to the $i$-th term.
\item Separate the resulting list from (2) into its even and odd parts, denoted by 
\[
\vec{e} = (2y_1<2y_2<\cdots < 2y_m) \text{ and } \vec{o} = (2y_1'+1< 2y_2'+1<\cdots < 2y_m'+1).
\]
\end{enumerate}
Finally, define 
\begin{eqnarray*}
A_\lambda &:=& \{0, y_1'+2, y_2'+3, \ldots, y_m'+m+1 \}, \text{ and } \\
 B_\lambda &:=& \{y_1+1, y_2+2, \ldots, y_m+m\}.
\end{eqnarray*}
The \emph{distinguished symbol corresponding to the partition $\lambda \in X_{2n}^C$} is $(A_\lambda, B_\lambda)$. We denote the similarity class in $\Psi_{2n}$ of $(A_\lambda, B_{\lambda})$ by $\Psi_{2n, \lambda}$.
\end{Def}

\begin{example}\label{ex.dist} Let $\lambda=[4\,2^3]\in X_{10}^C$. We begin with the ordered list $(2,2,2,4)$ which has even length. Adding $i-1$ to the $i$-th term we obtain $(2,3,4,7)$ which separates into
\[
\vec{e} = (2,4) = (2\cdot 1, 2\cdot 2) \;\; \text{ and } \;\; \vec{o} = (3,7)=(2\cdot 1 +1, 2\cdot 3+1)
\]
so $(A_{[4\,2^3]}, B_{[4\,2^3]}) = (\{0,3,6\},\{2,4\}) \in \Psi_{10}$. Since $0\leq 2\leq 3\leq 4\leq 6$, this symbol is distinguished, as expected.
\end{example}

Consider the set $C_\lambda = (A_\lambda \cup B_\lambda)-(A_\lambda \cap B_\lambda)$.  A nonempty subset $I$ of $C_\lambda$ is an \emph{interval} if it is of the form $\{i,i+1,\dots, j\}$ for some $i<j$ with $i-1, j+1\notin C_\lambda$. Let $I_0$ denote the interval containing $0$ in $C_\lambda$ and define $\mathcal I_\lambda$ to be the collection of all intervals in $C_\lambda$, excluding $I_0$.  

\begin{Rem}\label{rem.zero.interval} By definition, $0\in A_\lambda$ and the elements of the interval $I_0$ must alternate between the sets $A_\lambda$ and $B_\lambda$.  It follows that $|A_\lambda \cap I_0| = |B_\lambda\cap I_0|+1$ if $I_0$ is an odd-length interval and $|A_\lambda \cap I_0| = |B_\lambda\cap I_0|$ if $I_0$ is an even-length interval. 
\end{Rem}

For each subset $S$ of intervals in  $\mathcal{I}_\lambda$, we obtain a new ordered pair $(A_{\lambda, S}, B_{\lambda, S})$ in $\Psi_{2n}$ in the similarity class $\Psi_{2n, \lambda}$ of $(A_\lambda, B_\lambda)$.   In particular, for each subset $S$ of $\ci_\lambda$, we define
\begin{align}
\label{eqn.int.bijection.def}
A_{\lambda,S} &:= \left(\bigcup_{I\in S}(I\cap A)\right)\cup\left(\bigcup_{I\notin S} (I\cap B)\right)\cup(I_0\cap A)\cup(A\cap B),\, \text{ and }\\
\nonumber B_{\lambda,S} &:= \left(\bigcup_{I\in S}(I\cap B)\right)\cup\left(\bigcup_{I\notin S} (I\cap A)\right)\cup(I_0\cap B)\cup(A\cap B).
\end{align}
That is, we fix the values $I_0=\{0,1,2,\dots\}$ of $A_\lambda$ and $B_\lambda$ and fix all elements in the intervals $I$ in the set $S$.  Then, we swap the remaining elements between $A_\lambda$ and $B_\lambda$.  We will see below that every element of $\Psi_{2n}$ is of the form $(A_{\lambda, S}, B_{\lambda, S})$ for some $S\subseteq \mathcal{I}_\lambda$ and $\lambda\in X_{2n}^C$.

\begin{example}\label{ex.4222} Let $\lambda=[4\,2^3]$ as in Example~\ref{ex.dist}, where it was computed that $(A_{[4\,2^3]}, B_{[4\,2^3]}) = (\{0,3,6\},\{2,4\})$. Thus $C_{[4\,2^3]} = \{0,2,3,4,6\}$ and $\mathcal I_{[4\,2^3]} = \{ \{2,3,4\},\{6\} \}$. The table below displays the correspondence from~\eqref{eqn.int.bijection.def} between subsets of $\mathcal I_{[4\, 2^3]}$ and elements of the similarity class $\Psi_{10, [4\,2^3]}$.
\[
\begin{array}{c|c}
S\subseteq \mathcal I_{[4\, 2^3]} & (A_{[4\,2^3], S}, B_{[4\,2^3], S})\\ \hline
 \mathcal I_{[4\, 2^3]} &  (A_{[4\,2^3]}, B_{[4\,2^3]}) = (\{0,3,6\},\{2,4\}) \\
 \{\{6\} \} & (\{0,2,4,6\},\{3\}) \\
 \{\{2,3,4\}\} &  (\{0,3\},\{2,4,6\}) \\ 
 \emptyset & (\{0,2,4\},\{3,6\})
\end{array}
\]
\end{example}

\begin{example}\label{ex.C2} Let $n=2$. The following table displays each partition $\lambda$ in $X_4^C$, the corresponding symbol $(A_\lambda, B_\lambda)$ in $\Psi_4$, the set $C_\lambda = (A_\gl \cup B_\gl) - (A_\gl\cap B_\gl)$, and corresponding set of intervals.
\[
\begin{array}{c|c|c|c}
\lambda\in X_4^C & (A_\lambda, B_\lambda) & C_\lambda & \mathcal I_\lambda\\ \hline
[4] & (\{0,4\}, \{1\}) & \{0,1,4\} & \{\{4\}\} \\ \hline
[2^2] & (\{0,3\}, \{2\}) & \{0,2,3\} & \{\{2,3\}\} \\ \hline
[2\,1^2] & (\{0,3,5\},\{1,3\}) & \{0,1, 5\} & \{\{5\}\} \\ \hline
[1^4] &  (\{0,2,4\},\{2,4\}) & \{0\} & \emptyset
\end{array}
\]
Now we compute all elements of $\Psi_4$ below (using the fact that each is of the form $(A_{\lambda, S}, B_{\lambda, S})$; see Lemma~\ref{lem.int.bijection} below)
\[
\begin{array}{c|c|c}
\lambda\in X_4^C & S\subseteq \mathcal I_\lambda & (A_{\lambda, S}, B_{\lambda, S}) \\ \hline
[4] &  \{\{4\}\} & (\{0,4\}, \{1\}) \\  
 & \emptyset & (\{0\}, \{1,4\}) \\ \hline
[2^2] &  \{\{2,3\}\} &  (\{0,3\}, \{2\}) \\
 & \emptyset & (\{0,2\}, \{3\}) \\ \hline
 [2\,1^2] &  \{\{5\}\} & (\{0,3,5\},\{1,3\}) \\
 & \emptyset & (\{0,3\},\{1,3,5\}) \\ \hline
[1^4] & \emptyset  &  (\{0,2,4\},\{2,4\})
\end{array}
\]
\end{example}

The key facts about the set $\Psi_{2n}$, for our purposes, are summarized in the following lemma.  
Recall that $\Psi_{2n,\lambda}$ is the similarity class of $(A_\lambda, B_\lambda)$.  

\begin{Lem}[Lusztig]\label{lem.int.bijection} 
We have $\Psi_{2n} =  \bigsqcup_{\lambda\in X_{2n}^C} \Psi_{2n, \lambda}$.
Let $\lambda\in X_{2n}^C$, and denote the power set of $\ci_\lambda$ by $\mathscr{P}(I_\lambda)$. There are bijections
\begin{eqnarray}\label{eqn.int.bijection}
\{ (\co_\lambda, \cl) \mid \text{$\cl$ a $G$-equiv.~local system on $\co_\lambda$} \} \leftrightarrow \mathscr{P}(\ci_\lambda)  \leftrightarrow \Psi_{2n,\lambda} 
\end{eqnarray}
where the bijection $\mathscr{P}(\ci_\lambda)  \rightarrow \Psi_{2n,\lambda}$ is given by \eqref{eqn.int.bijection.def}. Furthermore, the elements of $\ci_\lambda$ are in bijection with the distinct even parts of $\lambda$, and if we associate the $j$-th interval $I_j$ to the $j$-th smallest distinct even part of $\lambda$, then the size of the interval is equal to the multiplicity of the corresponding even part. 
\end{Lem}

Consider $\lambda = [4\,2^3]$ as in Example~\ref{ex.4222}. The partition $\lambda$ has two distinct even parts: $2$ and $4$. The number $2$ occurs with multiplicity $3$, and 
the number $4$ occurs with multiplicity $1$. Since $\ci_\lambda = \{  \{2,3,4\},\{6\}\}$, this confirms the last part of the statement of Lemma~\ref{lem.int.bijection} in this example.

Recall that given a pair $(\co_\lambda, \cl)$, where $\cl$ is a $G$-equivariant local system on $\co_\lambda$, the corresponding simple perverse sheaf $IC(\overline{\co_\lambda}, \cl)$ appears as a summand in exactly one Lusztig sheaf $\A_\mathbf{c}$ indexed by the cuspidal datum $\mathbf{c}$.  We say that $\mathbf{c}$ is the cuspidal datum corresponding to $(\co_\lambda, \cl)$.  In light of the bijections appearing in~\eqref{eqn.int.bijection} of Lemma~\ref{lem.int.bijection}, we will also refer to $\mathbf{c}$ as the \emph{cuspidal datum corresponding to the symbol equal to the image of $(\co_\lambda, \cl)$} under~\eqref{eqn.int.bijection}.

We need one final definition before we move on to the application of all of this to our setting.

\begin{Def} The \emph{defect} of a symbol $(A,B) \in \Psi_{2n}$ is $|A|-|B|$.
\end{Def}

The defect of the distinguished symbol $(A_\lambda, B_\lambda)$ for $\lambda\in X_{2n}^C$ is $1$. Note that the defect of  any symbol is odd.  Indeed,  if $(A,B)\in \Psi_{2n}$ is similar to $(A',B')$, then their defects are the same mod 2, as follows from the definition of similarity and the congruence $|A| + |B| \equiv |A \cup B | + |A \cap B|$ mod 2.  Since $(A,B)$ is similar to a distinguished element, its defect is necessarily odd.

The defect of the symbol $(A_{\lambda, S}, B_{\lambda, S})$ determines the Levi subgroup that appears in its corresponding cuspidal datum. Given an integer $i \geq 0$, let  $L_i\subseteq G=\Sp(2n)$ be the standard Levi subgroup of type
\[
L_i \simeq \underbrace{GL_1\times \cdots  \times GL_1}_{\text{$i$ factors}} \times \Sp(2n-2i).
\]
The relative Weyl group $N_G(L_i)/L_i$ of $L_i$ is a type $C_i$ Weyl group. 

\begin{Thm}[Lusztig] \label{thm.Lusztig.TypeC}  Suppose the defect of $(A,B) \in \Psi_{2n}$ is $k$ (which is necessarily odd).  Then the cuspidal datum corresponding to $(A,B)$ has Levi factor $L_i$, where $i = n - \frac{k(k-1)}{2}$.  Moreover,  any Levi subgroup of $G$ that supports cuspidal perverse sheaves must be  conjugate to some $L_i$, where $i = n - \frac{k(k-1)}{2}$ for some odd integer $k$.  
\end{Thm}

\begin{example}\label{ex.defect1} Let $n=2$. Using Example~\ref{ex.C2} we see that there are 5 symbols with defect $1$ in $\Psi_4$ and 2 symbols with defect $-1$. 
By Theorem~\ref{thm.Lusztig.TypeC}, we see that the symbols with defect $1$ correspond to cuspidal data with Levi factor $L_2 \simeq GL_1\times GL_1 = T$, and the  two symbols of defect $-1$ have cuspidal data with Levi factor $L_1\simeq GL_1 \times \Sp(2)$.
\end{example}

\begin{example}\label{ex.defect2} Let $n=5$. Consider the symbols corresponding to the partition $[4\,2^3]$ computed as in Example~\ref{ex.4222}. The table below computes the defect of each symbol, and the Levi subgroup for the corresponding cuspidal datum. 
\[
\begin{array}{c|c|c}
\text{ Symbol } & \text{ Defect } & L\\ \hline
(\{0,3,6\},\{2,4\}) & 1 & T = GL_1\times GL_1 \times GL_1\times GL_1\times GL_1  \\
 (\{0,2,4,6\},\{3\}) & 3 & GL_1\times GL_1\times \Sp(6) \\
(\{0,3\},\{2,4,6\}) & -1 &  GL_1 \times GL_1\times GL_1\times GL_1\times \Sp(2) \\ 
(\{0,2,4\},\{3,6\}) & 1 &  T = GL_1\times GL_1 \times GL_1\times GL_1\times GL_1 
\end{array}
\]
\end{example}

\begin{Cor}\label{cor.defect} Let $\co$ be a nilpotent $G$-orbit and $\cl$ a $G$-equivariant local system on $\co$.  Let $(A,B)$ denote the symbol equal to the image of the pair $(\co, \cl)$ under the bijections appearing in~\eqref{eqn.int.bijection}.  Then the cuspidal datum $\mathbf{c}$ corresponding to $(A,B)$ is 
\begin{itemize}
\item $\mathbf{c}=(T, \{0\}, \cl_{triv})$ if and only if the defect of $(A,B)$ is $1$ (that is, $IC(\overline{\co}, \cl)$ appears in the Springer sheaf), and 
\item $\mathbf{c} = (L_\chi, \co_{\reg}, \cl^{L_\chi})$ where $L_\chi$ is the Levi subgroup of type $A_1$ corresponding to the character $\chi$ appearing in the type $C_n$ row of Table~\ref{fig.Levi.table} if and only if the defect of $(A,B)$ is $-1$.  
\end{itemize}
In particular, $IC(\overline{\co}, \cl)$ appears as a summand of a Lusztig sheaf in $\psi_*\shc_{\tilde{\cm}}[\dim \cn]$ if and only if the defect of the corresponding symbol is $\pm 1$.
\end{Cor}
\begin{proof} By Theorem~\ref{thm.Lusztig.TypeC}, a symbol with defect $1$ corresponds to a cuspidal datum $\mathbf{c}$ with Levi factor $L_n = T$, and a symbol with defect $-1$ corresponds to
\[
L_{n-1}\simeq \underbrace{GL_1\times \cdots \times GL_1}_{\text{$n-1$ factors}}\times  \Sp(2) = \underbrace{GL_1\times \cdots \times GL_1}_{\text{$n-1$ factors}}\times  A_1. 
\]
When the defect of $(A,B) \in \Psi_{2n}$ is  $k\neq \pm 1$, it is easy to confirm that the corresponding Levi subgroup $L_i$ as in Theorem~\ref{thm.Lusztig.TypeC} satisfies $L_i\not\simeq T$ and $L_i$ is not of type $A$. Since $\mathbf{c}$ is minimal if and only if the Levi factor of $\mathbf{c}$ is of type $A$, this completes the proof. 
\end{proof}

\begin{Prop}\label{prop:typecdefect}
Suppose $\lambda\in X_{2n}^C$.  The orbit  $\co_\lambda$ is Levi-regular if and only if the defect of every symbol in the similarity class of $(A_\lambda, B_\lambda)$ is $\pm 1$. 
\end{Prop}

\begin{proof} Let $\lambda\in X_{2n}^C$. By Lemma~\ref{lem.g2orbits}, $\co_{\lambda}\cap \fg_2\neq \emptyset$ if and only if $\lambda$ has at most one even part that occurs with odd multiplicity.  A symbol has defect $1$ or $-1$ exactly when the size of $A$ and the size of $B$ only differ by $1$.  

As noted above, by the construction of $A_\lambda$ and $B_\lambda$, the defect of $(A_\lambda, B_\lambda) $ is $1$.   Also, by definition of the symbol $(A,B)$, consecutive integers must alternate between the sets $A$ and $B$.  Hence, for any interval $I$ in $A\cup B$, if $|I|$ is even, then
$|I \cap A| = |I \cap B|$; if $|I|$ is odd, then $|I \cap A| = |I \cap B| \pm 1$.

Any symbol in the similarity class of $(A_\lambda, B_\lambda)$ is of the form $(A_{\lambda, S}, B_{\lambda, S})$ for some $S\in \mathscr{P}(\ci_\lambda)$.
The symbol $(A_{\lambda, S}, B_{\lambda, S})$ is obtained from $A_\lambda$ and $B_\lambda$ by swapping elements belonging to intervals $I$ such that $I\notin S$ between the two sets.  The paragraph above implies that swapping elements belonging to an even-length interval does not change the defect.  

By Lemma~\ref{lem.int.bijection}, odd-length intervals in $\mathscr{P}(\ci_\lambda)$ come from even parts of $\lambda$ which occur with odd multiplicities.  Hence, if no even part of $\lambda$ occurs with odd multiplicity, then there are no odd-length intervals $I$ in $\ci_{\lambda}$, so for all $S\in \mathscr{P}(\ci_\lambda)$, the defect of $(A_{\lambda, S}, B_{\lambda, S})$ is $1$.  If one even part of $\lambda$ occurs with odd multiplicity,  then there is at most one odd length interval $I$.   Then the defect of $(A_{\lambda,S}, B_{\lambda,S})$ will be $1$ if $I\in S$ and $-1$ if $I\notin S$. 

Suppose now that $\lambda$ has at least two even parts with odd multiplicity.  We claim that there exists an odd-length interval $I\in \ci_\lambda$ such that $|I \cap A_\lambda|+1 = |I \cap B_\lambda|$.  Indeed,
$$
1 = |A_{\gl} | - |B_{\gl} | = |A_\lambda\cap I_0| - |B_\lambda\cap I_0| + \sum_{ \substack{I \in \mathcal{I}_\lambda \\  | I | \textrm{odd}}}  ( |A_\lambda\cap I| - |B_\lambda\cap I| ).
$$
The difference $|A_\lambda\cap I_0| - |B_\lambda\cap I_0|$ is $0$ or $1$, our assumption on $\lambda$ implies there are at least two terms in the sum,  and each of the differences $|A_\lambda\cap I| - |B_\lambda\cap I|$ is $\pm 1$.  Hence
at least one of these differences must be $-1$, proving the claim.

Now consider $S = \mathcal I_\lambda - \{I\}$ where $I\in \ci_\lambda$ has odd length and $|I\cap A_\lambda|+1 = |I\cap B_\lambda|$.  By definition, 
\[
A_{\lambda, S} = (A_\lambda \cap I^c) \cup (B_\lambda \cap I) \text{ and } B_{\lambda, S} = (B_\lambda \cap I^c) \cup (A_\lambda \cap I)
\] 
so $|A_{\lambda,S} | = |A_{\lambda} |+1$ and $|B_{\lambda,S} | = |B_{\lambda} |-1$.  Thus, the defect of $(A_{\lambda,S},B_{\lambda,S})$ is equal to $3$.  Since this phenomenon only requires the existence of at least two even parts with odd multiplicity in $\lambda$, the result follows.
\end{proof}

\begin{proof}[Proof of Theorem~\ref{thm.support} in type $C_n$]
Let $\lambda\in X_{2n}^C$.   By Lemma~\ref{lem.int.bijection} the subset of $\Psi_{2n}$ corresponding to pairs $(\overline{\co}_\lambda, \cl)$ is precisely the similarity class $\Psi_{2n,\lambda}$ of $(A_\lambda, B_\lambda)$. 
By Proposition~\ref{prop:typecdefect}, $IC(\overline{\co}_\lambda, \cl)$ is a summand of $\psi_*\shc_{\tilde{\cm}}[\dim \cn]$ if and only if the defect of the symbol corresponding to the pair $(\co_\lambda, \cl)$ is $\pm1$.  Thus every simple perverse sheaf of the form $IC(\overline{\co}_\lambda, \cl)$ appears as a summand of $\psi_*\shc_{\tilde{\cm}}[\dim \cn]$ if and only if every symbol in $\Psi_{2n,\lambda}$ has defect $\pm 1$.  The result now follows immediately from Proposition~\ref{prop:typecdefect}.
\end{proof}


\subsection{The orthogonal groups $\SO(2n)$ and $\SO(2n+1)$} \label{sec.BDsymbols}

We begin with Lusztig's symbol constructions describing the generalized Springer correspondence for the special orthogonal group $G=\SO(N)$.

\begin{Def}
Fix $N\geq 3$. Define $\widetilde \Psi_{N}'$ to be the set of unordered pairs $\{A,B\}$ where $A$ and $B$ are finite subsets of $\{0,1,2,\dots\}$ such that condition (1) of Definition~\ref{def.typeCsymbols} holds (i.e., consecutive integers alternate as members of the sets $A$ and $B$) and so that 
\begin{eqnarray}\label{eqn.typeBD}
\sum_{a\in A} a + \sum_{b\in B} b = \frac{1}{2}N + \frac{1}{2} \left( \left(|A|+|B| -1 \right)^2 -1\right).
\end{eqnarray}
This equation implies that  $|A|+|B|$ has the same parity as $N$, since both sides of the equation are integers.  
Let $\Psi_{N}'$ be the set of equivalence classes on $\widetilde\Psi_{N}'$ with respect to the relation
\begin{eqnarray}\label{eqn.symbol.relBD}
\{ A,B \} \sim \{ \{0\}\cup (A+2), \{0\}\cup(B+2) \}.
\end{eqnarray}
The elements of $\Psi_{N}'$ are called \emph{type B/D symbols}.
\end{Def}

As in the section above, we say that the symbols $\{A,B\}$ and $\{A',B'\}$ in $\Psi_N'$ are \emph{similar} if $A\cup B = A'\cup B'$ and $A\cap B = A'\cap B'$. We call the set of all elements of $\Psi_N'$ similar to $\{A,B\}$ the \emph{similarity class} of $\{A,B\}$.  As before, each similarity class contains a unique \emph{distinguished element}, that is, a unique element $\{A,B\}$ such that 
\[
A=\{a_1<a_2<\cdots < a_{m'}\} \;\; \text{ and } \;\; B = \{b_1<b_2<\cdots <b_m\}
\]
and 
\begin{itemize}
\item If $N$ is odd, $m'=m+1$ and $a_1\leq b_1\leq a_2\leq b_2 \cdots \leq a_m \leq b_m \leq a_{m+1}$, and 
\item if $N$ is even $m'=m$ and $a_1\leq b_1\leq a_2\leq b_2 \cdots \leq a_m \leq b_m$.
\end{itemize}

Following Lusztig~\cite[\S 11.7]{Lu84}, we define a distinguished element $\{A_\lambda, B_\lambda\}$ of $\Psi_{N}'$ for any partition $\lambda\in X_{N}^{BD}$.  
\begin{Def} Let $\lambda\in X_{N}^{BD}$ and write $N = r_1+2r_2+\cdots + \ell r_\ell$, where $r_k$ is the number of parts of $\lambda$ equal to $k$. 
\begin{enumerate}
\item Consider the ordered list
\[
(1,\ldots, 1, 2, \ldots, 2, \ldots, \ell, \ldots, \ell),
\]
where $1$ appears $r_1$ times, $2$ appears $r_2$ times, and so on. The length $M$ of this list has the same parity as $N$. 
\item Add $i-1$ to the $i$-th term.
\item ($N$ even) Exactly $M/2$ of the terms in the list resulting from (2) are even and $M/2$ terms are odd. Separate the list into its even and odd parts, denoted by 
\begin{eqnarray*}
\vec{e} &=& (2y_1<2y_2<\cdots < 2y_{M/2})\\
\vec{o} &=& (2y_1'+1< 2y_2'+1<\cdots < 2y_{M/2}'+1).
\end{eqnarray*}
\item ($N$ odd) Exactly $(M-1)/2$ of the terms in the list resulting from (2) are even and $(M+1)/2$ terms are odd. Separate the list into its even and odd parts, denoted by 
\begin{eqnarray*}
\vec{e} &=& (2y_1<2y_2<\cdots < 2y_{(M-1)/2})\\ 
\vec{o} &=& (2y_1'+1< 2y_2'+1<\cdots < 2y_{(M+1)/2}'+1).
\end{eqnarray*}
\end{enumerate}
Finally, if $N$ is even set 
\begin{eqnarray*}
A_\lambda &:=& \left\{ y_1', y_2'+1, \ldots, y_{M/2}'+\tfrac{M}{2}-1 \right\}, \text{ and } \\
 B_\lambda &:=& \left\{y_1, y_2+1, \ldots, y_{M/2}+\tfrac{M}{2} -1 \right\}
\end{eqnarray*}
and if $N$ is odd set
\begin{eqnarray*}
A_\lambda &:=& \left\{ y_1', y_2'+1, \ldots, y_{(M+1)/2}'+\tfrac{M+1}{2}-1 \right\}, \text{ and } \\
 B_\lambda &:=& \left\{y_1, y_2+1, \ldots, y_{(M-1)/2}+\tfrac{M-1}{2}-1 \right\}.
\end{eqnarray*}
The \emph{distinguished symbol corresponding to the partition $\lambda\in X_N^{BD}$} is $\{A_\lambda, B_\lambda\}$. We denote the similarity class in $\Psi_{N}'$ of $\{A_\lambda, B_{\lambda}\}$ by $\Psi_{N,\lambda}'$.
\end{Def}

\begin{example} Let $N=7$ and consider $\lambda =[5\,1^2]\in X_7^{BD}$. Following steps (1), (2) and (4) in the definition above,
\begin{eqnarray*}
(1,1,5) &\longrightarrow& (1,2,7)\\ &\longrightarrow& \vec{e} = (2) = (2\cdot 1) \, \text{ and } \, \vec{o} = (1,7) = (2\cdot 0+1, 2\cdot 3+1).
\end{eqnarray*}
Thus $\{A_{\lambda}, B_{\lambda}\} = \{\{0,4\}, \{1\}\}$. Now consider $\mu=[3\,1^4]\in X_{7}^{BD}$. Once again, following the steps we get
\begin{eqnarray*}
(1,1,1,1,3) &\longrightarrow& (1,2,3,4,7)\\ &\longrightarrow& \vec{e}=(2,4) = (2\cdot 1, 2\cdot 2) \, \text{ and } \, \\
&&\qquad\qquad \vec{o} = (1,3,7) = (2\cdot 0+1, 2\cdot 1 + 1, 2\cdot 3 +1)
\end{eqnarray*}
so $\{A_\lambda, B_\lambda\} = \{\{ 0,2,5 \},\{1,3\}\}$. 
\end{example}

\begin{example}\label{ex.VE} Let $N=8$ and consider $\lambda = [2^4]\in X_8^{BD}$.  Recall that $\lambda$ indexes two nilpotent orbits, as this is an example of a very even partition. Following the steps to compute the associated symbol we get
\begin{eqnarray*}
(2,2,2,2)&\longrightarrow& (2,3,4,5)\\ &\longrightarrow& \vec{e} = (2,4) = (2\cdot 1, 2\cdot 2) \, \text{ and } \, \vec{o} = (3,5) = (2\cdot 1+1, 2\cdot 2+1)
\end{eqnarray*}
so $\{A_\lambda, B_\lambda\} = \{\{1,3\},\{1,3\}\}$. 
\end{example}

For each $\lambda \in X_N^{BD}$ set $C_\lambda = (A_\lambda \cup B_\lambda)- (A_\lambda\cap B_{\lambda})$. As in the previous section, a nonempty subset $I$ of $C_\lambda$ is an 
\emph{interval} if it is of the form $\{i, i+1, \ldots, j\}$ for some $i<j$ with $i-1, j+1 \notin C_\lambda$. Let $\ci_\lambda$ be the collection of all intervals in $C_\lambda$. Unlike the type $C_n$ case considered in the previous section, we do not exclude the interval containing~$0$.

For each subset $S$ of intervals in $\ci_\lambda$ the sets $A_{\lambda, S}$ and $B_{\lambda, S}$ defined as in~\eqref{eqn.int.bijection.def} define a symbol $\{A_{\lambda, S},B_{\lambda, S}\}\in \Psi_{N,\lambda}'$. One difference in this case is that symbols are defined by unordered pairs of sets, so a subset $S\subseteq \ci_\lambda$ and its complement $\ci_\lambda- S$ yield the same symbol.

\begin{Rem} Recall that $\lambda\in X_N^{BD}$ is called very even if it has only even parts that each occur with even multiplicity. It is easy to confirm that $\ci_\lambda = \emptyset$, as $A_\lambda = B_\lambda$. For example, if $N=8$ and $\lambda = [2^4]$ then the computation in Example~\ref{ex.VE} shows that $A_\lambda = B_\lambda = \{1,3\}$.  Lemma~\ref{lem.int.bijectionBD} below tells us that every element of the similarity class $\Psi_{N,\lambda}'$ is of the form $\{A_{\lambda, S}, B_{\lambda, S}\}$. Thus, if $\lambda$ is very even, the similarity class $\Psi_{N,\lambda}'$ consists of a single element, namely $\{A_\lambda, B_\lambda\}$ itself. 
\end{Rem}

The following is a type $B/D$ variant of Lemma~\ref{lem.int.bijection}.

\begin{Lem}[Lusztig] \label{lem.int.bijectionBD} 
We have $\Psi_N' = \sqcup_{\lambda\in X_N^{BD}} \Psi_{N,\lambda}'$. Let $\lambda\in X_{N}^{BD}$ and $\mathscr{P}(\ci_\lambda)'$ denote the power set of $\ci_\lambda$, modulo the relation that identifies each set with its complement. Then there is a bijection,
\begin{eqnarray}\label{eqn.int.bijectionBD}
\mathscr{P}(\ci_\lambda)'  \leftrightarrow \Psi_{N,\lambda} '
\end{eqnarray}
defined by $S\mapsto \{A_{\lambda, S}, B_{\lambda, S}\}$. Furthermore, there exists a map
\begin{eqnarray}\label{eqn.int.bijectionBD2}
\{ (\co_\lambda, \cl) \mid \text{$\cl$ a $G$-equiv.~local system on $\co_\lambda$} \} \rightarrow \Psi_{N,\lambda} '
\end{eqnarray}
that is a bijection when $\lambda$ is not very even and, when $\lambda$ is very even, such that the fiber over $\{A_\lambda, B_\lambda\}$ has two elements (the two $IC$-sheaves corresponding to the two orbits indexed by $\lambda$ defined using the trivial local system). Finally, the elements of $\ci_\lambda$ are in bijection with the set of distinct odd parts of $\lambda$, and if we associate the  $j$-th interval $I_j$ in $\ci_\lambda$ to the $j$-th smallest distinct odd part of $\lambda$,
then the size of the interval is equal to the multiplicity of the corresponding odd part. 
\end{Lem}

We refer to the cuspidal datum $\mathbf{c}$ such that $IC(\overline{\co}_\lambda, \cl)$ appears as a summand of the Lusztig sheaf $\A_{\mathbf{c}}$ as the \emph{cuspidal datum corresponding to the symbol equal to the image of $(\co_\lambda, \cl)$} under~\eqref{eqn.int.bijectionBD2}.

\begin{example}\label{ex.53221111} Consider $\lambda = [5\,3\,2^2\,1^4]\in X_{16}^{BD}$. To compute the corresponding symbol, we have
\begin{eqnarray*}
(1,1,1,1,2,2,3,5) &\longrightarrow& (1,2,3,4,6,7,9,12) \\
&\longrightarrow& \vec{e}=(2,4,6,12) \, \text{ and } \, \vec{o} = (1,3,7,9)
\end{eqnarray*}
so $\{A_\lambda, B_\lambda\} = \{\{0,2,5,7\}, \{1,3,5,9\}\}$. In this case $C_\lambda = \{0,1,2,3,7,9\}$ and $\ci_\lambda = \{\{0,1,2,3\},\{7\},\{9\}\}$.  Taking $S=\{\{0,1,2,3\}\}$ gives us 
\[
\{A_{\lambda,S}, B_{\lambda,S}\} = \{ \{0,2,5,9\},\{1,3,5,7\} \}.
\] 
If $S=\{\{7\}\}$ we get 
\[
\{A_{\lambda,S}, B_{\lambda,S}\} = \{\{1,3,5,7,9\}, \{0,2,5\}\}.
\]
Furthermore, $\lambda$ has exactly three distinct odd parts: $1, 3, 5$. The number $1$ appears with multiplicity $4$, which confirms the last part of the statement of Lemma~\ref{lem.int.bijectionBD} since the first interval in $\ci_\lambda$ has length $4$. 
\end{example}

\begin{Def} The \emph{defect} of a symbol $\{A, B\} \in \Psi_N'$ is $||A|-|B||$.
\end{Def}

The defect of the distinguished symbol $\{A_\lambda, B_\lambda\}$ for $\lambda\in X_N^{BD}$ is $1$ when $N$ is odd and $0$ when $N$ is even. In general, the defect of any type B/D symbol is always a nonnegative integer with the same parity as $N$. We denote finite Weyl groups of type $B_n$ and $D_n$ by $W(B_n)$ and $W(D_n)$, respectively (of course $W(B_n) \cong W(C_n)$).
Let $i=\frac{1}{2}(N-j^2)$ for some $j\geq 0$. Let  $L_i$ be the Levi subgroup of $G=\SO(N)$ of type
\begin{equation} \label{eqn:Li}
L_i \simeq \underbrace{GL_1\times \cdots \times GL_1}_{\text{$i$ factors}} \times \SO(N-2i)
\end{equation}
with relative Weyl group 
\[
N_G(L_i)/L_i \simeq \left\{\begin{array}{ll} W(D_{N/2}) & \text{ if $j^2=0$} \\ 
W(B_{(N-1)/2}) & \text{ if $j^2=1$} \\
W(B_i) & \text{ if $j^2\geq 4$}. \end{array} \right. 
\]

\begin{Thm}[Lusztig]\label{thm.Lusztig.TypeBD}
Suppose that the defect of $\{A,B\}\in \Psi_{N}'$ is $j$.  Then the cuspidal datum corresponding to $\{A,B\}$ has Levi factor $L_i$, where $i = \frac{1}{2}(N - j^2)$.  Moreover, any Levi subgroup of $G$ that supports cuspidal perverse sheaves must be conjugate
to $L_i$ for some $i = \frac{1}{2}(N - j^2)$.
\end{Thm}

\begin{example} Let $N=16$.  Consider the symbols corresponding to the partition $\lambda = [5\,3\,2^2\,1^4]$ (see Example~\ref{ex.53221111}). The table below computes the defect of each symbol, and the Levi subgroup for the corresponding cuspidal datum. 
\[
\begin{array}{c|c|c}
\text{ Symbol } & \text{ Defect } & L\\ \hline
(\{0,2,5,7\},\{1,3,5,9\}) & 0 & T = GL_1\times GL_1 \times GL_1\times \cdots \times GL_1  \\
(\{0,2,5,9\},\{1,3,5,7\}) & 0 & T = GL_1\times GL_1 \times GL_1\times \cdots \times GL_1   \\
(\{1,3,5 ,7, 9\},\{0,2,5\}) & 2 & GL_1\times GL_1 \times GL_1\times GL_1 \times GL_1\times GL_1\times \SO(4)   \\
(\{0 ,2 ,5 ,7, 9\},\{1,3,5\}) & 2 & GL_1\times GL_1 \times GL_1\times GL_1 \times GL_1\times GL_1\times \SO(4)   \\
\end{array}
\]
Note that the Levi subgroup in the last two rows has type $A_1\times A_1$. 
\end{example}

\begin{Cor}\label{cor.defectBD} Let $\co$ be a nilpotent $G$-orbit and $\cl$ a $G$-equivariant local system on $\co$.  Let $\{A,B\}$ denote the symbol equal to the image of the pair $(\co, \cl)$ under the map appearing in~\eqref{eqn.int.bijectionBD2}.  Then the cuspidal datum $\mathbf{c}$ corresponding to $\{A,B\}$ is: 
\begin{itemize}
\item $\mathbf{c}=(T, \{0\}, \cl_\mathrm{triv})$ if and only if the defect of $\{A,B\}$ is $0$ or $1$ (that is, $IC(\overline{\co}, \cl)$ appears in the Springer sheaf), and 
\item $\mathbf{c} = (L_\chi, \co_{\reg}, \cl^{L_\chi})$, where $L_\chi$ is the Levi subgroup of type $A_1\times A_1$ corresponding to the character $\chi$ appearing in the type $D_n$ rows of Table~\ref{fig.Levi.table}, if and only if the defect of $\{A,B\}$ is~$2$.  
\end{itemize}
In particular, $IC(\overline{\co}, \cl)$ appears as a summand of a Lusztig sheaf in $\psi_*\shc_{\tilde{\cm}}[\dim \cn]$ if and only if the defect of the corresponding symbol is $0$, $1$, or~$2$.
\end{Cor}
\begin{proof} 

By Theorem~\ref{thm.Lusztig.TypeBD}, a symbol with defect $0$ corresponds to a cuspidal datum $\mathbf{c}$ with the Levi factor
$L_{N/2} \simeq T$ and a symbol with defect $1$ corresponds to the Levi subgroup $L_{(N-1)/2} \simeq T$, since in this case $\SO(N-2i) = \SO(1)$. Note that a defect of $0$ only occurs when $N$ is even and a defect of $1$ only occurs when $N$ is odd, so these indices make sense.   Hence, in each case, the corresponding cuspidal
datum is $\mathbf{c}=(T, \{0\}, \cl_\mathrm{triv})$.   

If the defect is $2$, then $N$ must be even, that is, $G$ is of type $D_{N/2}$, and the corresponding Levi subgroup
is 
\[
L_{(N-4)/2}\simeq \underbrace{GL_1\times \cdots \times GL_1}_{\text{$(N-4)/2$ factors}}\times  \SO(4) = \underbrace{GL_1\times \cdots \times GL_1}_{\text{$(N-4)/2$ factors}}\times  A_1\times A_1.
\]
Finally, when the defect of $\{A,B\} \in \Psi_N'$ is greater than $2$, it is easy to confirm that the corresponding Levi subgroup $L_i$ as in Theorem~\ref{thm.Lusztig.TypeBD} satisfies $L_i\not\simeq T$ and $L_i$ is not of type $A$. Since $\mathbf{c}$ is minimal if and only if the Levi factor of $\mathbf{c}$ is of type $A$, this completes the proof. 
\end{proof}

\begin{Prop}\label{typeBDdefect} Suppose $\lambda\in X_{N}^{BD}$.  If $N$ is odd, the defect of every symbol in the similarity class of $\{A_\lambda, B_\lambda\}$ is 
 $1$ if and only if $\co_\lambda$ is Levi-regular. If $N$ is even, the defect of every symbol in the similarity class of $\{A_\lambda, B_\lambda\}$ is $0$ or $2$ if and only if either $\co_\lambda$ is Levi-regular or $\lambda$ has exactly two odd parts with odd multiplicity, both of which are greater than $1$.   
\end{Prop}

\begin{proof}  Let $\lambda\in X_N^{BD}$. By definition of a symbol, if $i$ and $i+1$ are elements of some interval $I$ in $\ci_\lambda$, then $i\in A_\lambda$ and $i+1\in B_\lambda$ or vice versa. Thus, 
\begin{eqnarray}\label{eqn.BDformula}
|I\cap A_\lambda| - |I\cap B_\lambda| = \left\{ \begin{array}{ll} \pm 1 & \text{ if $|I|$ is odd} \\
0 & \text{ if $|I|$ is even. } \end{array}\right.
\end{eqnarray}

We consider two possible cases, according to the parity of $N$.

\medskip 

\textit{Type B case: $N$ is odd.}  Let $N = 2n+1$ and $\lambda\in X_{2n+1}^{BD}$.  By Lemma~\ref{lem.g2orbits}, $\co_\lambda\cap \fg_2\neq \emptyset$ if and only if $\lambda$ has at most one odd part with odd multiplicity.  We aim to show that $\lambda$ satisfies this property if and only if every symbol in the similarity class of $\{A_\lambda, B_\lambda \}$ has defect $1$. Recall that the defect must be a positive odd integer in this case.

Suppose $\lambda$ has $k\geq 1$ distinct odd parts occurring with odd multiplicity.  Since $N$ is odd, 
$k$ must be an odd integer which is at least $1$.  By Lemma~\ref{lem.int.bijectionBD} there are exactly $k$ intervals in the set $\ci_\lambda$ with odd cardinality.  
Since the defect of $\{A_\lambda, B_\lambda\}$ is $1$, we have $| |A_\lambda| - |B_\lambda| | = 1$.  In the argument below, we can interchange 
the roles of $A_\lambda$ and $B_\lambda$ if necessary, and therefore we may assume $ |A_\lambda| - |B_\lambda| = 1$.  Then \eqref{eqn.BDformula} implies that
\begin{eqnarray}\label{eqn.Bformula}
1 = |A_\lambda| - |B_\lambda| =  \sum_{\substack{I\in \ci_\lambda \\ |I| \text{ odd}}} \left( |I\cap A_\lambda| - |I\cap B_\lambda| \right).
\end{eqnarray}

If $k = 1$, there is only one interval in $\ci_\lambda$ with odd length, call it $I$, and \eqref{eqn.Bformula} becomes
$$
1 = |A_\lambda| - |B_\lambda| =  |I\cap A_\lambda| - |I\cap B_\lambda|.
$$
If $I \in S$, then  $ |I \cap A_{\lambda, S}| = |I\cap A_\lambda|$ and $|I \cap B_{\lambda, S}| = |I\cap B_\lambda|$;
if $I \not\in S$, then $ |I \cap A_{\lambda, S}| = |I\cap B_\lambda|$ and $|I \cap B_{\lambda, S}| = |I\cap A_\lambda|$.
Hence for any $S \subset \ci_\lambda$, we have
$$
|A_{\lambda, S}| - |B_{\lambda, S}| = \pm ( |A_\lambda| - |B_\lambda| ) = \pm 1,
$$
so the defect of any $\{ A_{\lambda, S}, B_{\lambda, S} \}\in \Psi_{N,\lambda}'$ is $| |A_{\lambda, S}| - |B_{\lambda, S}| | = 1$.

Now suppose $k\geq 3$, which holds whenever $k\neq 1$ since $k$ must be odd. By~\eqref{eqn.BDformula} and \eqref{eqn.Bformula}, there must be some interval $I'\in \ci_\lambda$ with odd length whose intersection with $B_\lambda$ is larger than the intersection with $A_\lambda$, that is
\[
|I' \cap A_\lambda| - |I' \cap B_\lambda| = -1. 
\]
Now set $S=\ci_\lambda \setminus \{I'\}$ and consider the symbol $\{A_{\lambda, S}, B_{\lambda, S}\}$.  By definition of the symbol $\{A_{\lambda, S}, B_{\lambda, S}\}$ from~\eqref{eqn.int.bijection.def}, $ I'\cap A_{\lambda, S}  =I' \cap B_\lambda $ and $I' \cap B_{\lambda, S} =I'\cap  A_\lambda$.  We get 
\begin{eqnarray*}
|A_{\lambda, S}| - |B_{\lambda, S}| &=& \left(  \sum_{\substack{I \in S \\ |I| \text{ odd}}} \left( |I\cap A_\lambda| - |I\cap B_\lambda| \right) \right) + \left( |I'\cap B_\lambda| - |I'\cap A_\lambda| \right) \\
&=& \sum_{\substack{I \in \ci_\lambda \\ |I| \text{ odd}}}\left( |I\cap A_\lambda| - |I\cap B_\lambda| \right) + 2\left( |I'\cap B_\lambda| - |I'\cap A_\lambda| \right) \\
&=& 1+2=3.
\end{eqnarray*}
Thus $||A_{\lambda, S}| - |B_{\lambda, S}||=3$, so when $\lambda$ has more than $1$ odd part with odd multiplicity there is some symbol in the similarity class of $\{A_\lambda, B_\lambda\}$ with defect greater than $1$.

\medskip

\textit{Type D case: $N$ is even.} Let $N=2n$ and $\lambda \in X_{2n}^{BD}$. By Lemma~\ref{lem.g2orbits}, $\co_\lambda$ is Levi-regular
if and only if $\lambda$ has at most two odd parts with odd multiplicity, one of which is equal to $1$. Thus, we aim to show that $\lambda$ has at most two odd parts with odd multiplicity if and only if every symbol in the similarity class of $\{A_\lambda, B_\lambda\}$ has defect $0$ or $2$. Recall that the defect must be a nonnegative even integer in this case.

Suppose $\lambda$ has $k\geq 0$ distinct odd parts with odd multiplicity (note that $k$ must be even since $N$ is even). By Lemma~\ref{lem.int.bijectionBD} there are exactly $k$ intervals in the set $\ci_\lambda$ with odd cardinality.  
Since the defect of $\{A_\lambda, B_\lambda\}$ is $0$ and since $\{A_\lambda, B_\lambda\} =\{A_{\lambda,S}, B_{\lambda,S}\}$ for $S=\ci_\lambda$ we have that
\begin{eqnarray}\label{eqn.Dformula}
0 =  |A_\lambda| - |B_\lambda|  =  \sum_{\substack{I\in \ci_\lambda \\ |I| \text{ odd}}} \left( |I\cap A_\lambda| - |I\cap B_\lambda| \right).
\end{eqnarray}
Equation~\eqref{eqn.BDformula} tells us that the only changes in defect from $\{A_\lambda, B_\lambda\}$ to $\{A_{\lambda, S}, B_{\lambda, S}\}$ can occur from intervals whose size is odd.  Hence if $k=0$, the defect of any symbol $\{A_{\lambda, S}, B_{\lambda, S}\}$ in the similarity class of $\{A_\lambda, B_\lambda\}$ is $0$.

If $k=2$, let $I'$ and $I''$ denote the corresponding intervals in $\ci_\lambda$ with odd length. Equation~\ref{eqn.Dformula} becomes
\begin{eqnarray}\label{eqn.Dformula2}
0= \left( |I'\cap A_\lambda| - |I' \cap B_\lambda| \right) +(|I''\cap A_\lambda| - |I''\cap B_\lambda|).
\end{eqnarray}
Let $S\subseteq \ci_\lambda$. If $\{I',I''\} \subseteq S$ then the elements of these intervals are fixed when we flip elements between $A_\lambda$ and $B_\lambda$ to create the symbol $\{A_{\lambda, S}, B_{\lambda, S}\}$; see~\eqref{eqn.int.bijection.def}. The defect of $\{A_{\lambda, S}, B_{\lambda, S}\}$ is equal to that of $\{A_\lambda, B_\lambda\}$ in this case, so it is $0$. Similarly, if $S\cap \{I',I''\} = \emptyset$, then all elements of $I'$ and $I''$ are switched when we
go from  $\{ A_\lambda, B_\lambda \}$ to $\{A_{\lambda, S}, B_{\lambda, S}\}$, so \eqref{eqn.Dformula2}
implies that defect of $\{A_{\lambda, S}, B_{\lambda, S}\}$ is $0$.

Now suppose $I' \in S$ and $I''\notin S$. We get that 
\begin{eqnarray*}
|A_{\lambda, S}| - |B_{\lambda, S}| &=& \left( |I' \cap A_\lambda| - |I'\cap B_\lambda| \right) + \left( |I''\cap B_\lambda| - |I'' \cap A_\lambda| \right) \\
&=& 2 \left(  |I' \cap A_\lambda| - |I'\cap B_\lambda|  \right) = \pm 2,
\end{eqnarray*}
where the last two equalities follow from~\eqref{eqn.Dformula2} and~\eqref{eqn.BDformula}, respectively. Thus, when $\lambda$ has at most two odd parts with odd multiplicity, the defect of any symbol in the equivalence class of $\{A_\lambda, B_\lambda\}$ is at most $2$. 

Finally, suppose $k\geq 4$. In this case, there exist intervals $I',I'' \in \ci_\lambda$ of odd length such that 
\begin{eqnarray}\label{eqn.Dformula3}
|I\cap A_\lambda| - |I\cap B_\lambda| = -1 \; \text{ for $I=I'$ and $I=I''$}.
\end{eqnarray}
Let $S=\ci_\lambda \setminus \{I',I''\}$. The elements of $I'$ and $I''$ are swapped between $A_\lambda$ and $B_\lambda$ when creating the symbol $\{A_{\lambda, S} , B_{\lambda, S} \}$. Thus, 
\begin{eqnarray*}
|A_{\lambda, S}| - |B_{\lambda, S}| &=& \left( \sum_{\substack{I\in S\\ |I|\text{ odd}}} (|I\cap A_\lambda| - |I\cap B_\lambda|) \right) \\
&& \quad\quad + \left( |I'\cap B_\lambda| - |I'\cap A_\lambda| \right) + \left( |I''\cap B_\lambda|-|I''\cap A_\lambda| \right) \\
&=& \left( \sum_{\substack{I\in \ci_\lambda\\ |I|\text{ odd}}} (|I\cap A_\lambda| - |I\cap B_\lambda|) \right) \\
&& \quad\quad + 2\left( |I'\cap B_\lambda| - |I'\cap A_\lambda| \right) + 2\left( |I''\cap B_\lambda|-|I''\cap A_\lambda| \right) \\
&=& 4
\end{eqnarray*}
where the final equality follows by an application of~\eqref{eqn.Dformula2} and~\eqref{eqn.Dformula3}. Thus we have proven that if $\lambda$ has more than $2$ parts with odd multiplicity, then there is some symbol in the similarity class of $\{A_\lambda, B_\lambda\}$ with defect greater than $2$.
\end{proof}

\begin{Prop}\label{prop.orthogonal.odd} Suppose $G=\SO(2n+1)$. Every simple perverse sheaf of the form $IC(\overline{\co}, \cl)$  for $\cl$ a $G$-equivariant local system on $\co$ appears as a summand of $\psi_*\shc_{\tilde{\cm}}[\dim \cn]$ if and only if $\co$ is Levi-regular.
\end{Prop}

\begin{proof}
Let $\lambda\in X_{N}^{BD}$.   By Lemma~\ref{lem.int.bijectionBD} the subset of $\Psi_{N}'$ corresponding to pairs $(\co_\lambda, \cl)$ is precisely the similarity class $\Psi_{N,\lambda}'$ of $\{A_\lambda, B_\lambda\}$. 
By Corollary~\ref{cor.defectBD}, $IC(\overline{\co}_\lambda, \cl)$ is a summand of $\psi_*\shc_{\tilde{\cm}}[\dim \cn]$ if and only if the defect of the symbol corresponding to the pair $(\co_\lambda, \cl)$ is $0$, $1$, or $2$; since $N$ is odd, the defect must be $1$. 
By Proposition~\ref{typeBDdefect} this occurs precisely when  $\co_\lambda$ is Levi-regular.
\end{proof}

\begin{Prop}\label{prop.orthogonal.even} Suppose $G=\SO(2n)$. Every simple perverse sheaf of the form $IC(\overline{\co}, \cl)$  for $\cl$ a $G$-equivariant local system on $\co$ appears as a summand of $\psi_*\shc_{\tilde{\cm}}[\dim \cn]$ if and only if $\co$ is Levi-regular or if $\co=\co_\lambda$ for $\lambda\in X_{2n}^{BD}$ with exactly two odd parts with odd multiplicity, both of which are greater than $1$.
\end{Prop}
\begin{proof}  Let $N=2n$ and $\lambda\in X_{N}^{BD}$.   By Lemma~\ref{lem.int.bijectionBD} the subset of $\Psi_{N}'$ corresponding to pairs $(\co_\lambda, \cl)$ is precisely the similarity class $\Psi_{N,\lambda}'$ of $\{A_\lambda, B_\lambda\}$. 
By Corollary~\ref{cor.defectBD}, $IC(\overline{\co}_\lambda, \cl)$ is a summand of $\psi_*\shc_{\tilde{\cm}}[\dim \cn]$ if and only if the defect of the symbol corresponding to the pair $(\co_\lambda, \cl)$ is $0$, $1$, or $2$; since $N$ is even, the defect must be $0$ or $2$. 
By Proposition~\ref{typeBDdefect} this occurs precisely when $\lambda$ has at most two odd parts with odd multiplicity. Thus by Lemma~\ref{lem.g2orbits} we have that either $\co_\lambda$ is Levi-regular or that both odd parts of $\lambda$ with odd multiplicity are greater than $1$. 
\end{proof}

\subsection{The spin groups $\Spin(2n)$ and $\Spin(2n+1)$}  \label{ss.spin}
We now focus on the extension of the generalized Springer correspondence to the simply-connected spin groups in types $B_n$ and $D_n$.  This was studied by Lusztig in \cite[\S 14]{Lu84} and by Lusztig--Spaltenstein in \cite{LuSpal85}.  
Suppose that $G = \Spin(N)$. As discussed in Section~\ref{ss.strategy} above,  $\SO(N)= G/Z_0$, and every $\SO(N)$-equivariant local system on a nilpotent orbit $\co$ is automatically $G$-equivariant.  
If we fix a character $\chi$ of the center $Z$ of $G$ that acts trivially on $Z_0$, then any cuspidal datum corresponding to $\chi$ is the same for $G$ as for $\SO(N)$.
The simple perverse sheaves $IC(\overline{\co},\cl)$ which appear in the corresponding Lusztig sheaf are the same as for $\SO(N)$, and were determined in the previous section.  

Suppose instead that $\chi$ is a character of $Z$ such that $\chi|_{Z_0}$ is not trivial.  Such a character $\chi$ then is associated to one or more cuspidal data that are not obtained in the generalized Springer correspondence for $\SO(N)$.  Given an orbit $\co$, there is at most one simple perverse sheaf
$IC(\overline{\co}, \cl)$ with the central character $\chi$.  This perverse sheaf necessarily appears in exactly one cuspidal datum.  
The local system on $\co$ is uniquely determined by $\chi$, so
it is not
necessary to consider multiple local systems on $\co$.
This simplifies the discussion and avoids the need for symbol computations.  
The reader may wish to refer to Appendix \ref{s.overview}, which describes the generalized Springer correspondence in all types, with specific examples where $G = \Spin(N)$.

We now work to complete the proof of Theorem~\ref{thm.support} in types $B_n$ and $D_n$. Suppose $\lambda\in X_N^{BD}$, that is, $\lambda$ is a partition of $N$ such that every even part occurs with even multiplicity.  Lusztig and Spaltenstein~\cite{LuSpal85} characterize the nilpotent orbits such that additional sheaves $IC(\overline{\co}, \cl)$ occur for the generalized Springer correspondence of $G=\Spin(N)$ using a numerical invariant of the partition $\lambda$. In order to state their results, we define some notation.

Following~\cite[\S0.10]{LuSpal85}, given an integer $d$ such that $d\equiv N\mod 4$ and $d(2d-1)\leq N$ we define a Levi subgroup $L_d$ of $\Spin(N)$ as follows. If $N$ is odd, then $L_d$ is a standard Levi subgroup of type 
\begin{equation}\label{eqnBLd}
B_{(d-1)(2d+1)/2}\times A_1\times \cdots\times A_1
\end{equation}
with $\frac{1}{4}\left( N-d(2d-1)\right)$ factors of type $A_1$. Note that when $d=1$ (which only occurs when $N=2n+1$ for $n$ even) the Levi subgroup has associated Dynkin diagram:
\begin{center}
\begin{tikzpicture}[scale=1]
\node[circle, fill=black, draw, inner sep=2pt] (1) at (1,0) {};
\node[circle, draw, inner sep=2pt] (2) at (2,0) {};
\node[circle, fill=black, draw, inner sep=2pt] (3) at (3,0) {};
\node[circle, draw, inner sep=2pt] (4) at (4,0) {};
\node[circle, fill=black, draw, inner sep=2pt] (5) at (5,0) {};
\node[inner sep=-0.5pt, outer sep=-0.5pt] (7) at (5.6,0) {};
\node[circle, draw, inner sep=2pt] (6) at (6,0) {};
\draw (1)--(2);
\draw (2)--(3);
\draw (3)--(4);
\draw[dashed] (4)--(5);
\draw[double, double distance=1.5pt, ->] (5)--(7);
\draw[double, double distance=1.5pt] (7)--(6);
\end{tikzpicture} 
\end{center}
and when $d=-1$  (which only occurs when $N=2n+1$ for $n$ odd) the Levi subgroup has the following associated Dynkin diagram.
\begin{center}
\begin{tikzpicture}[scale=1]
\node[circle, fill=black, draw, inner sep=2pt] (1) at (1,0) {};
\node[circle, draw, inner sep=2pt] (2) at (2,0) {};
\node[circle, fill=black, draw, inner sep=2pt] (3) at (3,0) {};
\node[circle, draw, inner sep=2pt] (4) at (4,0) {};
\node[circle, draw, inner sep=2pt] (5) at (5,0) {};
\node[inner sep=-0.5pt, outer sep=-0.5pt] (7) at (5.6,0) {};
\node[circle, fill=black, draw, inner sep=2pt] (6) at (6,0) {};
\draw (1)--(2);
\draw (2)--(3);
\draw (3)--(4);
\draw[dashed] (4)--(5);
\draw[double, double distance=1.5pt, ->] (5)--(7);
\draw[double, double distance=1.5pt] (7)--(6);
\end{tikzpicture}
\end{center}

If $N$ is even, then $L_d$ is a standard Levi subgroup of type 
\begin{equation}\label{eqnDLd}
D_{d(2d-1)/2}\times A_1\times \cdots \times A_1
\end{equation}
with $\frac{1}{4}\left( N-d(2d-1)\right)$ factors of type $A_1$. This uniquely determines $L_d$ up to conjugacy, except when $d=0$ (which only occurs when $N=2n$ for $n$ even).  In this case, the Levi subgroup $L_0$ has type $A_1\times A_1\times \cdots \times A_1$ and the Dynkin diagram for $L_0$ has one of the following two forms.
\begin{center}
\begin{tikzpicture}[scale=1]
\node[circle, fill=black, draw, inner sep=2pt] (1) at (1,0) {};
\node[circle, draw, inner sep=2pt] (2) at (2,0) {};
\node[circle, fill=black, draw, inner sep=2pt] (3) at (3,0) {};
\node[circle, draw, inner sep=2pt] (4) at (4,0) {};
\node[circle, fill= black, draw, inner sep=2pt] (5) at (5,0) {};
\node[circle, draw, inner sep=2pt] (6) at (6,0) {};
\node[circle, fill=black, draw, inner sep=2pt] (7) at (7,0.3) {};
\node[circle, draw, inner sep=2pt] (8) at (7,-0.3) {};
\draw (1)--(2);
\draw (2)--(3);
\draw[dashed] (3)--(4);
\draw (4)--(5);
\draw (5)--(6);
\draw (6)--(7);
\draw (6)--(8);
\end{tikzpicture}\\
\end{center}
\begin{center}
\begin{tikzpicture}[scale=1]
\node[circle, fill=black, draw, inner sep=2pt] (1) at (1,0) {};
\node[circle, draw, inner sep=2pt] (2) at (2,0) {};
\node[circle, fill=black, draw, inner sep=2pt] (3) at (3,0) {};
\node[circle, draw, inner sep=2pt] (4) at (4,0) {};
\node[circle, fill= black, draw, inner sep=2pt] (5) at (5,0) {};
\node[circle, draw, inner sep=2pt] (6) at (6,0) {};
\node[circle,  draw, inner sep=2pt] (7) at (7,0.3) {};
\node[circle, fill=black, draw, inner sep=2pt] (8) at (7,-0.3) {};
\draw (1)--(2);
\draw (2)--(3);
\draw[dashed] (3)--(4);
\draw (4)--(5);
\draw (5)--(6);
\draw (6)--(7);
\draw (6)--(8);
\end{tikzpicture}
\end{center}
Which form appears depends on the central character $\chi$. 

\begin{example}\label{ex.Spin.Levi} Consider $G=\Spin(10)$.  The integers $d$ satisfying $d\equiv 10\mod 4$ and $d(2d-1)\leq 10$ are $d=\pm 2$. When $d=2$, we have $d(2d-1) = 6$ and when $d=-2$, we have $d(2d-1) =10$. Thus, 
\[
L_{2} \simeq D_3 \times A_1  \; \text{ and } \; L_{-2} \simeq D_{5}. 
\]
Note that $L_{-2} = \Spin(10)$ is not a type $A$ Levi subgroup of $G$, while $L_2$ is a type $A$ Levi since $D_3\simeq A_3$. 
\end{example}

The next lemma follows by inspection:

\begin{Lem}\label{lem.Spin.Levi} Let $d$ be an integer such that $d\equiv N\mod 4$ and $d(2d-1)\leq N$. The Levi subgroup $L_d$ is of type $A$ if and only if $d\in \{0,-1,1,2\}$. 
\end{Lem}

Let $K_N \subseteq X_N^{BD}$ be the subset of \emph{rather odd partitions}, that is, partitions such that all even parts occur with even multiplicity and such that each odd part occurs with multiplicity at most one. Note that very even partitions (the partitions in $X_N^{BD}$ with only even parts) are considered to be rather odd.  Recall that very even partitions only occur when $N=2n$ with $n$ even, and each such partition $\lambda$ corresponds to two nilpotent orbits, denoted below by $\co_\lambda^+$ and $\co_\lambda^-$.

Define an integer $d(\lambda)$ for each $\lambda= [\lambda_1\ \lambda_2\ \dots\ \lambda_l] \in K_N$ by 
\[
d(\lambda) = \displaystyle \sum_{1\leq i\leq l} d(\lambda_i),
\] 
where 
\begin{equation*}
d(\lambda_i)=\begin{cases} 0 & \text{if } \lambda_i\text{ is even} \\ 1 & \text{if } \lambda_i\equiv 1\mod 4 \\ -1 & \text{if }\lambda_i\equiv  -1 \mod 4.\end{cases} 
\end{equation*}
It is straightforward to check that  $d(\lambda)\equiv N \mod 4$ and $d(\lambda)(2d(\lambda)-1)\leq N$. Note that $d(\lambda)=0$ whenever $\lambda$ is very even. We can now summarize the key results of~\cite{LuSpal85}.

\begin{Thm}[Lusztig--Spaltenstein] \label{thm.Lusztig-Spaltenstein} Suppose $G=\Spin(N)$ and let $\chi$ be a fixed central character acting nontrivially on $Z_0 = \ker(G\to \SO(N))$. For each integer $d$ such that $d\equiv N \mod 4$ and $d(2d-1)\leq N$ there exists a unique cuspidal datum $\mathbf{c}_d:=(L_d, \co_{L_d}, \cl^{L_d})$, up to conjugacy, with central character $\chi$.  The simple perverse sheaves occurring as summands of $\mathbb{A}_{\mathbf{c}_d}$ are determined as follows.
\begin{enumerate}
\item If $IC(\overline{\co}_\lambda, \cl)$ is a summand of $\mathbb{A}_{\mathbf{c}_d}$, then $\lambda\in K_N$. 
\item When $\lambda\in K_N$ is not very even, there is exactly one perverse sheaf of the form  $IC(\overline{\co}_\lambda, \cl)$ occurring as a summand of $\mathbb{A}_{\mathbf{c}_d}$ for some $d$, namely $d=d(\lambda)$.
\item When $\lambda\in K_N$ is very even, there is exactly one perverse sheaf of the form either  $IC(\overline{\co}^+_\lambda, \cl)$ or $IC(\overline{\co}^-_\lambda, \cl)$ occurring as a summand of $\mathbb{A}_{\mathbf{c}_d}$ for some $d$, namely $d=d(\lambda)=0$. 
\end{enumerate}
\end{Thm}

\begin{example}\label{ex.Spin10} Let $N=10$ and consider $\lambda = [7 \ 3]$. Since $\lambda$ is rather odd, $\lambda\in K_{10}$ and we have $d(\lambda) = d(3)+d(7) = -2$. By Theorem~\ref{thm.Lusztig-Spaltenstein}, there is precisely one sheaf $IC(\overline{\co}_\lambda, \cl)$ occurring in the generalized Springer correspondence for $\Spin(10)$ for each central character $\chi$ acting nontrivially on $Z_0$, and the Levi factor of its cuspidal datum is $L_{d(\lambda)} =L_{-2}$. Since $L_{-2}$ is not type $A$, the cuspidal datum of $IC(\overline{\co}_\lambda, \cl)$ is not minimal, and it follows that not all simple perverse sheaves for the orbit $\co_\lambda$ occur as summands of $\psi_*\shc_{\tilde{\cm}}[\dim \cn]$. 
\end{example}

The next lemma establishes which of the partitions from $K_N$ correspond to Levi-regular orbits. The result follows directly from Lemma~\ref{lem.g2orbits}.

\begin{Lem}\label{Lem:partitionxn}
Suppose $\lambda\in K_N$. Then $\co_\lambda$ is Levi-regular if and only if 
\begin{itemize}
\item \textup{Type $B_n$:}  $N$ is odd and $\lambda$ has precisely one odd part with multiplicity one and all other parts of $\lambda$ are even and occur with even multiplicity.
\item \textup{Type $D_n$:}  $N$ is even and either $\lambda$ is very even or  $\lambda$ has precisely two odd parts, an odd part $p > 1$ with multiplicity one and an odd part equal to $1$ occurring with multiplicity one, and all other parts of $\lambda$ are even and occur with even multiplicity.
\end{itemize}
\end{Lem}

\begin{Lem}\label{lem.spinsoln}
Suppose $\lambda \in K_N$.  
\begin{enumerate}
\item If $N = 2n+1$ with $n$ even and $\lambda$ is Levi-regular, then $d(\lambda) = 1$.
\item If $N = 2n+1$ with $n$ odd and $\lambda$ is Levi-regular, then $d(\lambda) = -1$. 
\item If $N = 2n$ with $n$ even and $\lambda$ has at most two odd parts, then at least one odd part is congruent to $1$ modulo $4$ and  $d(\lambda) = 0$.
\item If $N = 2n$ with $n$ odd and $\lambda$ has at most two odd parts, then $d(\lambda) = 2$ if and only if at least one odd part is congruent to $1$ modulo $4$, otherwise $d(\lambda)=-2$. 
\end{enumerate}
\end{Lem}

\begin{proof}
First suppose $N = 2n+1$ and $\lambda$ is Levi-regular.  By Lemma~\ref{Lem:partitionxn}, $\lambda$ has one odd part with multiplicity one and all other parts are even.  Denote the single odd part by $p$.  Since $d(\lambda)=d(p)$, we need only determine whether $p\equiv 1 \mod 4$ or $p\equiv 3\mod 4$ to compute $d(\lambda)$.  The sum of even parts with even multiplicities will always be divisible by $4$.  Thus, $p \equiv N\mod 4$.  This proves parts (1) and (2).

Now suppose $N = 2n$.  If $\lambda$ is very even, $d(\lambda)=0$. Assume $\lambda$ has exactly two distinct odd parts. Denote these parts by $p_1=2m_1+1$ and $p_2=2m_2+1$ for some $m_1, m_2\in \Z$. Since every even part of $\lambda$ occurs with even multiplicity, there exists $k\in \Z$ such that 
\begin{align}\label{eqn.parity}
N = 4k + 2(m_1+m_2)+2 \Rightarrow n = 2k+ m_1+m_2+1.
\end{align}
If $n$ is even, then $m_1+m_2$ is odd and $m_1$ and $m_2$ must have opposite parity. It follows that at least one of the odd parts of $\lambda$ is congruent to $1$ modulo $4$ and $d(\lambda) = d(p_1)+d(p_2)=1+(-1) = 0$. This proves part (3) of the lemma.

Now assume $n$ is odd. Equation~\eqref{eqn.parity} now implies $m_1+m_2$ is even so $m_1$ and $m_2$ have the same parity. If one of $m_1$ or $m_2$ is even, then both are. Equivalently, if one odd part of $\lambda$ is congruent to $1$ modulo $4$, then both are and we have $d(\lambda) = d(p_1)+d(p_2)=1+1=2$. Both $m_1$ and $m_2$ are odd if and only if $p_1$ and $p_2$ are both congruent to $-1$ modulo $4$, i.e., $d(\lambda) = d(p_1)+d(p_2)=-1-1=-2$. This completes the proof of part (4).
\end{proof}

We are now ready to prove Theorem~\ref{thm.support} in type $B_n$.  When $N=2n+1$, the center of $G=\Spin(N)$ is $\Z_2$ and the nontrivial character $\chi\in \hat{Z}$ is the only one which remains nontrivial when restricted to $Z_0$. 

\begin{proof}[Proof of Theorem~\ref{thm.support} in type $B_n$] Let $G= \Spin(N)$ with $N=2n+1$. Using Theorem~\ref{t:main1} and  Proposition~\ref{prop.orthogonal.odd}, it suffices to show that if $\co$ is Levi-regular and $\cl$ is a $G$-equivariant local system on $\co$ which is not $\SO(2n+1)$-equivariant, then the cuspidal datum of $IC(\overline{\co}_\lambda, \cl)$ is minimal. By Theorem~\ref{thm.Lusztig-Spaltenstein} applied to the character $\chi$, there is only one such simple perverse sheaf, $\co=\co_\lambda$ for $\lambda\in K_{N}$, and the Levi factor of the cuspidal datum of $IC(\overline{\co}_\lambda, \cl)$ is $L_{d(\lambda)}$. Lemmas~\ref{lem.Spin.Levi} and~\ref{lem.spinsoln} together imply that $L_{d(\lambda)} = L_{\pm 1}$ is a type $A$ Levi (which is necessarily $L_\chi$). Since a cuspidal datum is minimal if and only if the Levi factor has type $A$, this completes the proof. 
\end{proof}

Next, we prove Theorem~\ref{thm.support} in type $D_n$.  Let $G=\Spin(N)$ with $N=2n$.  If $n$ is even, then $Z \simeq \Z_2\times \Z_2$; if $n$ is odd, then $Z\simeq \mathbb Z_4$.
In either case, we have four central characters to consider.  We will denote these characters using the notation established in Section~\ref{ss:Dn}.  By Lemma~\ref{lem.kernel}, the trivial character and $\chi_1$ are both trivial when restricted to $Z_0$, so all Lusztig sheaves for these central characters appear as part of the generalized Springer correspondence for $\SO(N)$.  The characters $\chi_{n-1}$ and $\chi_{n}$ each have nontrivial restrictions to $Z_0$.

\begin{proof}[Proof of Theorem~\ref{thm.support} in type $D_n$ for $n$ even] Let $G= \Spin(N)$ where $N=2n$ with $n$ even. Suppose $\lambda\in K_N$ is not very even and has two odd parts. By Theorem~\ref{thm.Lusztig-Spaltenstein} applied to each of the characters $\chi_{n-1}$ and $\chi_n$, there is a unique simple perverse sheaf of the form $IC(\overline{\co}_\lambda, \cl)$ with cuspidal datum $\mathbf{c}_{d(\lambda)}$ for each central character. By Lemma~\ref{lem.spinsoln}, at least one of the odd parts of $\lambda$ is congruent to $1$ modulo $4$. Thus, in order to complete the proof of Theorem~\ref{thm.support} for such $\co_\lambda$, it suffices by Proposition~\ref{prop.orthogonal.even} to show that the cuspidal datum of $IC(\overline{\co}_\lambda, \cl)$ for each of $\chi_{n-1}$ and $\chi_n$ has a type $A$ Levi factor. Lemmas~\ref{lem.Spin.Levi} and~\ref{lem.spinsoln} together imply that $L_{d(\lambda)} = L_{0}$ is indeed a type $A$ Levi (which is necessarily $L_{\chi_{n-1}}$ or $L_{\chi_n}$ depending on our choice of central character).

Now suppose $\lambda\in K_N$ is very even. By the same reasoning as above, there is a unique simple perverse sheaf either of the form $IC(\overline{\co}_\lambda^+, \cl)$ or of the form $IC(\overline{\co}_\lambda^-, \cl)$ to consider. Since $d(\lambda)=0$, the Levi part of the associated cuspidal datum has type $A$ by Lemma~\ref{lem.Spin.Levi} (which is, once again, necessarily $L_{\chi_{n-1}}$ or $L_{\chi_n}$ depending on our choice of central character). This concludes the proof.
\end{proof}

\begin{proof}[Proof of Theorem~\ref{thm.support} in type $D_n$ for $n$ odd]  Let $G=\Spin(2n)$ with $n$ odd. The assumption $n$ odd implies the set $X_N^{BD}$, and consequently $K_N$, does not contain any very even partitions. 

By Theorem~\ref{thm.Lusztig-Spaltenstein} applied to each of the characters $\chi_{n-1}$ and $\chi_n$, there is precisely one simple perverse sheaf $IC(\overline{\co}_{\lambda}, \cl)$  that occurs in the generalized Springer correspondence for $G$ but not for $SO(2n)$ for each $\lambda \in K_N$ and for each central character $\chi_{n-1}$ and $\chi_n$. The Levi component of the associated cuspidal datum $\mathbf{c}_{d(\lambda)}$ is $L_{d(\lambda)}$.  

Suppose $\lambda\in X_N^{BD}$ has two distinct odd parts with odd multiplicity. By Proposition~\ref{prop.orthogonal.even}, the proof of Theorem~\ref{thm.support} is complete for the orbit $\co_\lambda$ unless the multiplicity of both odd parts is $1$, in which case $\lambda \in K_N$. Suppose at least one odd part of $\lambda$ is congruent to $1$ modulo $4$. Then $d(\lambda)=2$ by Lemma~\ref{lem.spinsoln}. The Levi part $L_2$ of the cuspidal datum $\mathbf{c}_{d(\lambda)}$ is of type $A$ by Lemma~\ref{lem.Spin.Levi} (which is necessarily $L_{\chi_{n-1}}$ or $L_{\chi_n}$ depending on our choice of central character), and thus $IC(\overline{\co}_{\lambda}, \cl)$ is a summand of $\psi_*\shc_{\tilde{\cm}}[\dim \cn]$. Otherwise, the same lemmas tell us $d(\lambda)=-2$ and the cuspidal datum for $IC(\overline{\co}_{\lambda}, \cl)$ is not minimal, so we conclude $IC(\overline{\co}_{\lambda}, \cl)$ is not a summand of $\psi_*\shc_{\tilde{\cm}}[\dim \cn]$.
\end{proof}


\appendix

\section{Proof of Lemma \ref{lem.ZintypeD}} \label{appendix.ZintypeD}
Using the notational conventions of~\cite{HumphreysLieAlg}, we have
\[
\lambda_i = \alpha_1 + \cdots + (i-1)\alpha_{i-1}+ i( \alpha_i+ \cdots + \alpha_{n-2}) + \frac{i}{2}(\alpha_{n-1}+\alpha_n)
\]
for all $1\leq i \leq n-2$, and
\begin{eqnarray*}
\lambda_{n-1} &=&  \frac{1}{2} \left( \alpha_1 + 2\alpha_2 + \cdots + (n-2)\alpha_{n-2} + \frac{n}{2}\alpha_{n-1}+\frac{n-2}{2}\alpha_n \right) \\
\lambda_{n} &=&   \frac{1}{2} \left( \alpha_1 + 2\alpha_2 + \cdots + (n-2)\alpha_{n-2} + \frac{n-2}{2}\alpha_{n-1}+\frac{n}{2}\alpha_n \right).   \\
\end{eqnarray*}
In particular, we obtain the following.

\begin{Lem}\label{lem.Dformula} Writing each simple root of the type $D_n$ root system in terms of the fundamental weights we obtain
\begin{eqnarray*}
\alpha_1 &=& 2\lambda_1-\lambda_2\\
\alpha_i &=& -\lambda_{i-1}+2\lambda_i -\lambda_{i+1},\; 2\leq i \leq n-3\\
\alpha_{n-2} &=& -\lambda_{n-3}+2\lambda_{n-2}-\lambda_{n-1}-\lambda_n\\
\alpha_{n-1} &=& -\lambda_{n-2}+2\lambda_{n-1}\\
\alpha_n &=& -\lambda_{n-2}+2\lambda_{n}.
\end{eqnarray*}
\end{Lem}
The proof is straightforward and we omit it.

We now prove Lemma \ref{lem.ZintypeD}. First, by Lemma~\ref{lem.Dformula} any $(y_1, \ldots, y_n)$ in the kernel of the map $T \to T_{ad}$ must satisfy
\begin{eqnarray}
\label{eqn.D1}y_1^2 &=&y_2,\\ 
\label{eqn.D2}y_i^2&=&y_{i-1}y_{i+1}\, (2\leq i \leq n-3),\\ 
\label{eqn.D3}y_{n-2}^2&=&y_{n-3}y_{n-1}y_n,\\ 
\label{eqn.D4}y_{n-1}^2 &=& y_{n-2},\\ 
\label{eqn.D5}y_n^2&=& y_{n-2}.
\end{eqnarray}
We claim that any such kernel element must also satisfy the following:
\begin{enumerate}
\item $y_{n-k}=y_{n-1}^2 = y_{n}^2$ for all $k$ even such that $k\geq 2$,
\item $y_{n-k}= y_{n-1}y_{n}$ for all $k$ odd such that $k\geq 3$, 
\item $y_i^2=y_{i+1}^2$ for all $1\leq i \leq n-2$, and
\item $y_2=1$.
\end{enumerate}
We argue that (1) and (2) hold by strong induction on $k$. First, equations~\eqref{eqn.D4} and~\eqref{eqn.D5} imply $y_n^2=y_{n-1}^2$ and that Claim (1) holds for $k=2$.  Using~\eqref{eqn.D3} and Claim (1) for $k=2$ we obtain 
\[
y_{n-1}^2y_{n}^2 = y_{n-2}^2 = y_{n-3}y_{n-1}y_n \Rightarrow y_{n-3} = y_{n-1}y_n
\]
so Claim (2) holds for $k=3$. 

Now suppose that both Claim (1) and (2) hold for all $2 \leq k <j$. Applying~\eqref{eqn.D2} with $i=n-(j-1)$ we have that $y_{n-j}y_{n-(j-2)} = y_{n-(j-1)}^2 $. Now, if $j$ is even, then
\[
y_{n-j}y_{n-(j-2)} = y_{n-(j-1)}^2 \Rightarrow y_{n-j}y_n^2 = y_{n-1}^2y_n^2 \Rightarrow y_{n-j} = y_{n-1}^2
\]
where the second implication uses the induction hypothesis for (1) with $k=j-2$ and the induction hypothesis for (2) with $k=j-1$. When $j$ is odd, then 
\[
y_{n-j}y_{n-(j-2)} = y_{n-(j-1)}^2 \Rightarrow  y_{n-j} y_{n-1}y_n = y_{n-1}^2y_n^2 \Rightarrow y_{n-j} = y_{n-1}y_n
\]
where the second implication uses the induction hypothesis for (2) with $k=j-2$ and the induction hypothesis for (1) with $k=j-1$. 
This completes the proof of Claims (1) and (2).  Claim (3) follows immediately since for all $1\leq i\leq n-2$ we have $y_i^2 = y_{n-1}^2y_n^2$ by (1) and (2). Finally, applying (3) with $i=1$ and equation~\eqref{eqn.D1} yields
\[
y_2=y_1^2=y_2^2 \Rightarrow y_2=1
\]
which proves Claim (4).

If $n$ is even, then Claims (1) and (4) imply $y_k=1$ for $k$ even with $2\leq k \leq n-2$.  Equations~\eqref{eqn.D4} and~\eqref{eqn.D5} become
\[
y_{n-1}^2 = y_n^2 = 1\Rightarrow y_n, y_{n-1}\in \{1, -1\}.
\]
If $y_{n-1}$ and $y_n$ have opposite signs, then $y_k=-1$ for $k$ odd with $1\leq k \leq n-3$ by Claim~(2). We now see that $w_1$ and $w_2$ as in the statement of Lemma~\ref{lem.ZintypeD}~(1) are generators for $Z$, as desired. 

If $n$ is odd, then Claims (2) and (4) imply $1 = y_k= y_{n-1}y_n$ for $k$ even with $2\leq k \leq n-3$. In particular, $y_n$ and $y_{n-1}$ are multiplicative inverses of one another.  Furthermore, $y_{n-1}^4 = y_{n-1}^2 y_n^2 = y_2^2 = 1$ so $y_{n-1}$ is a fourth root of unity.  Finally, if $y_{n-1} = \omega$ then $y_n=\omega^3$ and Claim (1) implies $y_k = y_{n-1}^2 = -1$ for $k$ odd with $1\leq k \leq n-2$.  We now see that $z$ as in the statement of Lemma~\ref{lem.ZintypeD}~(2) is a generator for $Z$, as desired.


\section{Proof of Lemma \ref{lem:Bonnafe}} \label{appendix.Bonnafe}
In this section we prove Lemma \ref{lem:Bonnafe}, which relates the generators of the center $Z$ of $G$ defined in \cite{Bonnafe2004} to our description of $Z$ from Section~\ref{ss:Dn}.   
We begin by recalling some results and notation from \cite{Bonnafe2004}; although Bonnaf\'e works more generally, we assume that
$G$ is semisimple and simply connected, and the ground field is $\C$.  Then we specialize to the case of the lemma, where $G$ is of type
$D_n$ for $n$ even.  In this section we denote the character of $T$ corresponding to a root or weight $\nu$ simply by $\nu$, rather
than by $e^{\nu}$ as in previous sections.

Following \cite[\S 1.A]{Bonnafe2004}, taking ${\bf H} = T$, and ${\mathbb F} = \C$, let $X(T) = \cp$ denote the group of characters of $T$, and $Y(T)$ the group of cocharacters, i.e., homomorphisms $\C^* \to T$. 
There is a perfect pairing $\langle \cdot, \cdot \rangle: X(T) \times Y(T) \to \Z$, which for $x \in X(T)$, $y \in Y(T)$,
satisfies
$$
x( y(z) ) = z^{\langle x, y \rangle}.
$$
This extends to a pairing $X(T)_{\QQ} \times Y(T)_{\QQ} \to \QQ$, where the subscript $\QQ$ means $\otimes_{\Z} \QQ$.

Let $\mu_{\infty}$ denote the roots of unity in $\C^*$.  We define
$$
\iota: (\QQ/ \Z)_{tors} = \QQ/ \Z \to \mu_{\infty}, \ \ \ \iota(a/b) = e^{2 \pi i a/b} = (\zeta_b)^a,
$$
where $\zeta_b =  e^{2 \pi i /b}$. Define
\begin{align*}
\ibarT: Y(T)_\QQ \to T ; \; y \otimes_{\Z} r \mapsto y(\iota (r))
\end{align*}
for all $y\in Y(T)$ and $r\in \QQ$.

\begin{example} \label{ex.B1}
Let $T = \C^*$.  Suppose $X(T) = \Z \lambda$ and $Y(T) = \Z \varphi$ with $\langle \lambda, \varphi \rangle = 1$.
Concretely, $\lambda(t) = t$ and $\varphi(z) = z$ for $t \in T$ and $z \in \C^*$ (of course $T = \C^*$ but we use different letters to distinguish the source and target). 
To illustrate the definitions, we calculate $\ibarT\left(\frac{1}{3} \varphi\right)$.  By definition, $\frac{1}{3} \varphi = \varphi \otimes_{\Z} \frac{1}{3}$.
We have $\iota(1/3) = \zeta_3 = e^{2 \pi i/3}$, so
\begin{align*}
\ibarT\left(\frac{1}{3} \varphi\right) & =  \ibarT\left(\varphi \otimes_{\Z} \frac{1}{3}\right) = \varphi(\iota (1/3)) = \varphi(\zeta_3) = \zeta_3.
\end{align*}
\end{example}

In general, if $\lambda \in Y(T)_{\QQ}$, the value of $\ibarT(\lambda)$ depends on the integer $b$ such that $b \lambda \in Y(T)$, in other words, on what $b$ is needed to clear denominators.  
If $\lambda$ is already in $Y(T)$ then $\ibarT(\lambda) = 1$ (i.e., the identity in $T$).

Let $\gD$ denote the set of simple roots.
The fundamental weights $\varpi_{\ga}$ for $\ga \in \gD$ are a basis for $X(T)$; the corresponding dual basis of $Y(T)$ is the basis of simple coroots.  (In the rest of this paper, the simple roots have generally been denoted $\ga_1, \ldots, \ga_n$, and $\varpi_{\ga_i}$ has
been denoted $\lambda_i$. We change notation here to align with~\cite{Bonnafe2004}.)
The simple roots $\gD$ are a basis for $X(T)_{\QQ}$ (although not in general for $X(T)$); 
let $\{ \varpi^{\vee}_{\ga} \}$ denote the basis of $Y(T)_{\QQ}$ dual to the basis $\gD$ of $X(T)_{\QQ}$. 
The center $Z=\cz(G)$ is generated by the $\ibarT(\varpi^{\vee}_{\ga} )$ (see the remarks before Proposition~2.8 of \cite{Bonnafe2004}).

We note that the elements  $\varpi^{\vee}_{\ga}$ need not all be in $Y(T)$, since although the simple roots are in $X(T)$, they are only a basis over $\QQ$. 

\begin{example}
Let $G = \SL(2)$ with $T = \{ \diag(t, t^{-1}) \} \cong \C^*$.  We abuse notation and write $t$ for $\diag(t, t^{-1})$.  Then $X(T) = \Z \lambda$ and the coroot $\ga^{\vee}$ is $\varphi$, where $\lambda$ and  $\varphi$ are as defined in Example~\ref{ex.B1} above.
The simple root is $\ga = 2 \lambda$ and its dual element is $ \varpi^{\vee}_{\ga} = \frac{1}{2} \varphi$.  We have
\begin{align*}
\ibarT(\varpi^{\vee}_{\ga} )  & = \ibarT\left(\frac{1}{2} \varphi \right) =  \ibarT\left(\varphi \otimes_{\Z} \frac{1}{2}\right) = \varphi(\iota(1/2)) = \varphi(\zeta_2) = \zeta_2 = -1.
\end{align*}
Thus, the center $Z$ is generated by $-1 = \diag(-1, -1) \in T$.
\end{example}

\begin{example}
Let $G = \SL(3)$ with $T = \{ \diag(t_1, t_2, t_3) \mid t_1 t_2 t_3 = 1 \}$.  The simple roots are given by
$\ga_1 (t) = t_1/t_2$, $\ga_2(t) = t_2/t_3$.  The fundamental dominant weights are $\varpi_{\ga_1}(t) = t_1$, $\varpi_{\ga_2} = t_1 t_2$.  
The fundamental coroots are $\ga^{\vee}_1(z) = \diag(z, z^{-1}, 1)$ and $\ga^{\vee}_2(z) = \diag(1, z, z^{-1})$.  As expected, $\langle \ga_i, \ga_i^{\vee} \rangle = 2$ for $i=1,2$.
This can be seen directly: for example, for $\ga_1$ we have
$$
\ga_1(\ga^{\vee}_1(z)) = \ga_1( \diag(z, z^{-1}, 1)) = z^2 = z^{\langle \ga_1, \ga^{\vee}_1 \rangle}.
$$
We know that $ \varpi_{\ga_1}  = \frac{2}{3} \ga_1 + \frac{1}{3} \ga_2 $ and $\varpi_{\ga_2}  = \frac{1}{3} \ga_1 + \frac{2}{3} \ga_2$.  
The analogous equations hold for cocharacters, i.e.,  $\varpi^{\vee}_{\ga_1}  = \frac{2}{3} \ga^{\vee}_1 + \frac{1}{3} \ga^{\vee}_2$
and $\varpi^{\vee}_{\ga_2}  = \frac{1}{3} \ga^{\vee}_1 + \frac{2}{3} \ga^{\vee}_2$.
We can verify directly that
\begin{equation} \label{e:Z}
Z = \{1,  \ibarT(\varpi^{\vee}_{\ga_1} ) , \ibarT(\varpi^{\vee}_{\ga_2} ) \}.
\end{equation}
Indeed, it suffices to show that, writing $\zeta = \zeta_3 = e^{2 \pi i/3}$ and $1=(1,1,1)$ for the identity in $T$, we have
\begin{equation} \label{e:Z2}
\ibarT(\varpi^{\vee}_{\ga_1} ) = \zeta^2 \cdot 1 \mbox{   and   } \ibarT(\varpi^{\vee}_{\ga_2} )  = \zeta \cdot 1.
\end{equation}
To see this, observe that
\begin{align*}
(3\varpi^{\vee}_{\ga_1} )(z) & =  (2 \ga^{\vee}_1 + \ga^{\vee}_2) (z) = \ga^{\vee}_1(z^2) \ga^{\vee}_2 (z) \\
& = \diag(z^2, z^{-2}, 1) \cdot \diag(1, z, z^{-1}) = \diag(z^2, z^{-1}, z^{-1}).
\end{align*}
Hence 
\begin{align*}
\ibarT(\varpi^{\vee}_{\ga_1} )  & = \ibarT\left((3\varpi^{\vee}_{\ga_1} ) \otimes_{\Z} \frac{1}{3} \right) = (3\varpi^{\vee}_{\ga_1} ) (\iota(1/3))  = (3\varpi^{\vee}_{\ga_1} )(\zeta) \\ 
& = \diag(\zeta^2, \zeta^{-1}, \zeta^{-1}) = \diag(\zeta^2, \zeta^2, \zeta^2)  = \zeta^2 \cdot 1.
\end{align*}
This proves the first equation of \eqref{e:Z2}; the second is proved similarly.
\end{example}

For the remainder of this section, let $G$ be of type $D_{2r}$.  We use notation from Section~\ref{ss:Dn}.  
The minuscule weights are $\varpi_1:= \varpi_{\ga_1}, \varpi_{2r-1}:=\varpi_{\ga_{2r-1}}, \varpi_{2r}:=\varpi_{\ga_{2r}}$ (these are denoted respectively in Section~\ref{ss:Dn} by $\lambda_1$, $\lambda_{2r-1}$, and $\lambda_{2r}$).
Let $w_1, w_2$ be the elements of $Z$ defined in Section~\ref{ss:Dn}.
To prove Lemma \ref{lem:Bonnafe}, we must prove that if $r$ is odd (i.e., if $2r \equiv 2 \textup{ mod } 4$), then
\begin{equation} \label{e:Bonnafe-odd}
\ibarT(\varpi_1^{\vee})  = w_1 w_2 , \hspace{.2in} \ibarT(\varpi_{2r-1}^{\vee})  = w_2 , \hspace{.2in}  \ibarT(\varpi_{2r}^{\vee})  = w_1,
\end{equation}
while if $r$ is even (i.e., if $2r \equiv 0 \textup{ mod } 4$), then
\begin{equation} \label{e:Bonnafe-even}
\ibarT(\varpi_1^{\vee})   = w_1 w_2 , \hspace{.2in}   \ibarT(\varpi_{2r-1}^{\vee})  = w_1 , \hspace{.2in}  \ibarT(\varpi_{2r}^{\vee})  = w_2.
\end{equation}

As a preliminary, note that in terms of the $y_i$ coordinates as in Lemma~\ref{lem.ZintypeD}, we have
\begin{equation} \label{e:yi}
y_i(\ga_j^{\vee}(z)) = \varpi_{i}(\ga_j^{\vee}(z)) = z^{\langle \varpi_i, \ga_j^{\vee} \rangle} = \begin{cases}
z & \mbox{if } i = j \\
1 & \mbox{if } i \neq j.
\end{cases}
\end{equation}
We will determine $\ibarT(\varpi_1^{\vee}), \ibarT(\varpi_{2r-1}^{\vee}), \ibarT(\varpi_{2r}^{\vee})$.
Using the information in Table~\ref{fig.minimal.table} we have
\begin{align*}
\varpi_1^{\vee} &= \ga_1^{\vee} + \ga_2^{\vee} + \cdots+  \ga_{2r-2}^{\vee}   + \frac{1}{2} (\ga_{2r-1}^{\vee} + \ga_{2r}^{\vee}) \\
\varpi_{2r-1}^{\vee} & = \frac{1}{2}( \ga_1^{\vee} +  2 \ga_2^{\vee} + \cdots +  (2r-2)  \ga_{2r-2}^{\vee}) + \frac{r}{2} \ga_{2r-1}^{\vee} + \frac{r-1}{2}\ga_{2r}^{\vee} \\
\varpi_{2r}^{\vee} & = \frac{1}{2}( \ga_1^{\vee} +  2 \ga_2^{\vee} +  \cdots +  (2r-2)  \ga_{2r-2}^{\vee}) +  \frac{r-1}{2} \ga_{2r-1}^{\vee} + \frac{r}{2}  \ga_{2r}^{\vee}.
\end{align*}
We see that
\begin{align*}
\ibarT(\varpi_1^{\vee}) & = \ibarT\left(2 \varpi_1^{\vee} \otimes \frac{1}{2}\right) = (2 \varpi_1^{\vee} ) (\iota(1/2) ) = (2 \varpi_1^{\vee} ) (-1) \\
& =  \left( \prod_{i = 1}^{2r-2} \ga_i^\vee((-1)^2) \right) \ga_{2r-1}^{\vee}(-1) \ga_{2r}^{\vee}(-1)  =  \ga_{2r-1}^{\vee}(-1) \ga_{2r}^{\vee}(-1) .
\end{align*}
Using \eqref{e:yi}, we deduce that in terms of the coordinates $(y_1, \ldots, y_{2r})$, 
$$
\ibarT(\varpi_1^{\vee}) = (1, 1, \ldots, 1, -1, -1) = w_1 w_2,
$$
where the last equality follows from the formulas for $w_1$ and $w_2$ in Lemma \ref{lem.ZintypeD}. 

Now suppose $r$ is odd. Then we have 
\begin{align*}
\ibarT(\varpi_{2r-1}^\vee) &= \ibarT\left( 2 \varpi_{2r-1}^\vee   \otimes \frac{1}{2} \right) = (2 \varpi_{2r-1}^{\vee} ) (\iota(1/2) ) = (2 \varpi_{2r-1}^{\vee} ) (-1) \\
&=\ga_1^{\vee}(-1) \ga_3^{\vee}(-1) \cdots \ga_{2r-3}^{\vee}(-1)  \ga_{2r-1}^{\vee}(-1).
\end{align*}
A similar calculation shows that $\ibarT(\varpi_{2r}^{\vee})  =  \ga_1^{\vee}(-1) \ga_3^{\vee}(-1) \cdots \ga_{2r-3}^{\vee}(-1)  \ga_{2r}^{\vee}(-1)$.
Again using  \eqref{e:yi} and Lemma \ref{lem.ZintypeD}, we see that in 
terms of the coordinates $(y_1, \ldots, y_{2r})$,
\begin{align*}
\ibarT(\varpi_{2r-1}^{\vee}) = (-1,1 -1,\ldots, -1,1,-1,1) = w_2 \\
\ibarT(\varpi_{2r}^{\vee})  =  (-1,1,-1, \ldots, -1,1,1,-1) = w_1.
\end{align*}
This proves \eqref{e:Bonnafe-odd}, and so proves Lemma \ref{lem:Bonnafe} in the case where $r$ is odd.  The case where $r$ is even
is almost the same and is left to the reader.

\section{An overview of the generalized Springer correspondence} \label{s.overview}
In this appendix we give an overview of Lusztig's generalized Springer correspondence in classical types (see~\cite{Lu84}).  
In particular, we describe Lusztig's classification of cuspidal data for simple algebraic groups.  By doing so, we hope to illustrate
how our results (which relate to minimal cuspidal data) fit into the larger picture, and to provide some details
relevant to the discussion in Section \ref{ss.strategy} of our strategy of proof.

Let $G$ be a simple algebraic group and $L\subseteq G$ a Levi subgroup.  A triple $\vc = (L, \co_L, \cl^L)$ for $G$ is 
a \emph{cuspidal datum} 
if the local system $\cl$ is irreducible, and the simple perverse sheaf $IC(\overline{\co_L}, \cl_L)$ on $\cn_L$ is \emph{cuspidal} for $L$.  This means that for every proper parabolic subgroup $P_L\subset L$ with Levi factor $M$ we have $\mathrm{res}_{M}^L IC(\overline{\co_L}, \cl_L)=0$.  There are only finitely many standard Levi subgroups $L$ of $G$ which admit a cuspidal perverse sheaf, or equivalently, which occur as the Levi subgroup in a cuspidal datum for $G$. We are particularly interested in the case where $\co_L$ is the regular orbit, so we record the following (known) lemma.

\begin{Lem} \label{l:typeAcuspidal}
Suppose $L$ is a reductive algebraic group.  If $L$ is of type $A$, i.e., if the Lie algebra $[\fl, \fl]$ is a product of type $A$ factors, then the 
only nilpotent orbit $\co_L$ supporting a cuspidal perverse sheaf $IC(\overline{\co_L}, \cl_L)$ is the regular orbit $\co_L^{\reg}$.   Conversely, if $\co_L$ is the regular 
orbit, and $IC(\overline{\co_L}, \cl_L)$ is a cuspidal perverse
sheaf, then $L$ is of type
$A$.
\end{Lem}

\begin{proof}
The only nilpotent orbits supporting a cuspidal perverse sheaf are the distinguished nilpotent orbits, by \cite[Prop.~2.6]{AHJR17a}.  In type $A$, the only distinguished nilpotent orbit is the regular nilpotent orbit.  This implies the first assertion of the lemma.  The converse follows by inspection of Lusztig's classification of cuspidal perverse sheaves (see \cite{Lu84}).
\end{proof}

Lusztig classifies the cuspidal data that occur in each simple type in \cite{Lu84}, with some details in the case of the Spin groups addressed by Lusztig and Spaltenstein in \cite{LuSpal85}.   Key to this classification is the connection between central characters and cuspidal datum, so
we begin by discussing the connection between characters of the center $Z$ of $G$ and local systems on nilpotent orbits.  Both nilpotent $G$-orbits and nilpotent $L$-orbits of a Levi subgroup $L$ play a role in what follows.

Suppose $\co$ is a nilpotent $G$-orbit.  There is a natural map $Z = \cz(G) \to \pi_1^G(\co) = G^x/G^x_0$ where
$x$ is any element of $\co$. Let $(\co, \cl)$ be a pair consisting of a nilpotent orbit $\co$ and an irreducible $G$-equivariant local system $\cl$ on $\co$.
The local system corresponds to an irreducible representation of $\pi_1^G(\co) $ on a vector space $V$.
Associated to this pair is a central character $\chi = \chi_{(\co, \cl)}$ defined as the composition 
$$
\chi: Z \to \pi_1^G(\co) \to GL(V),
$$
where the second map is the representation of $\pi_1^G(\co)$ corresponding to $\cl$.  The map $\chi$ takes $Z$ to the group of
scalar matrices in $GL(V)$, since the representation $V$ is irreducible, and the action of $Z$ commutes with the action of
$\pi_1^G(\co)$.  Thus, this composition can be viewed as a map $\chi: Z \to \C^*$.

A character of $Z$ is also associated to an irreducible local system $\cl_L$ on a nilpotent orbit $\co_L$ of a Levi subgroup $L$ of $G$.
Indeed, by the above discussion, $\cl_L$ yields a character of $\cz(L)$, which pulls back to $Z$ via the natural map $Z\rightarrow \cz(L)$.  
If $\vc = (L, \co_L, \cl^L)$ is a cuspidal datum, we refer to this character of $Z$ as the central character of~$\vc$.

\begin{Thm}[Lusztig]\label{lu_thm_central}
Suppose $\mathbf{c}$ is a cuspidal datum of a simple group $G$.  Then
the central character of $\vc$ coincides
with the central character of any $\cl$ such that $IC(\overline{\co}, \cl)$
appears as a summand in $\mathbb A_\textbf{c}$.
\end{Thm}

The following result summarizes some results of Lusztig concerning the relationship between central
characters and cuspidal data.

\begin{Thm}[Lusztig]\label{lu_summary}
Let $\chi$ be a character of $Z$.
\begin{enumerate}
\item Any Levi subgroup supports at most one cuspidal datum with central character $\chi$.  
\item There is exactly one cuspidal datum with 
central character $\chi$ with Levi subgroup
of type $A$.
\item The minimal cuspidal datum is the unique cuspidal
datum with central character $\chi$ for which the nilpotent orbit $\co_L$ is regular.
This is the cuspidal datum with Levi subgroup of type $A$.
\end{enumerate}
\end{Thm}

Parts (1) and (2) of this result follow from Lusztig's classification of cuspidal data, which we discuss
in more detail below in classical types.  Part (3) follows from the first two parts and Lemma \ref{l:typeAcuspidal}.
As noted in Section \ref{ss.strategy}, Theorem \ref{lu_summary} implies that a cuspidal
datum is minimal if and only if its Levi subgroup $L$ is of type $A$.

There 
can be multiple Levi subgroups which support
a cuspidal datum for $G$ with a given central character $\chi$, but Theorem \ref{lu_summary}
implies that to specify a cuspidal datum $\vc$ corresponding to a fixed $\chi$,
it suffices to specify the Levi subgroup $L$ corresponding to $\vc$, since then the pair $(\co_L, \cl_L)$
is determined.  For the minimal cuspidal datum, the fact that the local system on 
$\co_L^{reg}$ is uniquely
determined by the requirement that its central character is $\chi$ can be seen directly.  The reason is
that for this $L$, the map $Z \to \cz(L) = \pi_1^L(\co_L^{reg})$ is surjective,
so the map $\hat{\cz(L)} \to \hat{Z}$ is injective.

We now discuss the cuspidal data in each of the classical types,
following \cite{Lu84, LuSpal85}.  In each type, we consider the central characters individually, and for
each central character, describe the Levi subgroups which occur as part of a cuspidal datum $\vc =  (L, \co_L, \cl_L)$ with this character.
A description of the orbit $\co_L$ of a cuspidal datum can be found in \cite{Lu84}, \cite{LuSpal85}.
This is not needed for our discussion, so we do not recall this here.
For the trivial central character, the minimal cuspidal datum is always $\vc = (T, \{0 \}, \cl_{triv})$; the corresponding Lusztig sheaf is the Springer sheaf.

In the discussion below, we refer to a Levi subgroup by the corresponding subdiagram of the Dynkin diagram of $G$.  Note
that for small values of $k$, the groups $B_k$, $C_k$, and $D_k$ coincide with groups of type $A$.

In type $A_{n-1}$, there is a unique cuspidal datum associated to any central character $\chi$ of $Z$. Hence the only cuspidal data
are the minimal cuspidal data, which are described as follows.
The center $Z$ of  $G$ is isomorphic to $\mathbb Z_n$.   If $\chi$ is of order $n$, then the Lusztig sheaf $A_{\vc_{\chi}}$  is a single simple perverse sheaf. More precisely, $A_{\vc_{\chi}}$ is the simple perverse sheaf $IC(\co^{\reg},\cl_\chi)$ where $\cl_\chi$ is the  local system on the principal orbit  with central character $\chi$.  If $\chi_d$ is a character of order $d$, where $d<n$ is a divisor of $n$, then
$\vc_{\chi_d} = (L_{d}, \co_{L_d}^{\reg}, \cl^{L_{d}})$ where $L_d = L_{\chi_d}$ is a Levi subgroup with subdiagram $A_{d-1}\times \cdots \times A_{d-1}$, with $n/d$ factors (as in the first row of Table~\ref{fig.Levi.table}). 
The map $Z \to \cz(L_d) = \pi_1^{L_d}(\co_{L_d}^{\reg})$ is onto, and the local system $\cl^L_{\chi_d}$ corresponds to the unique character of $\cz(L)$ which pulls back to $\chi_d$ . Taking $d = 1$ yields the trivial central character $\chi_1$, and $A_{\vc_{\chi_1}}$ is the Springer sheaf. 

In type $C_n$, the center $Z$ of the simply connected group $G = \Sp(2n)$ is $\Z_2$, so there are only two central characters.  If $\chi$ is the trivial character, the associated cuspidal data are the minimal cuspidal datum $\vc = (T, \{ 0 \}, \cl_{triv})$ and a cuspidal datum corresponding to each Levi subgroup with subdiagram $C_k$, where $1 \leq k \leq n$ and $k$ is an even triangular number (recall that a triangular number is a number of the form
$\frac{1}{2} m(m+1)$ for some integer $m$).  If $\chi$ is the nontrivial character, then a cuspidal datum occurs for each Levi of type $C_k$ for $1\leq k\leq n$ where $k$ is an odd triangular number.  The minimal cuspidal datum occurs with $k=1$ and has a Levi subgroup of type $C_1 = A_1$.

In type $B_n$, the center $Z$ of the simply connected group $G = \Spin(2n+1)$ is $\Z_2$, so again, there are two central characters.  If $\chi$ is the trivial character, then there
is a cuspidal datum corresponding to each Levi subgroup with subdiagram $B_k$, where $0 \leq k \leq n$ and $2k+1$ is a perfect square. These appear in equation \eqref{eqn:Li} as $L_i$ where $i = \frac{1}{2}(2n+1-j^2)$, if $j$ is the defect of a symbol.  In that case, $2k+1 = j^2$.   The minimal 
cuspidal datum $\vc = (T, \{ 0 \}, \cl_{triv})$ corresponds to $k=0$.  If $\chi$ is the nontrivial character, there is a cuspidal datum for each Levi of the form $A_1\times \cdots\times A_1\times B_k$ where $0 \leq k\leq n$, $2k+1 = \frac{1}{2}m(m+1)$ for some integer $m$, and $k\equiv n \pmod{2}$.  These appear as $L_d$ in equation \eqref{eqnBLd} and
Theorem \ref{thm.Lusztig-Spaltenstein}, where $2d-1 = m$, or equivalently, $k = \frac{1}{2}(d-1)(2d+1)$.   In this case, the number of $A_1$ components is taken as large as possible.  
The subdiagram is a product of type $A$ factors only when $k = 0$ or $k = 1$, so these give the Levi subgroups for the minimal cuspidal data.

In type $D_n$, the center of the simply connected group $G = \Spin(2n)$ is $Z \cong \mathbb Z_2\times \mathbb Z_2$ if $n$ is even and $Z \cong \mathbb Z_4$ if $n$ is odd. 
For the trivial character of $Z$, the Levi subgroups which appear in cuspidal data are:
\begin{enumerate}
\item The torus $T$ (which yields the minimal cuspidal datum); and
\item Any Levi subgroup with diagram $D_k$ where $2k$ is a square, $k\leq n$, and $k\equiv 0\pmod 4$.  
\end{enumerate}
These Levi subgroups appear in equation \eqref{eqn:Li} as $L_i$ where $i = \frac{1}{2}(2n-j^2)$, if $j$ is the defect of a symbol (so $2k=j^2$ is a square).  The torus $T$ corresponds to $j=0$.

There are three nontrivial characters of $Z$, denoted by $\chi_1$, $\chi_{n-1}$, and $\chi_n$ (see Section \ref{ss:Dn}). Of these, only $\chi_1$ is trivial on the kernel of the map from $\Spin(2n)$ to $\SO(2n)$.
The cuspidal data with central character $\chi_1$ correspond to the following Levi subgroups: 
\begin{enumerate}
\item A subgroup of type $A_1\times A_1 \cong D_2$ (which yields the minimal cuspidal datum); and
\item Any Levi subgroup with diagram $D_k$ where $2k$ is a square, $k\leq n$, and $k\equiv 2\pmod 4$.  
\end{enumerate}
These Levi subgroups appear in equation \eqref{eqn:Li} as $L_i$ where $i = \frac{1}{2}(2n-j^2)$, if $j$ is the defect of a symbol and $2k=j^2$.  Note that (1) occurs when $j = 2$.

The two remaining characters $\chi_{n-1}$ and $\chi_n$ have cuspidal data in which the same Levi subgroups occur for each character: 
\begin{enumerate}
\item The entire group with diagram $D_n$, if $2n$ is a triangular number; and
\item The Levi subgroups with diagram $A_1\times\dots \times A_1\times D_k$ where $k=\frac{1}{2}m(2m-1)$ for some $m\equiv 2n\pmod 4$, $k<n$; the number of $A_1$ components is $\tfrac{1}{2}(n-k)$. 
\end{enumerate}
(Note that (1) is technically superfluous, since it is the case of (2) where $k = n$, but we separate it for emphasis.)
These Levi subgroups appear as $L_d$ in Equation \eqref{eqnDLd} and Theorem \ref{thm.Lusztig-Spaltenstein}, where $d=m$, or equivalently, $k = \frac{1}{2}d(2d-1)$.
The minimal cuspidal data for these characters are as follows.   If $n$ is even, then $k = 0$ and the Levi subgroup is a product of $\tfrac{n}{2}$ components of type $A_1$.
If $n$ is odd, then $k = 3$ and the Levi subgroup is a product of $\tfrac{n-3}{2}$ components of type $A_1$ and a component of type $D_3=A_3$.
Note that since $\chi_1$ and the trivial character are the characters of $Z$ which are trivial on the kernel of 
$\Spin(2n) \to \SO(2n)$, the $IC(\overline{\co}_L, \cl^L)$ which are equivariant for the group $\SO(2n)$ appear in the Lusztig sheaves with these central characters.

The following examples illustrate the generalized Springer correspondence for the spin groups. These examples were computed with the assistance of the CHEVIE package in Gap3 \cite{GH96}, but can be confirmed by hand using the results from \cite{Lu84} and \cite{LuSpal85}.

\begin{example}\label{ex.Spin.odd}
We describe the cuspidal data and Lusztig sheaves in two examples for $\Spin(N)$ with $N$ odd.
First suppose $N = 7$ (type $B_3$).  The center of $\SO(7)$ is trivial, so the only central character for $\SO(7)$ is
trivial.  For $N=7$, the only Levi subgroup supporting a cuspidal datum corresponding to the trivial character is the torus $T$ and the associated Lusztig sheaf is the Springer sheaf.   For $\Spin(7)$, which has center $\Z_2$, the only cuspidal datum corresponding to the nontrivial central character $\chi$ is the minimal cuspidal datum $\mathbf{c}_\chi$. The corresponding Levi subgroup $L_\chi$ is of type $A_1 \times A_1$. The Lusztig sheaf for this cuspidal datum contains simple perverse sheaves supported on the closures of the orbits with partition classifications $[7]$ and $[3\ 2^2]$ only.  Note that 
these orbits are rather odd, which is consistent with Theorem \ref{thm.Lusztig-Spaltenstein}.

For a larger example where $N$ is odd, suppose $N = 13$ (type $B_6$).   
The center of $\SO(13)$ is trivial.  There are two cuspidal data with the trivial central character: the cuspidal datum where $L = T$, which yields the Springer sheaf, and a cuspidal datum where $L$ corresponds to the Dynkin subdiagram of $B_4$.  The nilpotent orbit in $L$ for the latter datum is the orbit corresponding to the partition $[5 \ 3 \ 1]$.  Thus, there are two Lusztig sheaves corresponding to the trivial central character, and only the first is a summand of $\psi_*\shc_{\tilde{\cm}}[\dim \cn]$.   For $G = \Spin(13)$,
the center is $\Z_2$.  Corresponding to the nontrivial central character $\chi$ there is only one cuspidal datum, namely the minimal datum $\mathbf{c}_\chi$ with Levi subgroup $L_\chi$ of type $A_1 \times A_1 \times A_1$.
The simple perverse sheaves in the corresponding Lusztig sheaf $\mathbb{A}_{\mathbf{c}_\chi}$ are each of the form  $IC(\bar{\co}, \cl)$, where $\co$ is a rather odd orbit, and $\cl$ is the unique local system on $\co$ with central character $\chi$.
\end{example}

\begin{example} \label{ex.Spin.even} We describe the cuspidal data and Lusztig sheaves for $\Spin(8)$ (type $D_4$). We use the notation established in Section~\ref{ss:Dn} for characters of $Z=\cz(G)$.   The center of $\SO(8)$ is  $\Z_2$.  The central characters for $SO(8)$ are the characters of $Z$ which are trivial on $Z_0$, namely, the trivial
central character and $ \chi_1$ (see Lemma \ref{lem.kernel}).   
There is a unique cuspidal datum corresponding to the trivial character; the Levi subgroup is the torus $T$ and the associated Lusztig sheaf is the Springer sheaf.
There is also a unique cuspidal datum for the nontrivial character $\chi_1$ of $Z$, namely the minimal cuspidal datum $\mathbf{c}_{\chi_1}$ which has Levi subgroup with diagram as pictured below.
\begin{center}
\begin{tikzpicture}[scale=1]
\node[circle, draw, inner sep=2pt] (1) at (5,0) {};
\node[circle, draw, inner sep=2pt] (2) at (6,0) {};
\node[circle, fill=black, draw, inner sep=2pt] (3) at (7,0.3) {};
\node[circle, fill=black, draw, inner sep=2pt] (4) at (7,-0.3) {};
\draw (1)--(2);
\draw (2)--(3);
\draw (2)--(4);
\end{tikzpicture}\\
\end{center}
For $G =\Spin(8)$, the center is $\mathbb Z_2\times \mathbb Z_2$, and we have two additional nontrivial central characters, $\chi_{n-1}$ and $\chi_n$.  There is only one cuspidal datum corresponding to each of them, which is necessarily the minimal cuspidal datum. The corresponding Levi subgroups are $L_{\chi_{n-1}}$ and $L_{\chi_n}$; their
diagrams are pictured below.
\begin{center}
\begin{tikzpicture}[scale=1]
\node[circle, fill=black, draw, inner sep=2pt] (1) at (5,0) {};
\node[circle, draw, inner sep=2pt] (2) at (6,0) {};
\node[circle, fill=black, draw, inner sep=2pt] (3) at (7,0.3) {};
\node[circle, draw, inner sep=2pt] (4) at (7,-0.3) {};
\draw (1)--(2);
\draw (2)--(3);
\draw (2)--(4);
\end{tikzpicture}
\hspace*{1in}
\begin{tikzpicture}[scale=1]
\node[circle, fill=black, draw, inner sep=2pt] (1) at (5,0) {};
\node[circle, draw, inner sep=2pt] (2) at (6,0) {};
\node[circle, draw, inner sep=2pt] (3) at (7,0.3) {};
\node[circle, fill=black, draw, inner sep=2pt] (4) at (7,-0.3) {};
\draw (1)--(2);
\draw (2)--(3);
\draw (2)--(4);
\end{tikzpicture}\\
\end{center}
Again, the nilpotent orbits supporting simple perverse sheaves which appear in the Lusztig sheaves $\mathbb A_{\mathbf{c}_{\chi_{n-1}}}$ and $\mathbb A_{\mathbf{c}_{\chi_{n}}}$ are precisely the rather odd orbits,
including the very even orbits.   Each very even partition indexes two nilpotent orbits, distinguished by their weighted Dynkin diagram classifications, and each of these orbits has precisely one nontrivial local system. Thus there are two simple perverse sheaves of the form $IC(\overline{\co}_\lambda, \cl)$ when $\lambda$ is very even, one appearing as a summand of  $\mathbb A_{\mathbf{c}_{\chi_{n-1}}}$ and the other as a summand of $\mathbb A_{\mathbf{c}_{\chi_{n}}}$. Table~\ref{fig.D4} illustrates the generalized Springer correspondence for $\Spin(8)$.  Each column represents a Lusztig sheaf $\mathbb{A}_{\mathbf{c}_\chi}$ and the table entries indicate which simple perverse sheaves $IC(\overline{\co}, \cl)$ occur as a summand.
In particular, the entries reference the order of the character of $\pi_1(\co)$ corresponding to the local system $\cl$ on $\co$ defining the sheaf $IC(\overline{\co},\cl)$. So, for example, if $\pi_1(\co)\simeq  \Z_2\times \Z_2$, there are three distinct characters of order $2$ and therefore three distinct local systems.  In this example, all cuspidal data are minimal, but this is not true for larger values of $N$.  
\end{example}

\begin{table}[h]
\begin{tabular}{|c|c|c|c|c|c|}
\hline
Orbits & $\pi_1(\co)$ & $1$ & $\chi_{1}$ & $\chi_{n-1}$ & $\chi_{n}$\\
\hline
$[7\ 1]$ & $\mathbb Z_2\times \mathbb Z_2$ & triv & order 2 & order 2 & order 2\\
$[5\ 3]$ &  $\mathbb Z_2\times \mathbb Z_2$ & triv & order 2 & order 2 & order 2\\
$[5\ 1^3]$ & $ \mathbb Z_2$ & triv & order 2 & - & -\\
$[4^2]^+$ &  $\mathbb Z_2$ & triv & - & order 2 & -\\
$[4^2]^- $&  $\mathbb Z_2$ & triv & - & - & order 2\\
$[3^2\ 1^2]$ &$ \mathbb Z_2$ & triv, order 2 & - & - & -\\
$[3\ 2^2\ 1]$ & $\mathbb Z_2\times \mathbb Z_2$ & triv & order 2 & order 2 & order 2\\
$[3\ 1^5]$ & $ \mathbb Z_2$ & triv & order 2 & - & -\\
$[2^4]^+$ &  $\mathbb Z_2$ & triv & - & order 2 & -\\
$[2^4]^- $&  $\mathbb Z_2$ & triv & - & - & order 2\\
$[2^2\ 1^4]$ & triv & triv & - & - & -\\
$[1^8]$ & triv & triv & - & -& -\\
\hline
\end{tabular}
\medskip
\caption{The generalized Springer correspondence for $G=\Spin(8)$.} \label{fig.D4}
\end{table}


\ifx\undefined\bysame
\newcommand{\bysame}{\leavevmode\hbox to3em{\hrulefill}\,}
\fi

\end{document}